\documentclass[12pt,a4paper]{article}

\usepackage{amsmath,amstext,amssymb,amscd}
\usepackage{graphicx}

\newcommand{\Xcomment}[1]{}

\newtheorem{theorem}{Theorem}[section]
\newtheorem{lemma}[theorem]{Lemma}
\newtheorem{corollary}[theorem]{Corollary}

\newtheorem{prop}[theorem]{Proposition}

\makeatletter \@addtoreset{equation}{section} \makeatother

\newenvironment{proof}{\noindent{\bf Proof}\/}%
{\hfill$\qed$\medskip}

\def\qed{ \ \vrule width.1cm height.3cm depth0cm}

\newenvironment{numitem1}{\refstepcounter{equation}\begin{enumerate}%
\item[(\thesection.\arabic{equation})]}{\end{enumerate}}

\newcommand{\refeq}[1]{(\ref{eq:#1})}  

 \makeatletter
\renewcommand{\section}{\@startsection{section}{1}{0pt}%
{-3.5ex plus -1ex minus -.2ex}{2.3ex plus .2ex}%
{\normalfont\Large}}
 \makeatother

 \makeatletter
\renewcommand{\subsection}{\@startsection{subsection}{2}{0pt}%
{-3.0ex plus -1ex minus -.2ex}{-1.5ex plus .2ex}%
{\normalfont\normalsize\bf}}
 \makeatother

\newcommand{\SEC}[1]{\ref{sec:#1}}  
\newcommand{\SSEC}[1]{\ref{ssec:#1}}  

\def\Rset{{\mathbb R}}
\def\Kset{{\mathbb K}}
\def\Qset{{\mathbb Q}}
\def\Zset{{\mathbb Z}}

\def\Dscr{{\cal D}}
\def\Escr{{\cal E}}
\def\Fscr{{\cal F}}

\def\Iscr{{\cal I}}
\def\Kscr{{\cal K}}
\def\Lscr{{\cal L}}
\def\Mscr{{\cal M}}
\def\Nscr{{\cal N}}
\def\Oscr{{\cal O}}
\def\Pscr{{\cal P}}

\def\Rscr{{\cal R}}
\def\Sscr{{\cal S}}

\def\Zscr{{\cal Z}}

\def\tilde{\widetilde}

\def\bar{\overline}
\def\eps{\epsilon}

\def\Path{{\rm Path}}
\def\Succ{{\rm ISuc}}
\def\Open{{\rm Open}}
\def\Sign{{\rm sgn}}
\def\reve{{\rm rev}}
\def\Inter{{\rm Int}}

\def\bfC{{\bf C}}
\def\bfM{{\bf M}}

\def\Iw{I^\circ}
\def\Jw{J^\circ}
\def\Ib{I^\bullet}
\def\Jb{J^\bullet}
\def\Kw{{K^\circ}}
\def\Lw{{L^\circ}}
\def\Kb{{K^\bullet}}
\def\Lb{{L^\bullet}}

\def\Yr{Y^{\rm r}}
\def\Yc{Y^{\rm c}}
\def\Xr{X^{\rm r}}
\def\Xc{X^{\rm c}}
\def\Mr{M^{\rm r}}
\def\Mc{M^{\rm c}}
\def\Mrc{M^{\rm rc}}

\def\circal{\circlearrowleft}

\def\Pstan{\Pscr^{\rm st}}

\def\precast{\prec^\ast}

\def\horvert{\;\hbox{\unitlength=1mm\begin{picture}(2,3)%
\put(0,3){\line(1,0){2}}\put(2,0){\line(0,1){3}}%
\end{picture}}\;}

\def\verthor{\;\hbox{\unitlength=1mm\begin{picture}(2,3)%
\put(0,0){\line(0,1){3}}\put(0,0){\line(1,0){2}}%
\end{picture}}\;}

  \Xcomment{

  }

\begin{document}

\baselineskip=15pt
\parskip=2pt

\title{On universal quadratic identities for minors of quantum matrices}

\author{
Vladimir~I.~Danilov\thanks{Central Institute of Economics and Mathematics of
the RAS, 47, Nakhimovskii Prospect, 117418 Moscow, Russia; emails:
danilov@cemi.rssi.ru}
  \and
Alexander~V.~Karzanov\thanks{Institute for System Analysis at FRC Computer
Science and Control of the RAS, 9, Prospect 60 Let Oktyabrya, 117312 Moscow,
Russia; email: sasha@cs.isa.ru. Corresponding author.}
  }

\date{}

\maketitle

 \begin{abstract}
We give a complete combinatorial characterization of homogeneous quadratic
relations of ``universal character'' valid for minors of quantum matrices (more
precisely, for minors in the quantized coordinate ring
$\Oscr_q(\Mscr_{m,n}(\Kset))$ of $m\times n$ matrices over a field $\Kset$,
where $q\in\Kset^\ast$). This is obtained as a consequence of a study of
quantized minors of matrices generated by paths in certain planar graphs,
called \emph{SE-graphs}, generalizing the ones associated with Cauchon
diagrams. Our efficient method of verifying universal quadratic identities for
minors of quantum matrices is illustrated with many appealing examples.
 \medskip


\emph{MSC-class}: 16T99, 05C75, 05E99
 \end{abstract}

\parskip=3pt


\section{\Large Introduction}  \label{sec:intr}

The idea of quantization has proved its importance to bridge the commutative
and noncommutative versions of certain algebraic structures and promote better
understanding various aspects of the latter versions. One popular structure
studied for the last three decades is the quantized coordinate ring
$\Rscr=\Oscr_q(\Mscr_{m,n}(\Kset))$ of $m\times n$ matrices over a field
$\Kset$, where $q$ is a nonzero element of $\Kset$; it is usually called the
\emph{algebra of $m\times n$ quantum matrices}. Here $\Rscr$ is the
$\Kset$-algebra generated by the entries (indeterminates) of an $m\times n$
matrix $X$ subject to the following (quasi)commutation relations due to
Manin~\cite{man}: for $1\le i<\ell\le m$ and $1\le j<k\le n$,
   \begin{gather}
   x_{ij}x_{ik}=qx_{ik}x_{ij},\qquad x_{ij}x_{\ell j}=qx_{\ell j}x_{ij},
                                              \label{eq:xijkl}\\
   x_{ik}x_{\ell j}=x_{ \ell j}x_{ik}\quad \mbox{and}\quad
    x_{ij}x_{\ell k}-x_{\ell k}x_{ij}=(q-q^{-1})x_{ik}x_{\ell j}.  \nonumber
    \end{gather}

This paper is devoted to quadratic identities for minors of quantum matrices
(usually called quantum minors or quantized minors or $q$-minors). For
representative cases, aspects and applications of such identities, see,
e.g.,~\cite{fio,good,KL,LR,LZ,TT} (where the list is incomplete). We present a
novel, and rather transparent, combinatorial method which enables us to
completely characterize and efficiently verify homogeneous quadratic identities
of universal character that are valid for quantum minors.

The identities of our interest can be written as
  \begin{equation} \label{eq:q_ident}
  \sum(s_i q^{\delta_i} [I_i|J_i]_q\, [I'_i|J'_i]_q \colon i=1,\ldots,N)=0,
  \end{equation}
where $\delta_i\in\Zset$, ~$s_i\in\{+1,-1\}$, and $[I|J]_q$ denotes the quantum
minor whose rows and columns are indexed by $I\subseteq[m]$ and
$J\subseteq[n]$, respectively. (Hereinafter, for a positive integer $n'$, we
write $[n']$ for $\{1,2,\ldots,n'\}$.) The homogeneity means that each of the
sets $I_i\cup I'_i,\, I_i\cap I'_i,\, J_i\cup J'_i,\, J_i\cap J'_i$ is
invariant of $i$, and the term ``universal'' means that~\refeq{q_ident} should
be valid independently of $\Kset,q$ and a $q$-matrix (a matrix whose entries
obey Manin's relations and, possibly, additional ones). Note that any cortege
$(I|J,I'|J')$ may be repeated in~\refeq{q_ident} many times.

Our approach is based on two sources. The first one is the \emph{flow-matching
method} elaborated in~\cite{DKK} to characterize quadratic identities for usual
minors (viz. for $q=1$). In that case the identities are viewed simpler
than~\refeq{q_ident}, namely, as
  \begin{equation} \label{eq:com_ident}
  \sum(s_i[I_i|J_i]\, [I'_i|J'_i] \colon i=1,\ldots,N)=0.
  \end{equation}
(In fact,~\cite{DKK} deals with natural analogs of~\refeq{com_ident} over
commutative semirings, e.g. the tropical semiring $(\Rset,+,\max)$.) In the
method of~\cite{DKK}, each cortege $S=(I|J,I'|J')$ is associated with a certain
set $\Mscr(S)$ of \emph{feasible matchings} on the set $(I\triangle I')\sqcup
(J\triangle J')$ (where $A\triangle B$ denotes the symmetric difference
$(A-B)\cup(B-A)$, and $A\sqcup B$ the disjoint union of sets $A,B$). The main
theorem in~\cite{DKK} asserts that~\refeq{com_ident} is valid (universally) if
and only if the families $\Iscr^+$ and $\Iscr^-$ of corteges $S_i$ with signs
$s_i=+$ and $s_i=-$, respectively, are \emph{balanced}, in the sense that the
total families of feasible matchings for corteges occurring in $\Iscr^+$ and in
$\Iscr^-$ are equal.

The main result of this paper gives necessary and sufficient conditions for the
quantum version (in Theorems~\ref{tm:nec_q_bal} and~\ref{tm:suff_q_bal}). It
says that~\refeq{q_ident} is valid (universally) if and only if the families of
corteges $\Iscr^+$ and $\Iscr^-$ along with the function $\delta$ are
$q$-\emph{balanced}, which now means the existence of a bijection between the
feasible matchings for $\Iscr^+$ and $\Iscr^-$ that is agreeable with $\delta$
in a certain sense. The proof of necessity (Theorem~\ref{tm:nec_q_bal})
considers non-$q$-balanced $\Iscr^+,\Iscr^-,\delta$ and explicitly constructs a
certain graph determining a $q$-matrix for which~\refeq{q_ident} is violated
when $\Kset$ is a field of characteristic 0 and $q$ is transcendental over
$\Qset$.

The second source of our approach is the path method due to
Casteels~\cite{cast1,cast2}. He associated with an $m\times n$ Cauchon diagram
$C$ of~\cite{cauch} a directed planar graph $G=G_C$ with $m+n$ distinguished
vertices $r_1,\ldots,r_m,c_1,\ldots,c_n$ in which the remaining vertices
correspond to white cells $(i,j)$ in the diagram $C$ and are labeled as
$t_{ij}$. An example is illustrated in the picture.

\vspace{0cm}
\begin{center}
\includegraphics[scale=0.8]{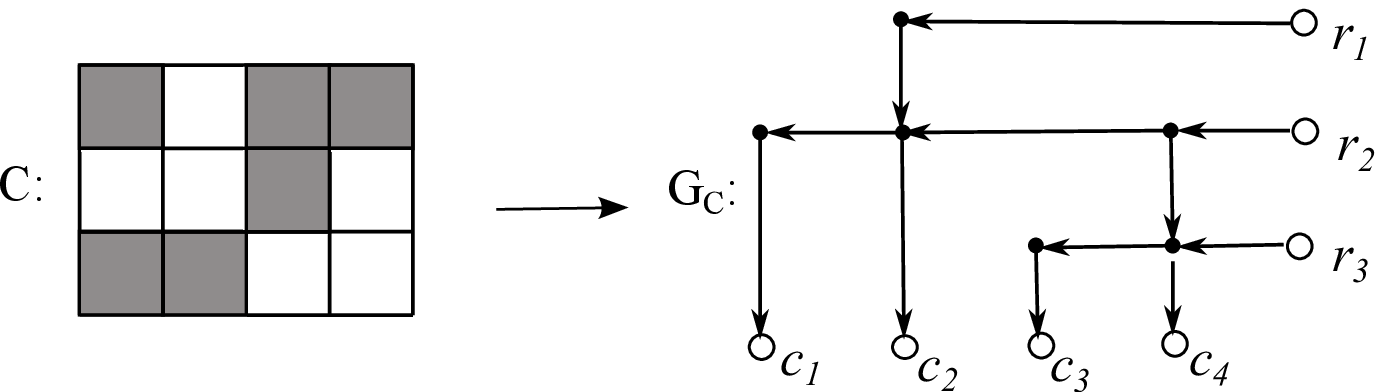}
\end{center}
\vspace{0cm}

\noindent The labels $t_{ij}$, regarded as indeterminates, are assumed to
(quasi)commute as
  \begin{eqnarray}
  t_{ij}t_{i'j'}=&qt_{i'j'}t_{ij}&\quad \mbox{if either $i=i'$ and $j<j'$,
     or $i<i'$ and $j=j'$},  \label{eq:trelat} \\
     =&t_{i'j'}t_{ij}& \quad\mbox{otherwise} \nonumber
  \end{eqnarray}
(which is viewed ``simpler'' than~\refeq{xijkl}). These labels determine
weights of edges and, further, weights of paths of $G$. The latter give rise to
the \emph{path matrix} $P_G$ of size $m\times n$, of which $(i,j)$-th entry is
the sum of weights of paths starting at $r_i$ and ending at $c_i$.

The path matrix $P_G=(p_{ij})$ has three important properties. (i) It is a
$q$-matrix, and therefore, $x_{ij}\mapsto p_{ij}$ gives a homomorphism of
$\Rscr$ to the corresponding algebra $\Rscr_G$ generated by the $p_{ij}$. (ii)
$P_G$ admits an analog of Lindstr\"om's Lemma~\cite{lind}: for any
$I\subseteq[m]$ and $J\subseteq[n]$ with $|I|=|J|$, the minor $[I|J]_q$ of
$P_G$ can be expressed as the sum of weights of systems of \emph{disjoint
paths} from $\{r_1\colon i\in I\}$ to $\{c_j\colon j\in J\}$ in $G$. (iii) From
Cauchon's Algorithm~\cite{cauch} interpreted in graph terms
in~\cite{cast1,cast2} it follows that: if the diagram $C$ is maximal (i.e., has
no black cells), then $P_G$ becomes a \emph{generic $q$-matrix}, see
Corollary~3.2.5 in~\cite{cast2}.

In this paper we consider a more general class of planar graphs $G$ with
horizontal and vertical edges, called \emph{SE-graphs}, and show that they
satisfy the above properties~(i)--(ii) as well. Our goal is to characterize
quadratic identities just for the class of path matrices of SE-graphs $G$.
Since this class contains a generic $q$-matrix, the identities are
automatically valid in $\Rscr$.

We take an advantage from the representation of $q$-minors of path matrices via
systems of disjoint paths, or \emph{flows} in our terminology, and the desired
results are obtained by applying a combinatorial machinery of handling flows in
SE-graphs. Our method of establishing or verifying one or another identity
admits a rather transparent implementation and we illustrate the method by
enlightening graphical diagrams.
  \smallskip

The paper is organized as follows. Section~\SEC{prelim} contains basic
definitions and backgrounds. Section~\SEC{flows} defines flows and path
matrices for SE-graphs and states Lindstr\"om's type theorem for them.
Section~\SEC{double} is devoted to crucial ingredients of the method. It
describes \emph{exchange operations} on \emph{double flows} (pairs of flows
related to corteges $(I|J,I'|J')$) and expresses such operations on the
language of planar matchings. The main working tool of the whole proof, stated
in this section and proved in Appendix~B, is Theorem~\ref{tm:single_exch}
giving a $q$-relation between double flows before and after an ordinary
exchange operation. Using this, Section~\SEC{q_relat} proves the sufficiency in
the main result: \refeq{q_ident} is valid if the corresponding
$\Iscr^+,\Iscr^-,\delta$ are $q$-balanced (Theorem~\ref{tm:suff_q_bal}).

Section~\SEC{examples} is devoted to illustrations of our method. It explains
how to obtain, with the help of the method, rather transparent proofs for
several representative examples of quadratic identities, in particular: (a) the
pure commutation of $[I|J]_q$ and $[I'|J']_q$ when $I'\subset I$ and $J'\subset
J$; (b) a quasicommutation of \emph{flag $q$-minors} $[I]_q$ and $[J]_q$ as in
Leclerc-Zelevinsky's theorem~\cite{LZ}; (c) identities on flag $q$-minors
involving triples $i<j<k$ and quadruples $i<j<k<\ell$; (d) Dodgson's type
identity; (e) two general quadratic identities on flag $q$-minors
from~\cite{LR,TT} occurring in descriptions of quantized Grassmannians and flag
varieties. In Section~\SEC{necess} we prove the necessity of the
$q$-balancedness condition for validity of quadratic identities
(Theorem~\ref{tm:nec_q_bal}); here we rebuild a corresponding construction
from~\cite{DKK} to obtain, in case of the non-$q$-balancedness, an SE-graph $G$
such that the identity for its path matrix is false (in a special case of
$\Kset$ and $q$). Section~\SEC{concl} poses the problem: when an identity in
the commutative case, such as~\refeq{com_ident}, can be turned, by choosing an
appropriate $\delta$, into the corresponding identity for the quantized case?
For example, this is impossible for the trivial identity $[I]\,[J]=[J]\,[I]$
with usual flag minors when $I,J$ are not weakly separated, as is shown
in~\cite{LZ}. Also this section gives a generalization of Leclerc-Zelevinsky's
quasicommutation theorem to arbitrary (non-flag) quantum minors
(Theorem~\ref{tm:generalLZ}) and discusses additional results.

Finally, Appendix~A exhibits several auxiliary lemmas needed to us and proves
the above-mentioned Lindstr\"om's type result for SE-graphs, and Appendix~B
gives the proof of Theorem~\ref{tm:single_exch} (which is rather technical).


\section{\Large Preliminaries}  \label{sec:prelim}


\subsection{Paths in graphs.}  \label{ssec:basterm}
Throughout, by a \emph{graph} we mean a directed graph. A \emph{path} in a
graph $G=(V,E)$ (with vertex set $V$ and edge set $E$) is a sequence
$P=(v_0,e_1,v_1,\ldots,e_k,v_k)$ such that each $e_i$ is an edge connecting
vertices $v_{i-1},v_i$. An edge $e_i$ is called \emph{forward} if it is
directed from $v_{i-1}$ to $v_i$, denoted as $e_i=(v_{i-1},v_i)$, and \emph{
backward} otherwise (when $e_i=(v_i,v_{i-1})$). The path $P$ is called {\em
directed} if it has no backward edge, and {\em simple} if all vertices $v_i$
are different. When $k>0$ and $v_0=v_k$, ~$P$ is called a \emph{cycle}, and
called a \emph{simple cycle} if, in addition, $v_1,\ldots,v_k$ are different.
When it is not confusing, we may use for $P$ the abbreviated notation via
vertices: $P=v_0v_1\ldots v_k$, or edges $P=e_1e_2\ldots e_k$.

Also, using standard terminology in graph theory, for a directed edge
$e=(u,v)$, we say that $e$ \emph{leaves} $u$ and \emph{enters} $v$, and that
$u$ is the \emph{tail} and $v$ is the \emph{head} of $e$.


\subsection{Quantum matrices.}  \label{ssec:quant_matr}
~It will be convenient for us to visualize matrices in the Cartesian form: for
an $m\times n$ matrix $A=(a_{ij})$, the row indices $i=1,\ldots,m$ are assumed
to increase upwards, and the column indices $j=1,\ldots,n$ from left to right.

As mentioned above, we deal with the \emph{quantized coordinate ring}
$\Rscr=\Oscr_q(\Mscr_{m,n}(\Kset))$ generated by indeterminates $x_{ij}$
satisfying relations~\refeq{xijkl}, shortly called the algebra of $m\times n$
\emph{quantum matrices}. A somewhat ``simpler'' object is the \emph{quantum
affine space} $\bar\Rscr$, the $\Kset$-algebra generated by indeterminates
$t_{ij}$ ($i\in[m],\, j\in[n]$) subject to relations~\refeq{trelat}.


\subsection{$q$-minors.}  \label{ssec:quant_minor} For an $m\times n$
matrix $A=(a_{ij})$, we denote by $A(I|J)$ the submatrix of $A$ whose rows are
indexed by $I\subseteq[m]$, and columns by $J\subseteq[n]$. Let $|I|=|J|=:k$,
and let $I$ consist of $i_1<\cdots<i_k$ and $J$ consist of $j_1<\cdots<j_k$.
Then the $q$-\emph{determinant} of $A(I|J)$, or the $q$-\emph{minor} of $A$ for
$I|J$, is defined as
  \begin{equation} \label{eq:qminor}
  [I|J]_{A,q}:=\sum_{\sigma\in S_k} (-q)^{\ell(\sigma)}
    \prod_{d=1}^{k} a_{i_dj_{\sigma(d)}},
    \end{equation}
where, in the noncommutative case, the product under $\prod$ is ordered (from
left to right) by increasing $d$, and $\ell(\sigma)$ is the \emph{length}
(number of inversions) of a permutation $\sigma$. The terms $A$ and/or $q$ in
$[I|J]_{A,q}$ may be omitted when they are clear from the context.


\subsection{SE-graphs.}  \label{ssec:SE}
~A graph $G=(V,E)$ of this sort (also denoted as $(V,E;R,C)$) satisfies the
following conditions:

(SE1) $G$ is planar (with a fixed layout in the plane);

(SE2) $G$ has edges of two types: \emph{horizontal} edges, or \emph{H-edges},
which are directed to the right, and \emph{vertical} edges, or \emph{V-edges},
which are directed downwards (so each edge points to either \emph{south} or
\emph{east}, justifying the term ``SE-graph'');

(SE3) $G$ has two distinguished subsets of vertices: set $R=\{r_1,\ldots,r_m\}$
of \emph{sources} and set $C=\{c_1,\ldots,c_n\}$ of \emph{sinks}; moreover,
$r_1,\ldots,r_m$ are disposed on a vertical line, in this order upwards, and
$c_1,\ldots,c_n$ are disposed on a horizontal line, in this order from left to
right; the sources (sinks) are incident only with H-edges (resp. V-edges);

(SE4) each vertex of $G$ belongs to a directed path from $R$ to $C$.

We denote by $W=W_G$ the set $V-(R\cup C)$ of \emph{inner} vertices of $G$. An
example of SE-graphs with $m=3$ and $n=4$ is drawn in the picture:

\vspace{0cm}
\begin{center}
\includegraphics[scale=0.8]{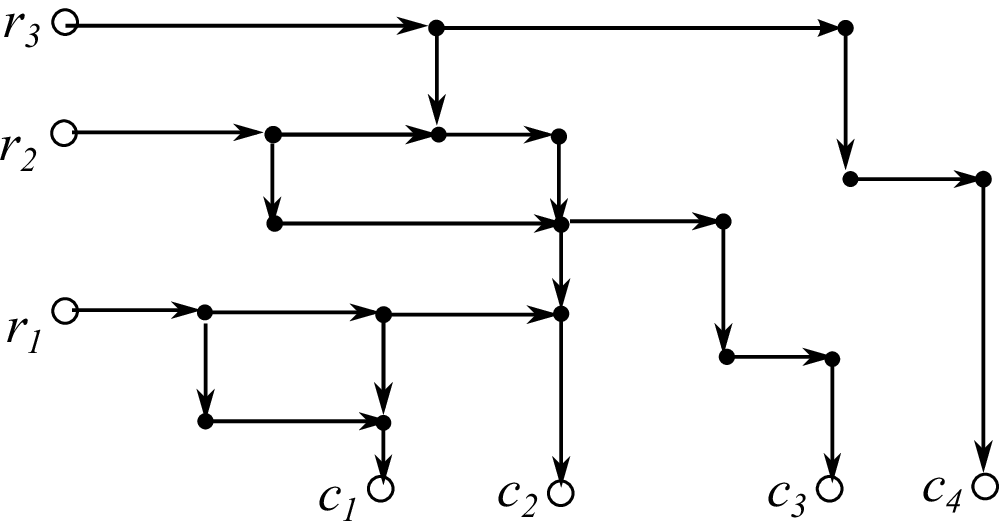}
\end{center}
\vspace{0cm}

\noindent\underline{\em Remark 1.} ~A special case of SE-graphs is formed by
those corresponding to \emph{Cauchon graphs} introduced in~\cite{cast1} (which
are associated with Cauchon diagrams~\cite{cauch}). In this case,
$R=\{(0,i)\colon i\in[m]\}$, $C=\{(j,0)\colon j\in[n]\}$, and $W\subseteq
[m]\times [n]$. (The correspondence with the definition in~\cite{cast1} is
given by $(i,j)\mapsto (m+1-i,n+1-j)$ and $q\mapsto q^{-1}$.) When $W=[m]\times
[n]$ (equivalently: when the Cauchon diagram has no black cells), we refer to
such a graph as the \emph{extended $(m,n)$-grid} and denote it by
$\Gamma_{m,n}$.
\medskip

Each inner vertex $v\in W$ is regarded as a \emph{generator}, and we assign the
weight $w(e)$ to each edge $e=(u,v)\in E$ in a way similar to that for Cauchon
graphs in~\cite{cast1}, namely:
  \begin{numitem1} \label{eq:edge_weight}
  \begin{itemize}
  \item[(i)] $w(e):=v$ if $u\in R$;
  \item[(ii)] $w(e):=u^{-1}v$ if $e$ is an H-edge and $u,v\in W$;
  \item[(iii)] $w(e):=1$ if $e$ is a V-edge.
    \end{itemize}
  \end{numitem1}

This gives rise to defining the weight $w(P)$ of a directed path
$P=e_1e_2\ldots e_k$ (written in the edge notation) in $G$, to be the ordered
(from left to right) product
  \begin{equation} \label{eq:wP}
  w(P)=w(e_1)w(e_2)\cdots w(e_k).
  \end{equation}

Then $w(P)$ is a Laurent monomial in elements of $W$. Note that when $P$ begins
in $R$ and ends in $C$, its weight can also be expressed in the following
useful form; cf.~\cite[Prop.~3.1.8]{cast2}. Let $u_1,v_1,u_2,v_2,\ldots,
u_{d-1},v_{d-1},u_d$ be the sequence of vertices where $P$ makes turns; namely,
$P$ changes the horizontal direction to the vertical one at each $u_i$, and
conversely at each $v_i$. Then (due to the ``telescopic effect'' caused
by~\refeq{edge_weight}(ii)),
  \begin{equation} \label{eq:telescop}
  w(P)=u_1v_1^{-1}u_2v_2^{-1}\cdots u_{d-1}v_{d-1}^{-1} u_d.
  \end{equation}

We assume that the generators $W$ obey (quasi)commutation laws somewhat similar
to those in~\refeq{trelat}; namely, for distinct $u,v\in W$,
  \smallskip

  (G1) if there is a directed \emph{horizontal} path from $u$ to $v$ in $G$, then
$uv=qvu$;

  (G2) if there is a directed \emph{vertical} path from $u$ to $v$ in $G$, then
$vu=quv$;

  (G3) otherwise $uv=vu$.


\section{\Large Path matrix and flows}  \label{sec:flows}

As mentioned in the Introduction, it is shown in~\cite{cast1} that the path
matrix associates with a Cauchon graph $G$ has a nice property of Lindstr\"om's
type, saying that $q$-minors of this matrix correspond to appropriate systems
of disjoint paths in $G$. We will show that this property is extended to the
SE-graphs.

Let $G=(V,E)$ be an SE-graph with sources $R=(r_1,\ldots,r_m)$ and sinks
$C=(c_1,\ldots,c_n)$, and let $w=w_G$ denote the edge weights in $G$ defined
by~\refeq{edge_weight}.
  \medskip

\noindent\textbf{Definition.} The \emph{path matrix} $\Path=\Path_G$ associated
with $G$ is the $m\times n$ matrix whose entries are defined by
  \begin{equation} \label{eq:Mat}
  \Path(i|j):=\sum\nolimits_{P\in\Phi_G(i|j)} w(P), \qquad (i,j)\in [m]\times [n],
  \end{equation}
where $\Phi_G(i|j)$ is the set of directed paths from $r_i$ to $c_j$ in $G$. In
particular, $\Path(i|j)=0$ if $\Phi_G(i|j)=\emptyset$.
  \smallskip

Thus, the entries of $\Path_G$ belong to the $\Kset$-algebra $\Lscr_G$ of
Laurent polynomials generated by the set $W$ of inner vertices of $G$ subject
to relations~(G1)--(G3).
\smallskip

\noindent\textbf{Definition.} Let $\Escr^{m,n}$ denote the set of pairs $(I|J)$
such that $I\subseteq [m]$, $J\subseteq [n]$ and $|I|=|J|$. Borrowing
terminology from~\cite{DKK}, for $(I|J)\in\Escr^{m,n}$, a set $\phi$ of
pairwise disjoint directed paths from the source set $R_I:=\{r_i\colon i\in
I\}$ to the sink set $C_J:=\{c_j\colon j\in J\}$ in $G$ is called an
$(I|J)$-\emph{flow}.
  \medskip

The set of $(I|J)$-flows $\phi$ in $G$ is denoted by $\Phi(I|J)=\Phi_G(I|J)$.
We usually assume that the paths forming a flow $\phi$ are ordered by
increasing the source indices. Namely, if $I$ consists of $i(1)<i(2)<\cdots<
i(k)$ and $J$ consists of $j(1)<j(2)<\cdots<j(k)$, then $\ell$-th path $P_\ell$
in $\phi$ begins at $r_{i(\ell)}$, and therefore, $P_\ell$ ends at
$c_{j(\ell)}$ (which easily follows from the planarity of $G$, the ordering of
sources and sinks in the boundary of $G$ and the fact that the paths in $\phi$
are disjoint). We write $\phi=(P_1,P_2,\ldots,P_k)$ and (similar to path
systems in~\cite{cast1}) define the weight of $\phi$ to be the ordered product
  \begin{equation} \label{eq:w_phi}
  w(\phi)=w(P_1)w(P_2)\cdots w(P_k).
  \end{equation}

Then the desired $q$-analog of Lindstr\"om's Lemma expresses $q$-minors of path
matrices via flows as follows.
  \begin{theorem} \label{tm:Linds}
For the path matrix $\Path=\Path_G$ of an $(m,n)$ SE-graph $G$ and for any
$(I|J)\in \Escr^{m,n}$, there holds
  \begin{equation} \label{eqn:Lind}
[I|J]_{\Path,q}=\sum\nolimits_{\phi\in\Phi(I|J)} w(\phi).
  \end{equation}
  \end{theorem}

A proof of this theorem, which is close to that in~\cite{cast1}, is given in
Appendix~A.

An important fact is that the entries of $\Path_G$ obey the (quasi)commutation
relations similar to those for the canonical generators $x_{ij}$ of the quantum
algebra $\Rscr$ given in~\refeq{xijkl}. It is exhibited in the following
assertion, which is known for the path matrices of Cauchon graphs due
to~\cite{cast1} (where it is proved by use of the ``Cauchon's deleting
derivation algorithm in reverse''~\cite{cauch}).
  \begin{theorem} \label{tm:Path-A}
~For an SE-graph $G$, the entries of its path matrix $\Path_G$ satisfy Manin's
relations.
 \end{theorem}

We will show this in Section~\SSEC{abcd} as an easy application of our
flow-matching method. This assertion implies that the map $x_{ij}\mapsto
\Path_G(i|j)$ determines a homomorphism of $\Rscr$ to the subalgebra $\Rscr_G$
of $\Lscr_G$ generated by the entries of $\Path_G$, i.e., $\Path_G$ is a
$q$-matrix for any SE-graph $G$. In an especial case of $G$, a sharper result,
attributed to Cauchon and Casteels, is as follows.
  \begin{theorem}[\cite{cauch,cast2}] \label{tm:Cach-Cast}
If $G=\Gamma_{m,n}$ (the extended $m\times n$-grid defined in Remark~1), then
$\Path_G$ is a generic $q$-matrix, i.e., $x_{ij}\mapsto \Path_G(i|j)$ gives an
injective map of $\Rscr$ to $\Lscr_G$.
\end{theorem}

Due to this important property, the quadratic relations that are valid
(universally) for $q$-minors of path matrices of SE-graphs turn out to be
automatically valid for the algebra $\Rscr$ of quantum matrices, and vice
versa.


\section{\Large Double flows, matchings, and exchange operations}  \label{sec:double}

Quadratic identities of our interest in this paper involve products of quantum
minors of the form $[I|J][I'|J']$, where $(I|J),(I'|J')\in \Escr^{m,n}$. This
leads to a proper study of ordered pairs of flows $\phi\in \Phi(I|J)$ and
$\phi'\in \Phi(I'|J')$ in an SE-graph $G$ (in light of Theorem~\ref{tm:Linds}).

We need some definitions and conventions, borrowing terminology
from~\cite{DKK}. Given $I,J,I',J',\phi,\phi'$ as above, we call the pair
$(\phi,\phi')$ a \emph{double flow} in $G$. Let
   \begin{gather}
   \Iw:=I-I',\quad \Jw:=J-J',\quad \Ib:=I'-I,\quad \Jb:=J'-J,
                          \label{eq:white-black} \\
   \Yr:=\Iw\cup\Ib\quad \mbox{and}\quad \Yc:=\Jw\cup\Jb. \nonumber
   \end{gather}

Note that $|I|=|J|$ and $|I'|=|J'|$ imply that $|\Yr|+|\Yc|$ is even and
  \begin{equation} \label{eq:balancIJ}
 |\Iw|-|\Ib|=|\Jw|-|\Jb|.
  \end{equation}

We refer to the quadruple $(I|J,I'|J')$ as above as a \emph{cortege}, and to
$(\Iw,\Ib,\Jw,\Jb)$ as the \emph{refinement} of $(I|J,I'|J')$, or as a
\emph{refined cortege}.

It is convenient for us to interpret $\Iw$ and $\Ib$ as the sets of
\emph{white} and \emph{black} elements of $\Yr$, respectively, and similarly
for $\Jw,\Jb,\Yc$, and visualize these objects by use of a \emph{circular
diagram} $D$ in which the elements of $\Yr$ (resp. $\Yc$) are disposed in the
increasing order from left to right in the upper (resp. lower) half of a
circumference $O$. For example if, say, $\Iw=\{3\}$, $\Ib=\{1,4\}$,
$\Jw=\{2',5'\}$ and $\Jb=\{3',6',8'\}$, then the diagram is viewed as in the
left fragment of the picture below. (Sometimes, to avoid a possible mess
between elements of $\Yr$ and $\Yc$, and when it leads to no confusion, we
denote elements of $\Yc$ with primes.)

\vspace{0cm}
\begin{center}
\includegraphics[scale=0.9]{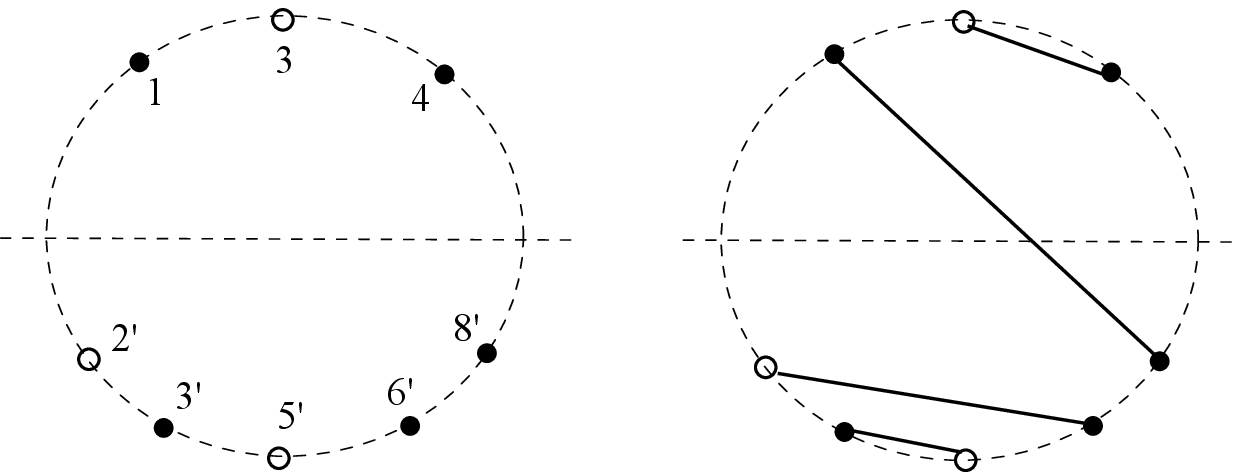}
\end{center}
\vspace{0cm}

Let $M$ be a partition of $\Yr\sqcup \Yc$ into 2-element sets (recall that
$A\sqcup B$ denotes the disjoint union of sets $A,B$). We refer to $M$ as a
\emph{perfect matching} on $\Yr\sqcup \Yc$, and to its elements as
\emph{couples}. More specifically, we say that $\pi\in M$ is: an
$R$-\emph{couple} if $\pi\subseteq \Yr$, a $C$-\emph{couple} if $\pi\subseteq
\Yc$, and an $RC$-\emph{couple} if $|\pi\cap \Yr|=|\pi\cap \Yc|=1$ (as though
$\pi$ ``connects'' two sources, two sinks, and one source and one sink,
respectively). \smallskip

\noindent\textbf{Definition.} A (perfect) matching $M$ as above is called a
\emph{feasible} matching for $(\Iw,\Ib,\Jw,\Jb)$ (and for $(I|J,I'|J')$) if:
  \begin{numitem1} \label{eq:feasM}
  \begin{itemize}
\item[(i)] for each $\pi=\{i,j\}\in M$, the elements $i,j$ have different
colors if $\pi$ is an $R$- or $C$-couple, and have the same color if $\pi$ is
an $RC$-couple;
\item[(ii)] $M$ is \emph{planar}, in the sense that the chords connecting the
couples in the circumference $O$ are pairwise non-intersecting.
  \end{itemize}
  \end{numitem1}

The set of feasible matchings for $(\Iw,\Ib,\Jw,\Jb)$ is denoted by
$\Mscr_{\Iw,\Ib,\Jw,\Jb}$ and may also be denoted as $\Mscr(I|J,I'|J')$. This
set is nonempty unless $\Yr\sqcup\Yc=\emptyset$. (A proof: a feasible matching
can be constructed recursively as follows. Let for definiteness
$|\Iw|\ge|\Ib|$. If $\Ib\ne\emptyset$, then choose $i\in\Iw$ and $j\in\Ib$ with
$|i-j|$ minimum, form the $R$-couple $\{i,j\}$ and delete $i,j$. And so on
until $\Ib$ becomes empty. Act similarly for $\Jw$ and $\Jb$. Eventually, in
view of~\refeq{balancIJ}, we obtain $\Ib=\Jb=\emptyset$ and $|\Iw|=|\Jw|$. Then
we form corresponding white $RC$-couples.)

The right fragment of the above picture illustrates an instance of feasible
matchings.

Return to a double flow $(\phi,\phi')$ as above. Our aim is to associate to it
a feasible matching for $(\Iw,\Ib,\Jw,\Jb)$.

To do this, we write $V_\phi$ and $E_\phi$, respectively, for the sets of
vertices and edges of $G$ occurring in $\phi$, and similarly for $\phi'$. An
important role will be played by the subgraph $\langle U\rangle$ of $G$ induced
by the set of edges
   $$
   U:=E_\phi\triangle E_{\phi'}
   $$
(where $A\triangle B$ denotes $(A-B)\cup(B-A)$). Note that a vertex $v$ of
$\langle U\rangle$ has degree 1 if $v\in R_{\Iw}\cup R_{\Ib}\cup C_{\Jw}\cup
C_{\Jb}$, and degree 2 or 4 otherwise. We slightly modify $\langle U\rangle$ by
splitting each vertex $v$ of degree 4 in $\langle U\rangle$ (if any) into two
vertices $v',v''$ disposed in a small neighborhood of $v$ so that the edges
entering (resp. leaving) $v$ become entering $v'$ (resp. leaving $v''$); see
the picture.

\vspace{0cm}
\begin{center}
\includegraphics{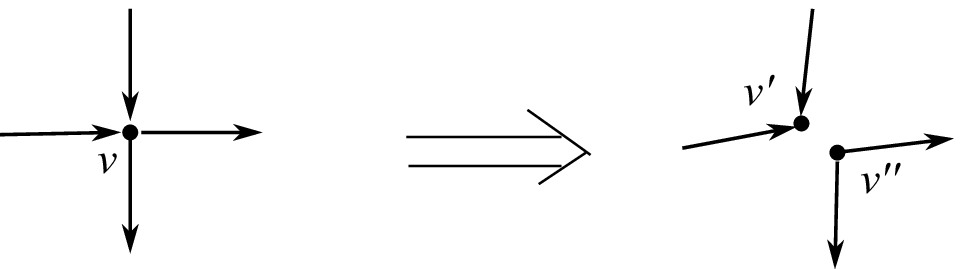}
\end{center}
\vspace{0cm}

The resulting graph, denoted as $\langle U\rangle '$, is planar and has
vertices of degree only 1 and 2. Therefore, $\langle U\rangle'$ consists of
pairwise disjoint (non-directed) simple paths $P'_1,\ldots,P'_k$ (considered up
to reversing) and, possibly, simple cycles $Q'_1,\ldots,Q'_d$. The
corresponding images of $P'_1,\ldots,P'_k$ (resp. $Q'_1,\ldots,Q'_d$) give
paths $P_1,\ldots,P_k$ (resp. cycles $Q_1,\ldots,Q_d$) in $\langle U\rangle$.
When $\langle U\rangle$ has vertices of degree 4, some of the latter paths and
cycles may be self-intersecting and may ``touch'', but not ``cross'', each
other.
  \begin{lemma} \label{lm:P1Pk}
{\rm (i)} $k=(|\Iw|+|\Ib|+|\Jw|+|\Jb|)/2$;

{\rm(ii)} the set of endvertices of $P_1,\ldots,P_k$ is $R_{\Iw\cup\Ib}\cup
C_{\Jw\cup\Jb}$; moreover, each $P_i$ connects either $R_{\Iw}$ and $R_{\Ib}$,
or $C_{\Jw}$ and $C_{\Jb}$, or $R_{\Iw}$ and $C_{\Jw}$, or $R_{\Ib}$ and
$C_{\Jb}$;

{\rm(iii)} in each path $P_i$, the edges of $\phi$ and the edges of $\phi'$
have different directions (say, the former edges are all forward, and the
latter ones are all backward).
  \end{lemma}
  \begin{proof}
(i) is trivial, and (ii) follows from (iii) and the fact that the sources $r_i$
(resp. sinks $c_j$) have merely leaving (resp. entering) edges. In its turn,
(iii) easily follows by considering a common inner vertex $v$ of a directed
path $K$ in $\phi$ and a directed path $L$ in $\phi'$. Let $e,e'$ (resp.
$u,u'$) be the edges of $K$ (resp. $L$) incident to $v$. Then: if
$\{e,e'\}=\{u,u'\}$, then $v$ vanishes in $\langle U\rangle$. If $e=u$ and
$e'\ne u'$, then either both $e',u'$ enter $v$, or both $e',u'$ leave $v$;
whence $e',u'$ are consecutive and differently directed edges of some path
$P_i$ or cycle $Q_j$. A similar property holds when
$\{e,e'\}\cap\{u,u'\}=\emptyset$, as being a consequence of splitting $v$ into
two vertices as described.  \end{proof}

Thus, each $P_i$ is represented as a concatenation $P_i^{(1)}\circ
P_i^{(2)}\circ\ldots\circ P_i^{(\ell)}$ of forwardly and backwardly directed
paths which are alternately contained in $\phi$ and $\phi'$, called the
\emph{segments} of $P_i$. We refer to $P_i$ as an \emph{exchange path} (by a
reason that will be clear later). The endvertices of $P_i$ determine, in a
natural way, a pair of elements of $\Yr\sqcup \Yc$, denoted by $\pi_i$. Then
$M:=\{\pi_1,\ldots,\pi_k\}$ is a perfect matching on $\Yr\sqcup \Yc$. Moreover,
it is a feasible matching, since~\refeq{feasM}(i) follows from
Lemma~\ref{lm:P1Pk}(ii), and~\refeq{feasM}(ii) is provided by the fact that
$P'_1,\ldots,P'_k$ are pairwise disjoint simple paths in $\langle U\rangle'$.

We denote $M$ as $M(\phi,\phi')$, and for $\pi\in M$, denote the exchange path
$P_i$ corresponding to $\pi$ (i.e., $\pi=\pi_i$) by $P(\pi)$.

\begin{corollary} \label{cor:Mphiphip}
$M(\phi,\phi')\in\Mscr_{\Iw,\Ib,\Jw,\Jb}$.
  \end{corollary}

\begin{figure}[htb]
\vspace{0.3cm}
\begin{center}
\includegraphics{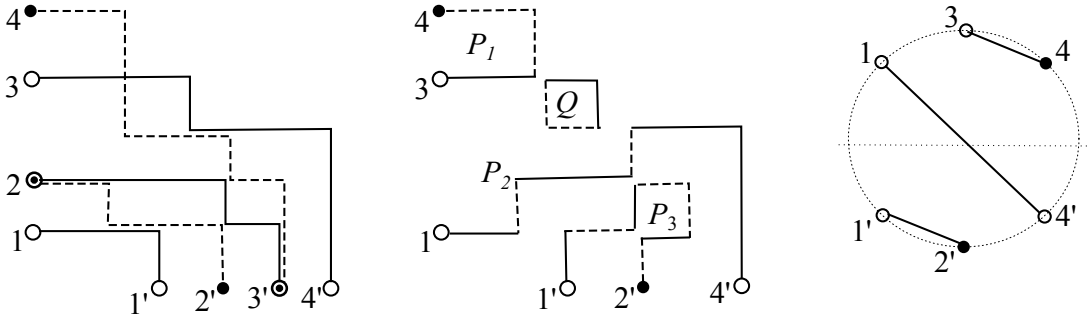}
\end{center}
\vspace{-0.3cm}
 \caption{ flows $\phi$ and $\phi'$ (left);  $\langle E_\phi\triangle
 E_{\phi'}\rangle$ (middle); $M(\phi,\phi')$ (right)}
 \label{fig:phi}
  \end{figure}

Figure~\ref{fig:phi} illustrates an instance of $(\phi,\phi')$ for
$I=\{1,2,3\}$, $J=\{1',3',4'\}$, $I'=\{2,4\}$, $J'=\{2',3'\}$. Here $\phi$ and
$\phi'$ are drawn by solid and dotted lines, respectively (in the left
fragment), the subgraph $\langle E_\phi\triangle E_{\phi'}\rangle$ consists of
three paths and one cycle (in the middle), and the circular diagram illustrates
$M(\phi,\phi')$ (in the right fragment).

\noindent \textbf{Flow exchange operation.} It rearranges a given double flow
$(\phi,\phi')$ for $(I|J,\,I'|J')$ into another double flow $(\psi,\psi')$ for
some cortege $(\tilde I|\tilde J,\,\tilde I'|\tilde J')$, as follows. Fix a
submatching $\varPi\subseteq M(\phi,\phi')$, and combine the exchange paths
concerning $\varPi$, forming the set of edges
   $$
   \Escr:= \cup (E_{P(\pi)}\colon \pi\in\varPi).
   $$
(where $E_P$ denotes the set of edges in a path $P$).
  \begin{lemma} \label{lm:phi-psi}
Let $V_{\varPi}:=\cup(\pi\in\varPi)$. Define
  $$
  \tilde I:=I\triangle (V_{\varPi}\cap \Yr),\;\; \tilde I':=I'\triangle (V_{\varPi}\cap\Yr), \;\;
  \tilde J:=J\triangle (V_{\varPi}\cap \Yc),\;\; \tilde J':=J'\triangle (V_{\varPi} \cap\Yc).
  $$
Then the subgraph $\psi$ induced by $E_\phi\triangle\Escr$ gives a $(\tilde
I|\tilde J)$-flow, and the subgraph $\psi'$ induced by $E_{\phi'}\triangle
\Escr$ gives a $(\tilde I'|\tilde J')$-flow in $G$. Furthermore, $E_\psi\cup
E_{\psi'}= E_\phi\cup E_{\phi'}$,\; $E_\psi\triangle E_{\psi'}= E_\phi\triangle
E_{\phi'}$ ($=U$), and $M(\psi,\psi')=M(\phi,\phi')$.
  \end{lemma}
  \begin{proof}
Consider a path $P=P(\pi)$ for $\pi\in\varPi$, and let $P$ consist of segments
$P^{(1)},P^{(2)},\ldots,P^{(\ell)}$. Let for definiteness the segments
$P^{(d)}$ with $d$ odd concern $\phi$, and denote by $v_d$ the common endvertex
of $P^{(d)}$ and $P^{(d+1)}$. Under the operation $E_\phi\mapsto
E_\phi\triangle E_P$ the pieces $P^{(1)},P^{(3)},\ldots$ in $\phi$ are replaced
by $P^{(2)},P^{(4)},\ldots$. In its turn, $E_{\phi'}\mapsto E_{\phi'}\triangle
E_P$ replaces the pieces $P^{(2)},P^{(4)},\ldots$ in $\phi'$ by
$P^{(1)},P^{(3)},\ldots$.

By Lemma~\ref{lm:P1Pk}(iii), for each $d$, the edges of $P^{(d)},P^{(d+1)}$
incident to $v_d$ either both enter or both leave $v_d$. Also each intermediate
vertex of any segment $P^{(d)}$ occurs in exactly one flow among $\phi,\phi'$.
These facts imply that under the above operations with $P$ the flow $\phi$
(resp. $\phi'$) is transformed into a set of pairwise disjoint directed paths
(a flow) going from $R_{I\triangle(\pi\cap \Yr)}$ to $C_{J\triangle(\pi\cap
\Yc)}$ (resp. from $R_{I'\triangle(\pi\cap \Yr)}$ to $C_{J'\triangle(\pi\cap
\Yc)}$).

Doing so for all $P(\pi)$ with $\pi\in\varPi$, we obtain flows $\psi,\psi'$
from $R_{\tilde I}$ to $C_{\tilde J}$ and from $R_{\tilde I'}$ to $C_{\tilde
J'}$, respectively. The equalities in the last sentence of the lemma are easy.
  \end{proof}

We call the transformation $(\phi,\phi')\stackrel{\varPi}\longmapsto
(\psi,\psi')$ in this lemma the \emph{flow exchange operation} for
$(\phi,\phi')$ using $\varPi\subseteq M(\phi,\phi')$ (or using $\{P(\pi)\colon
\pi\in \varPi\}$). Clearly the exchange operation applied to $(\psi,\psi')$
using the same $\varPi$ returns $(\phi,\phi')$. The picture below illustrates
flows $\psi,\psi'$ obtained from $\phi,\phi'$ in Fig.~\ref{fig:phi} by the
exchange operations using the single path $P_2$ (left) and the single path
$P_3$ (right).

\vspace{0.cm}
\begin{center}
\includegraphics{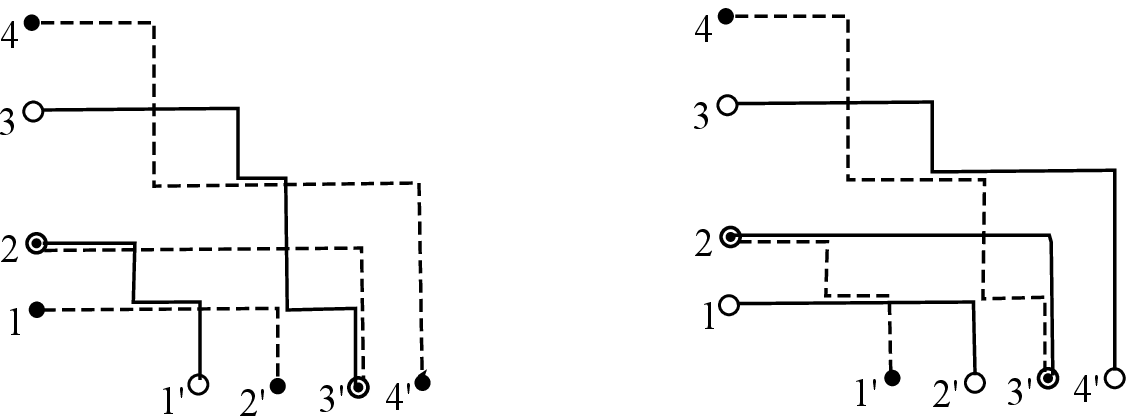}
\end{center}
\vspace{-0.3cm}

So far our description has been close to that given for the commutative case
in~\cite{DKK}. From now on we will essentially deal with the quantum version.
The next theorem will serve the main working tool in our arguments; its proof
appealing to a combinatorial techniques on paths and flows is given in
Appendix~B.
  \begin{theorem} \label{tm:single_exch}
Let $\phi$ be an $(I|J)$-flow, and $\phi'$ an $(I'|J')$-flow in $G$. Let
$(\psi,\psi')$ be the double flow obtained from $(\phi,\phi')$ by the flow
exchange operation using a single couple $\pi=\{f,g\}\in M(\phi,\phi')$. Then:

{\rm(i)} when $\pi$ is an $R$- or $C$-couple and $f<g$,
   \begin{eqnarray*}
   w(\phi)w(\phi')=qw(\psi)w(\psi') \quad &\mbox{in case}&\;\; f\in I\cup J; \\
 w(\phi)w(\phi')=q^{-1}w(\psi)w(\psi') \quad &\mbox{in case}&\;\;
                                                   f\in I'\cup J';
   \end{eqnarray*}

{\rm(ii)} when $\pi$ is an $RC$-couple, $w(\phi)w(\phi')=w(\psi)w(\psi')$.
  \end{theorem}

An immediate consequence from this theorem is the following
  \begin{corollary} \label{cor:gen_exch}
For an $(I|J)$-flow $\phi$ and an $(I'|J')$-flow $\phi'$, let $(\psi,\psi')$ be
obtained from $(\phi,\phi')$ by the flow exchange operation using a set
$\varPi\subseteq M(\phi,\phi')$. Then
  \begin{equation} \label{eq:gen_exch}
   w(\phi)w(\phi')=q^{\zeta^\circ-\zeta^\bullet} w(\psi)w(\psi'),
   \end{equation}
where $\zeta^\circ=\zeta^\circ(I|J,I'|J';\varPi)$ (resp.
$\zeta^\bullet=\zeta^\bullet(I|J,I'|J';\varPi)$) is the number of $R$- or
$C$-couples $\pi=\{f,g\}\in \varPi$ such that $f<g$ and $f\in I\cup J$ (resp.
$f\in I'\cup J'$).
  \end{corollary}
Indeed, the flow exchange operation using the whole $\varPi$ reduces to
performing, step by step, the exchange operations using single couples $\pi\in
\varPi$ (taking into account that for any current double flow $(\eta,\eta')$
occurring in the process, the sets $E_\eta\cup E_{\eta'}$ and $E_\eta\triangle
E_{\eta'}$, as well as the matching $M(\eta,\eta')$, do not change; cf.
Lemma~\ref{lm:phi-psi}). Then~\refeq{gen_exch} follows from
Theorem~\ref{tm:single_exch}.


\section{\Large Quadratic relations}  \label{sec:q_relat}

As before, we consider an SE-graph $G=(V,E;R,C)$ and the weight function $w$
which is initially defined on the edges of $G$ by~\refeq{edge_weight} and then
extends to paths and flows according to~\refeq{wP} and~\refeq{w_phi}. This
gives rise to the minor function on the set $\Escr^{m,n}=\{(I|J)\colon
I\subseteq [m],\; J\in[n],\; |I|=|J|\}$. In this section, based on
Corollary~\ref{cor:gen_exch} describing the transformation of the weights of
double flows under the exchange operation, and developing a $q$-version of the
flow-matching method elaborated for the commutative case in~\cite{DKK}, we
establish sufficient conditions on quadratic relations for $q$-minors of the
matrix $\Path_G$, to be valid independently of $G$ (and some other objects, see
Remark~2 below). Relations of our interest are of the form
  \begin{equation} \label{eq:gen_QR}
  \sum\nolimits_\Iscr q^{\alpha(I|J,I'|J')} [I|J][I'|J'] =
     \sum\nolimits_\Kscr q^{\beta(K|L,K'|L')} [K|L][K'|L'],
     \end{equation}
where $\alpha, \beta$ are integer-valued, $\Iscr$ is a family of corteges
$(I|J,I'|J')\in \Escr^{m,n}\times \Escr^{m,n}$ (with possible multiplicities),
and similarly for $\Kscr$. Cf.~\refeq{q_ident}. We usually assume that $\Iscr$
and $\Kscr$ are \emph{homogeneous}, in the sense that for any
$(I|J,I'|J')\in\Iscr$ and $(K|L,K'|L')\in\Kscr$,
   \begin{equation} \label{eq:homogen}
I\cup I'=K\cup K',\quad J\cup J'=L\cup L',\quad I\cap I'= K\cap K',\quad J\cap
J'=L\cap L'.
  \end{equation}
Moreover, we shall see that only the refinements $(\Iw,\Ib,\Jw,\Jb)$ and
$(\Kw,\Kb,\Lw,\Lb)$ are important, whereas the sets $I\cap I'$ and $J\cap J'$
are, in fact, indifferent. (As before, $\Iw$ means $I-I'$, $\Ib$ means $I'-I$,
and so on.)

To formulate the validity conditions, we need some definitions and notation.

$\bullet$ We say that a tuple $(I|J,I'|J';M)$, where $(I|J,I'|J')\in \Iscr$ and
$M\in\Mscr_{\Iw,\Ib,\Jw,\Jb}$ (cf.~\refeq{feasM}), is a \emph{configuration}
for $\Iscr$. The family of all configurations for $\Iscr$ is denoted by
$\bfC(\Iscr)$. Similarly, we define the family $\bfC(\Kscr)$ of configurations
for $\Kscr$.

$\bullet$ Define $\bfM(\Iscr)$ to be the family of all matchings $M$ occurring
in the members of $\bfC(\Iscr)$, respecting multiplicities (i.e., $\bfM(\Iscr)$
is a multiset). Define $\bfM(\Kscr)$ similarly.
 \medskip

\noindent\textbf{Definition.} Families $\Iscr$ and $\Kscr$ are called
\emph{balanced} (borrowing terminology from~\cite{DKK}) if there exists a
bijection $(I|J,I'|J';M)\stackrel{\gamma}\longmapsto (K|K',L|L';M')$ between
$\bfC(\Iscr)$ and $\bfC(\Kscr)$ such that $M=M'$. In other words, $\Iscr$ and
$\Kscr$ are balanced if $\bfM(\Iscr)=\bfM(\Kscr)$.
 \smallskip

\noindent\textbf{Definition.} We say that families $\Iscr$ and $\Kscr$ along
with functions $\alpha:\Iscr\to\Zset$ and $\beta:\Kscr\to\Zset$ are
$q$-\emph{balanced} if there exists a bijection $\gamma$ as above such that,
for each $(I|J,I'|J';M)\in\bfC(\Iscr)$ and for
$(K|K',L|L';M)=\gamma(I|J,I'|J';M)$, there holds
  \begin{equation} \label{eq:q_balan}
  \beta(K|K',L|L')-\alpha(I|J,I'|J')=\zeta^\circ-
                            \zeta^\bullet.
  \end{equation}
(In particular, $\Iscr,\Kscr$ are balanced.) Here $\zeta^\circ,\zeta^\bullet$
are defined according to Corollary~\ref{cor:gen_exch}. Namely,
$\zeta^\circ=\zeta^\circ(I|J,I'|J';\varPi)$ and
$\zeta^\bullet=\zeta^\bullet(I|J,I'|J';\varPi)$, where $\varPi$ is the set of
couples $\pi\in M$ such that the white/black colors of the elements of $\pi$ in
the refined corteges $(\Iw,\Ib,\Jw,\Jb)$ and $(\Kw,\Kb,\Lw,\Lb)$ are different.
(Then $\zeta^\circ$ ($\zeta^\bullet$) is the number of $R$- and $C$-couples
$\{f,g\}\in\varPi$ with $f<g$ and $f\in \Iw\cup\Jw$ (resp. $f\in \Ib\cup\Jb)$.)
We say that $(\Kw,\Kb,\Lw,\Lb)$ is obtained from $(\Iw,\Ib,\Jw,\Jb)$ by the
\emph{index exchange operation} using $\varPi$, and may write
$\zeta^\circ(\Iw,\Ib,\Jw,\Jb;\varPi)$ for $\zeta^\circ$, and
$\zeta^\bullet(\Iw,\Ib,\Jw,\Jb;\varPi)$ for $\zeta^\bullet$.

  \begin{theorem} \label{tm:suff_q_bal}
Let $\Iscr$ and $\Kscr$ be homogeneous families on
$\Escr^{m,n}\times\Escr^{m,n}$, and let $\alpha:\Iscr\to\Zset$ and
$\beta:\Kscr\to\Zset$. Suppose that $\Iscr,\Kscr,\alpha,\beta$ are
$q$-balanced. Then for any SE-graph $G=(V,E;R,C)$, relation~\refeq{gen_QR} is
valid for $q$-minors of $\Path_G$.
  \end{theorem}
  \begin{proof} ~It is close to the proof for the commutative
case in~\cite[Proposition~3.2]{DKK}.

We fix $G$ and denote by $\Dscr(I|J,I'|J')$ the set of double flows for
$(I|J,I'|J') \in\Iscr\cup\Kscr$ in $G$. A summand concerning
$(I|J,J'|J')\in\Iscr$ in the L.H.S. of~\refeq{gen_QR} can be expressed via
double flows as follows, ignoring the factor of $q^{\alpha(\cdot)}$:
  \begin{multline}  \label{eq:ff}
  [I|J][I'|J']=\left( \sum\nolimits_{\phi\in\Phi_G(I|J)} w(\phi)\right)
     \times\left(\sum\nolimits_{\phi'\in\Phi_G(I'|J')} w(\phi')\right) \\
  =\sum\nolimits_{(\phi,\phi')\in\Dscr(I|J,I'|J')} w(\phi)w(\phi')
                       \qquad\qquad\qquad\qquad\qquad\qquad\\
  =\sum\nolimits_{M\in\Mscr_{\Iw,\Ib,\Jw,\Jb}}
  \sum\nolimits_{(\phi,\phi')\in\Dscr(I|J,I'|J')\,:\, M(\phi,\phi')=M} w(\phi)w(\phi').
  \end{multline}

The summand for $(K|L,K'|L')\in\Kscr$ in the R.H.S. of~\refeq{gen_QR} is
expressed similarly.

Consider a configuration $S=(I|J,I'|J';M)\in\bfC(\Iscr)$ and suppose that
$(\phi,\phi')$ is a double flow for $(I|J,I'|J')$ with $M(\phi,\phi')=M$ (if
such a double flow in $G$ exists). Since $\Iscr,\Kscr,\alpha,\beta$ are
$q$-balanced, $S$ is bijective to some configuration
$S'=(K|L,K'|L';M)\in\bfC(\Kscr)$ satisfying~\refeq{q_balan}. As explained
earlier, the cortege $(K|L,K'|L')$ is obtained from $(I|J,I'|J')$ by the index
exchange operation using some $\varPi\subseteq M$. Then the flow exchange
operation applied to $(\phi,\phi')$ using this $\varPi$ results in a double
flow $(\psi,\psi')$ for $(K|L,K'|L')$ which satisfies relation~\refeq{gen_exch}
in Corollary~\ref{cor:gen_exch}. Comparing~\refeq{gen_exch}
with~\refeq{q_balan}, we observe that
  $$
  q^{\alpha(I|J,I'|J')} w(\phi)w(\phi')= q^{\beta(K|K',L|L')} w(\psi)w(\psi').
  $$
Furthermore, such a map $(\phi,\phi')\mapsto(\psi,\psi')$ gives a bijection
between all double flows concerning configurations in $\bfC(\Iscr)$ and those
in $\bfC(\Kscr)$. Now the desired equality~\refeq{gen_QR} follows by comparing
the last term in expression~\refeq{ff} and the corresponding term in the
analogous expression concerning $\Kscr$.
 \end{proof}

As a consequence of Theorems~\ref{tm:Cach-Cast} and~\ref{tm:suff_q_bal}, the
following result is obtained.
  \begin{corollary} \label{cor:QR_quant_matr}
~If $\Iscr,\Kscr,\alpha,\beta$ as above are $q$-balanced, then
relation~\refeq{gen_QR} is valid for the corresponding minors in the algebra
$\Rscr$ of quantum $m\times n$ matrices.
  \end{corollary}

\noindent\underline{\em Remark 2.} When speaking of a \emph{universal quadratic
identity} of the form~\refeq{gen_QR} with homogeneous $\Iscr$ and $\Kscr$,
abbreviated as a \emph{UQ identity}, we mean that it depends neither on the
graph $G$ nor on the field $\Kset$ and element $q\in \Kset^\ast$, and that the
index sets can be modified as follows. Given $(I|J,I'|J')\in\Iscr$, let
$A:=I\triangle I'$, $B:=J\triangle J'$, $S:=I\cap I'$ and $T:=J\cap J'$ (by the
homogeneity, these sets do not depend on $(I|J,I'|J)\in\Iscr\cup\Kscr$). Take
arbitrary $\tilde m\ge |A|$ and $\tilde n\ge |B|$ and replace $A,B,S,T$ by
disjoint sets $\tilde A,\tilde S\subseteq [\tilde m]$ and disjoint sets $\tilde
B,\tilde T\subseteq [\tilde n]$ such that $|\tilde A|=|A|$, $|\tilde B|=|B|$
and $|\tilde S|-|\tilde T|=|S|-|T|$. Let $\nu:A\to \tilde A$ and
$\mu:B\to\tilde B$ be the order preserving maps. Transform each $(I|J,I'|J')\in
\Iscr$ into $(\tilde I|\tilde J, \tilde I'|\tilde J')$, where
   $$
   \tilde I:=\tilde S \cup \nu(I-S), \;\; \tilde I':=\tilde S \cup\nu(I'-S),
   \;\; \tilde J:=\tilde T \cup \mu(J-T), \;\; \tilde J':=\tilde T
   \cup\mu(J'-T),
   $$
forming a new family $\tilde \Iscr$ on $\Escr^{\tilde m,\tilde n}\times
\Escr^{\tilde m,\tilde n}$. Transform $\Kscr$ into $\tilde K$ in a similar way.
One can see that if $\Iscr,\Kscr,\alpha,\beta$ are $q$-balanced, then so are
$\tilde \Iscr,\tilde \Kscr$, keeping $\alpha,\beta$. Therefore,
if~\refeq{gen_QR} is valid for $\Iscr,\Kscr$, then it is valid for $\tilde
\Iscr,\tilde\Kscr$ as well.
 \medskip

Thus, the condition of $q$-balancedness is sufficient for validity of
relation~\refeq{gen_QR} for minors of any $q$-matrix. In Section~\SEC{necess}
we shall see that this condition is necessary as well
(Theorem~\ref{tm:nec_q_bal}).

One can say that identity~\refeq{gen_QR}, where all summands have positive
signs, is written in the \emph{canonical form}. Sometimes, however, it is more
convenient to consider equivalent identities having negative summands in one or
both sides (e.g. of the form~\refeq{q_ident}). Also one may simultaneously
multiply all summands in~\refeq{gen_QR} by the same degree of $q$.
\medskip

\noindent\underline{\em Remark 3.} A useful fact is that once we are given an
instance of~\refeq{gen_QR}, we can form another identity by changing the
white/black coloring in all refined corteges. More precisely, for a cortege
$S=(I|J,I'|J')$, let us say that the cortege $S^\reve:=(I'|J',I|J)$ is
\emph{reversed} to $S$. Given a family $\Iscr$ of corteges, the \emph{reversed}
family $\Iscr^\reve$ is formed by the corteges reversed to those in $\Iscr$.
Then the following property takes place.
   \begin{prop} \label{pr:reversed}
Suppose that $\Iscr,\Kscr,\alpha,\beta$ are $q$-balanced. Then
$\Iscr^\reve,\Kscr^\reve,-\alpha,-\beta$ are q-balanced as well. Therefore (by
Theorem~\ref{tm:suff_q_bal}),
  \begin{equation} \label{eq:reverQR}
  \sum_{(I|J,I'|J')\in\Iscr} q^{-\alpha(I|J,I'|J')} [I'|J'][I|J] =
     \sum_{(K|L,K'|L')\in \Kscr} q^{-\beta(K|L,K'|L')} [K'|L'][K|L].
     \end{equation}
  \end{prop}
  \begin{proof}
~Let $\gamma:\bfC(\Iscr)\to \bfC(\Kscr)$  be a bijection in the definition of
$q$-balancedness. Then $\gamma$ induces a bijection of $\bfC(\Iscr^\reve)$ to
$\bfC(\Kscr^\reve)$ (also denoted as $\gamma$). Namely, if $\gamma(S;M)=(T;M)$
for $S=(I|J,I'|J')\in\Iscr$ and $T=(K|L,K'|L')\in\Kscr$, then we define
$\gamma(S^\reve;M):=(T^\reve;M)$. When coming from $S$ to $S^\reve$, each $R$-
or $C$-couple $\{i,j\}$ in $M$ changes the colors of both elements $i,j$. This
leads to swapping $\zeta^\circ$ and $\zeta^\bullet$, i.e.,
$\zeta^\circ(S^\reve;\varPi)=\zeta^\bullet (S;\varPi)$ and
$\zeta^\bullet(S^\reve;\varPi)=\zeta^\circ (S;\varPi)$ (where $\varPi$ is the
submatching in $M$ involved in the exchange operation). Now~\refeq{reverQR}
follows from relation~\refeq{q_balan}.
   \end{proof}

Another useful equivalent transformation is given by swapping row and column
indices. Namely, for a cortege $S=(I|J,I'|J')$, the \emph{transposed} cortege
is $S^\top:=(J|I,J'|I')$, and the family $\Iscr^\top$ \emph{transposed} to
$\Iscr$ consists of the corteges $S^\top$ for $S\in\Iscr$, and similarly for
$\Kscr$. One can see that the corresponding values $\zeta^\circ$ and
$\zeta^\bullet$ preserve when coming from $\Iscr$ to $\Iscr^\top$ and from
$\Kscr$ to $\Kscr^\top$, and therefore~\refeq{q_balan} implies the identity
  \begin{equation} \label{eq:transQR}
  \sum_{(I|J,I'|J')\in\Iscr} q^{\alpha(I|J,I'|J')} [J|I][J'|I'] =
     \sum_{(K|L,K'|L')\in \Kscr} q^{\beta(K|L,K'|L')} [L|K][L'|K'].
     \end{equation}
(Note also that~\refeq{transQR} immediately follows from the known fact that
any $q$-minor satisfies the symmetry relation $[J|I]_q=[J|I]_q$.)

We conclude this section with a rather simple algorithm which has as the input
a corresponding quadruple $\Iscr,\Kscr,\alpha,\beta$ and recognizes the
$q$-balanced for it. Therefore, in light of Theorems~\ref{tm:suff_q_bal}
and~\ref{tm:nec_q_bal}, the algorithm decides whether or not the given
quadruple determines a UQ identity of the form~\refeq{gen_QR}.
 \smallskip

\noindent\textbf{Algorithm.} Compute the set $\Mscr_{\Iw,\Ib,\Jw,\Jb}$ of
feasible matchings $M$ for each $(I|J,I'|J')\in\Iscr$, and similarly for
$\Kscr$. For each instance $M$ occurring there, we extract the family
$\bfC_M(\Iscr)$ of all configurations concerning $M$ in $\bfC(\Iscr)$, and
extract a similar family $\bfC_M(\Kscr)$ in $\bfC(\Kscr)$. If
$|\bfC_M(\Iscr)|\ne |\bfC_M(\Kscr)|$ for at least one instance $M$, then
$\Iscr$ and $\Kscr$ are not balanced at all. Otherwise for each $M$, we seek
for a required bijection $\gamma_M: \bfC_M(\Iscr)\to \bfC_M(\Kscr)$ by solving
the maximum matching problem in the corresponding bipartite graph $H_M$. More
precisely, the vertices of $H_M$ are the tuples $(I|J,I'|J';M)$ and
$(K|L,K'|L';M)$ occurring in $\bfC_M(\Iscr)$ and $\bfC_M(\Kscr)$, and such
tuples are connected by edge in $H_M$ if they obey~\refeq{q_balan}. Find a
maximum matching $N$ in $H_M$. (There are many fast algorithms to solve this
classical problem; for a survey, see, e.g.~\cite{schr}.) If
$|N|=|\bfC_M(\Iscr)|$, then $N$ determines the desired $\gamma_M$ in a natural
way. Taking together, these $\gamma_M$ give a bijection between $\bfC(\Iscr)$
and $\bfC(\Kscr)$ as required, implying that $\Iscr,\Kscr,\alpha,\beta$ are
$q$-balanced And if $|N|<|\bfC_M(\Iscr)|$ for at least one instance $M$, then
the algorithm declares the non-$q$-balancedness.


\section{\Large Examples of universal quadratic identities}  \label{sec:examples}

The flow-matching method described above is well adjusted to prove, relatively
easily, classical or less known quadratic identities. In this section we give a
number of appealing illustrations.

Instead of circular diagrams as in Section~\SEC{double}, we will use more
compact, but equivalent, \emph{two-level diagrams}. Also when dealing with a
flag pair $(I|J)$, i.e., when $I$ consists of the elements $1,2,\ldots, |I|$,
we may use an appropriate \emph{one-level diagrams}, which leads to no loss of
generality. For example, the refined cortege $(\Iw=\{3,4\},\;
\Ib=\emptyset,\;\Jw=\{1',3',4',6'\},\;\Jb=\{2',5'\})$ with the feasible
matching $\{1'2',4'5',33',46'\}$ can be visualized in three possible ways as:

\vspace{0.cm}
\begin{center}
\includegraphics{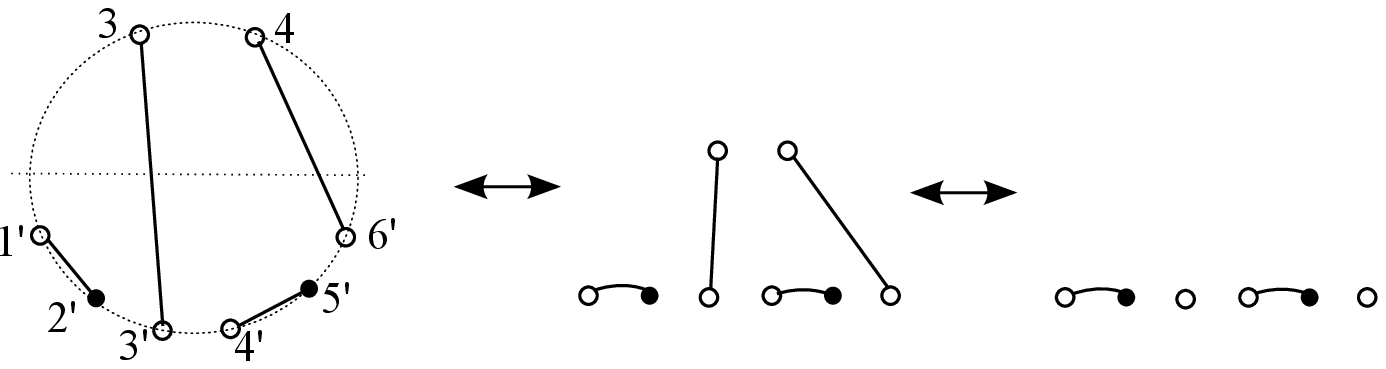}
\end{center}
\vspace{0.cm}

A couple $\{i,j\}$ may be denoted as $ij$. Also for brevity we write $Xi\ldots
j$ for $X\cup\{i,\ldots,j\}$, where $X$ and $\{i,\ldots,j\}$ are disjoint.

As before, we use notation $[I|J]$ for the corresponding $q$-minor of the path
matrix $\Path_G$ (defined in Section~\SEC{flows}). In the flag case $[I|J|$ is
usually abbreviated to $[J]$ (in view of $I=\{1,\ldots,|J|\}$).


\subsection{Commuting minors.}  \label{ssec:com_min}
We start with a simple illustration of our method by showing that $q$-minors
$[I|J]$ and $[I'|J']$ ``purely'' commute when $I'\subset I$ and $J'\subset J$.
(This matches the known fact that a minor of a $q$-matrix commutes with any of
its subminors, or that the $q$-determinant of a square $q$-matrix is a central
element of the corresponding algebra.)

Let $\Iw=I-I'$ consist of $i_1<\ldots<i_k$, and $\Jw=J-J'$ consist of
$j_1<\ldots <j_k$. Since $\Ib=I'-I=\emptyset$ and $\Jb=J'-J=\emptyset$, there
is only one feasible matching $M$ for $(\Iw,\Ib,\Jw,\Jb)$; namely, the one
formed by the $RC$-couples $\pi_\ell=i_\ell j_\ell$, $\ell=1,\ldots,k$. The
index exchange operation applied to $(I|J,I'|J')$ using the whole $M$ produces
the cortege  $(K|L,K'|L')$ for which $\Kw=\Ib=\emptyset$, $\Kb=\Iw$,
$\Lw=\Jb=\emptyset$, $\Lb=\Jw$ (and $K\cap K'=I\cap I'$, $L\cap L'=J\cap J'$).
Since $M$ consists of $RC$-couples only, we have $\zeta^\circ
(\Iw,\Ib,\Jw,\Jb;M)=\zeta^\bullet(\Iw,\Ib,\Jw,\Jb; M)=0$. So the (one-element)
families $\Iscr=\{(I|J,I'|J')\}$ and $\Kscr=\{(K|L,K'|L')\}$ along with
$\alpha=\beta=0$ are $q$-balanced, and Theorem~\ref{tm:suff_q_bal} gives the
desired equality $[I|J][I'|J']=[I'|J'][I|J]$.

This is illustrated in the picture with two-level diagrams (in case $k=5$).
Hereinafter we indicate by crosses the couples that are involved in the index
exchange operation that is applied (i.e., the couples where the colors of
elements are changed).

\vspace{0.cm}
\begin{center}
\includegraphics{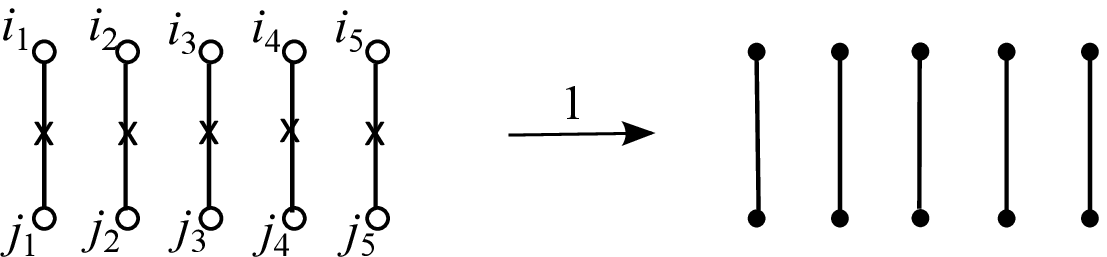}
\end{center}
\vspace{-0.5cm}


\subsection{Quasicommuting minors.}  \label{ssec:quasicommut}
Recall that two sets $I,J\subseteq[n]$ are called \emph{weakly separated} if,
up to renaming $I$ and $J$, there holds: $|I|\ge |J|$, and $J-I$ has a
partition $J_1\cup J_2$ such that $J_1<I-J<J_2$ (where we write $X<Y$ if $x<y$
for any $x\in X$ and $y\in Y$). Leclerc and Zelevinsky proved the following
  \begin{theorem}[\cite{LZ}] \label{tm:LZ}
Two flag minors $[I]$ and $[J]$ of a quantum matrix \emph{quasicommute}, i.e.,
satisfy
  \begin{equation} \label{eq:quasiIJ}
  [I][J] = q^c [J][I]
  \end{equation}
for some $c\in\Zset$, if and only if the column sets $I,J$ are weakly
separated. Moreover, when $|I|\ge |J|$ and $J_1\cup J_2$ is a partition of
$J-I$ with $J_1<I-J<J_2$, the number $c$ in~\refeq{quasiIJ} is equal to
$|J_2|-|J_1|$.
  \end{theorem}
(In case $I\cap J=\emptyset$, ``if'' part is due to Krob and
Leclerc~\cite{KL}). We explain how to obtain ``if'' part of Theorem~\ref{tm:LZ}
by use of the flow-matching method.

Let $A:=\{1,\ldots,|I|\}$, $B:=\{1,\ldots,|J|\}$, and define
  $$
  A^\circ:=A-B,\quad B^\bullet:=B-A\; (=\emptyset),\quad
    \Iw:=I-J,\quad \Jb:=J-I.
  $$
One can see that $(A^\circ,B^\bullet,\Iw,\Jb)$ has exactly one feasible
matching $M$; namely, $J_1$ is coupled with the first $|J_1|$ elements of
$\Iw$, ~$J_2$ is coupled with the last $|J_2|$ elements of $\Iw$ (forming all
$C$-couples), and the rest of $\Iw$ is coupled with $A^\circ$ (forming all
$RC$-couples).

Observe that the index exchange operation applied to $(A|I,B|J)$ using the
whole $M$ swaps $A|I$ and $B|J$ (since it changes the colors of all elements in
$A^\circ$, $\Iw$ and $\Jb$). Also $M$ consists of $|J_1|+|J_2|$ ~$C$-couples
and $|A^\circ|$ ~$RC$-couples. Moreover, the $C$-couples are partitioned into
$|J_1|$ couples $ij$ with $i<j$ and $i\in J_1$, and $|J_2|$ couples $ij$ with
$i<j$ and $j\in J_2$. This gives $\zeta^\circ=|J_2|$ and $\zeta^\bullet=|J_1|$.
Hence the (one-element) families $\{(A|I,B|J)\}$ and $\{(B|J,A|I)\}$ along with
$\alpha(A|I,B|J)=0$ and $\beta(B|J,A|I)=|J_2|-|J_1|$ are $q$-balanced. Now
Theorem~\ref{tm:suff_q_bal} implies~\refeq{quasiIJ} with $c=|J_2|-|J_1|$.

The picture with two-level diagrams illustrates the case $|I-J|=5$, $|J-I|=3$,
$|J_1|=1$ and $|J_2|=2$.

\vspace{0.cm}
\begin{center}
\includegraphics{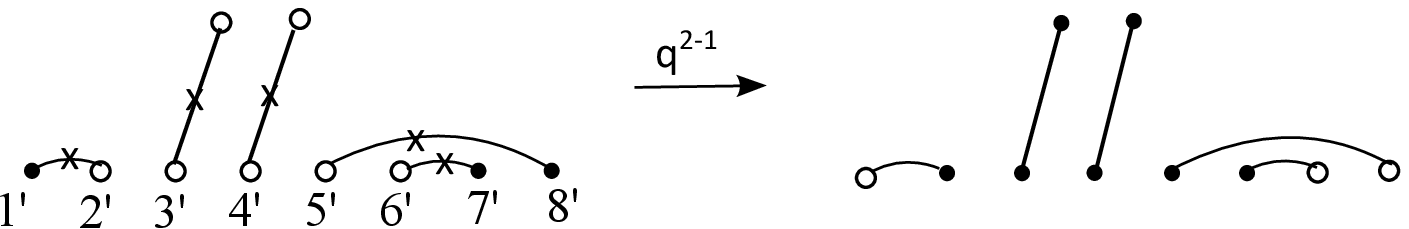}
\end{center}
\vspace{0cm}

``Only if'' part of Theorem~\ref{tm:LZ} will be discussed in
Section~\SEC{concl}. Also we will give there a generalization of this theorem
that characterizes the set of all pairs of quasicommuting $q$-minors (not
necessarily flag ones).


\subsection{Manin's relations in path matrices.}  \label{ssec:abcd}
We prove Theorem~\ref{tm:Path-A}.

(a) Consider entries $[i|j]$ and $[i|j']$ with $j<j'$ in $\Path_G$. The cortege
$S=(i|j,i|j')$ admits a unique feasible matching; it consists of the single
$C$-couple $\pi=jj'$. The index exchange operation using $\pi$ transforms $S$
into $T=(i|j',i|j)$; see the picture with one-level diagrams:

\vspace{0.cm}
\begin{center}
\includegraphics{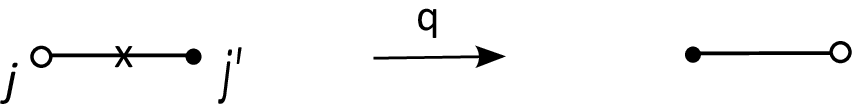}
\end{center}
\vspace{-0.3cm}

We observe that $\{S\}$ and $\{T\}$ along with $\alpha=0$ and $\beta=1$
($=\zeta^\circ-\zeta^\bullet$) are $q$-balanced, and
Theorem~\ref{tm:suff_q_bal} yields $[i|j][i|j']=q[i|j'][i|j]$, as required.

(b) For a $2\times 1$ submatrix of $\Path_G$, the argument is similar.
\smallskip

(c) Consider a $2\times 2$ submatrix $\genfrac{(}{)}{0pt}{1}{c\;d}{a\;b}$ of
$\Path_G$, where $a=[i|j]$, $b=[i|j']$, $c=[i'|j]$, $d=[i'|j']$ (then $i<i'$
and $j<j'$). Let $\Iscr$ consist of two corteges $S_1=(i|j,i'|j')$,
$S_2=(i|j',i'|j)$, and $\Kscr$ consist of two corteges $T_1=(i|j',i'|j)$,
$T_2=(i'|j',i|j)$ (note that $S_2=T_1$). Observe that $S_1$ admits 2 feasible
matchings, namely, $M=\{ii',jj'\}$ and $N=\{ij,i'j'\}$, while $S_2$ admits only
one feasible matching $M$. In their turn, $\Mscr(T_1)=\{M\}$ and
$\Mscr(T_2)=\{M,N\}$. Hence we can form the bijection between $\bfC(\Iscr)$ and
$\bfC(\Kscr)$ that sends $(S_1;M)$ to $(T_1;M)$, ~$(S_1,N)$ to $(T_2;N)$, and
$(S_2,M)$ to $(T_2;M)$. This bijection is illustrated in the picture (where, as
before, we indicate the submathings involved in the exchange operations with
crosses).

\vspace{0.cm}
\begin{center}
\includegraphics{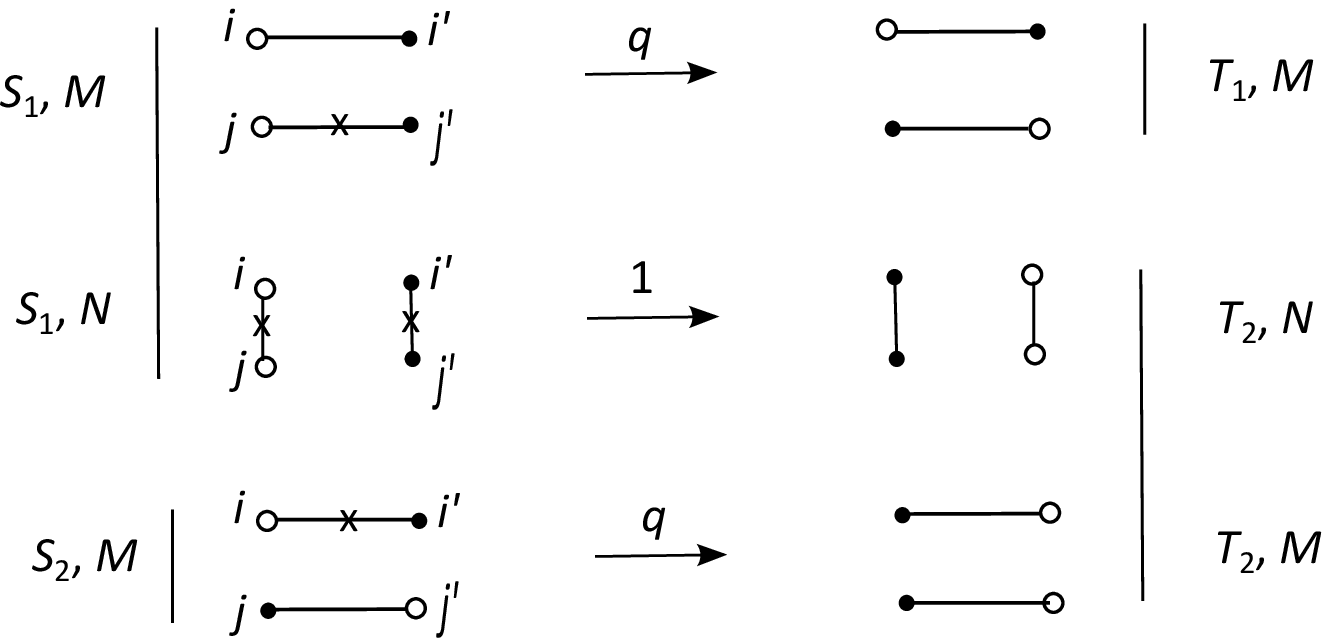}
\end{center}
\vspace{0cm}

Assign $\alpha(S_1)=0$, $\alpha(S_2)=-1$, $\beta(T_1)=1$ and $\beta(T_2)=0$.

One can observe from the above diagrams that $\Iscr,\Kscr,\alpha,\beta$ are
$q$-balanced. We obtain
   $$
 [i|j][i'|j']+q^{-1}[i|j'][i'|j]=q[i|j'][i'|j]+[i'|j'][i|j],
   $$
yielding $ad-da=(q-q^{-1}) bc$, as required.

Finally, to see $bc=cb$, take the 1-element families $\{S'=(i|j',i'|j)\}$ and
$\{T'=(i'|j,i|j')\}$; then $\{ii',jj'\}$ is the only feasible matching for each
of $S',T'$. The above families along with $\alpha=\beta=0$ are $q$-balanced, as
is seen from the picture:

\vspace{0.cm}
\begin{center}
\includegraphics{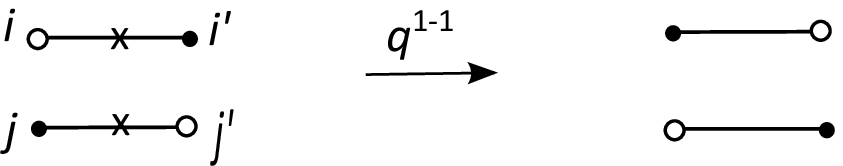}
\end{center}
\vspace{-0.3cm}

This gives $[i|j'][i'|j]=[i'|j][i|j']$, or $bc=cb$, as required.


\subsection{Relations with triples and quadruples.}  \label{ssec:P3P4}
In the commutative case (when dealing with the commutative coordinate ring of
$m\times n$ matrices over a field), the simplest examples of quadratic
identities on flag minors are presented by the classical Pl\"ucker relations
involving 3- and 4-element sets of columns. More precisely, for
$A\subseteq[n]$, let $g(A)$ denote the flag minor with the set $A$ of columns.
Then for any three elements $i<j<k$ in $[n]$ and a set $X\subseteq
[n]-\{i,j,k\}$, there holds
  \begin{equation} \label{eq:P3}
    g(Xik)g(Xj)=g(Xij)g(Xk)+g(Xjk)g(Xi),
    \end{equation}
and for any $i<j<k<\ell$ and $X\subseteq[n]-\{i,j,k,\ell\}$,
  \begin{equation} \label{eq:P4}
    g(Xik)g(Xj\ell)=g(Xij)g(Xk\ell)+g(Xj\ell)g(Xjk),
    \end{equation}

There are two quantized counterparts of~\refeq{P3} (concerning flag $q$-minors
of the matrix $\Path_G$). One of them is viewed as
   \begin{equation} \label{eq:qP3a}
   [Xj][Xik]=[Xij][Xk]+[Xjk][Xi],
   \end{equation}
and the other as
   \begin{equation} \label{eq:qP3b}
   [Xik][Xj]=q^{-1}[Xij][Xk]+q[Xjk][Xi].
   \end{equation}

To see~\refeq{qP3a}, associate to $Xj$ the white pair
$(\Iw,\Jw)=(\emptyset|\{j\})$, and to $Xik$ the black pair
$(\Ib|\Jb)=(\{p\}|\{i,k\})$, where $p$ is the last row index for $[Xik]$ (i.e.,
$p=|X|+2$). Then $\Mscr_{\Iw,\Ib,\Jw,\Jb}$ consists of two feasible matchings:
$M=\{pi,jk\}$ and $N=\{ij,pk\}$. Now~\refeq{qP3a} is seen from the following
picture with two-level diagrams, where we write $S$ for the cortege
$([p-1]\,|Xj, [p]\,|Xik)$, ~$T_1$ for $([p]\,|Xij, [p-1]\,|Xk)$, and $T_2$ for
$([p]\,|Xjk, [p-1]\,|Xi)$:

\vspace{0.cm}
\begin{center}
\includegraphics{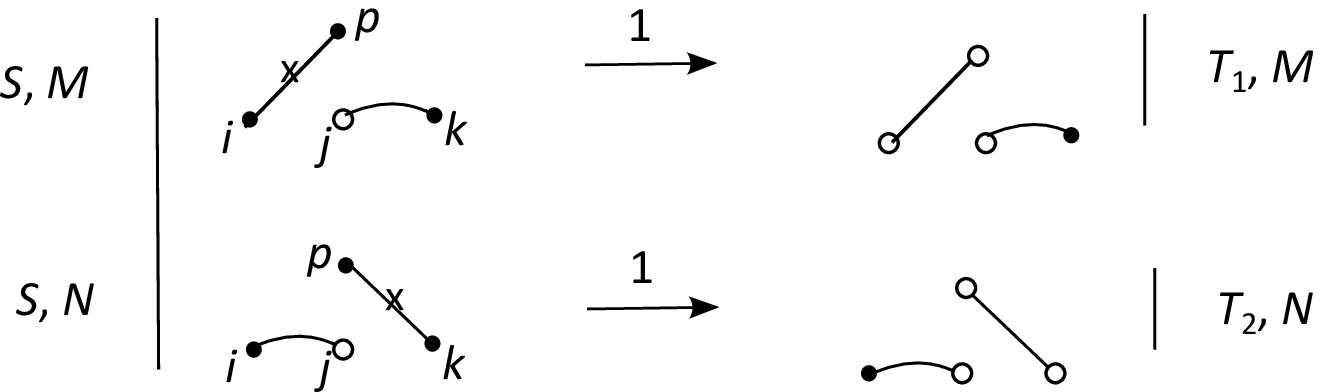}
\end{center}
\vspace{-0.3cm}

As to~\refeq{qP3b}, it suffices to consider one-level diagrams (as we are not
going to use $RC$-couples in the exchange operations). Now the ``white'' object
is the column set $\Jw=\{i,k\}$ and the ``black'' object is $\Jb=\{j\}$. Then
$\Mscr_{\{p\},\emptyset,\Jw,\Jb}$ consists of two feasible matchings, one using
the $C$-couple $\pi=jk$, and the other using the $C$-couple $\mu=ij$.
Now~\refeq{qP3b} can be seen from the picture, where we write $S$ for the flag
cortege $(Xik,Xj)$, ~$T_1$ for $(Xij,Xk)$, and $T_2$ for $(Xjk,Xi)$.

\vspace{-0.3cm}
\begin{center}
\includegraphics{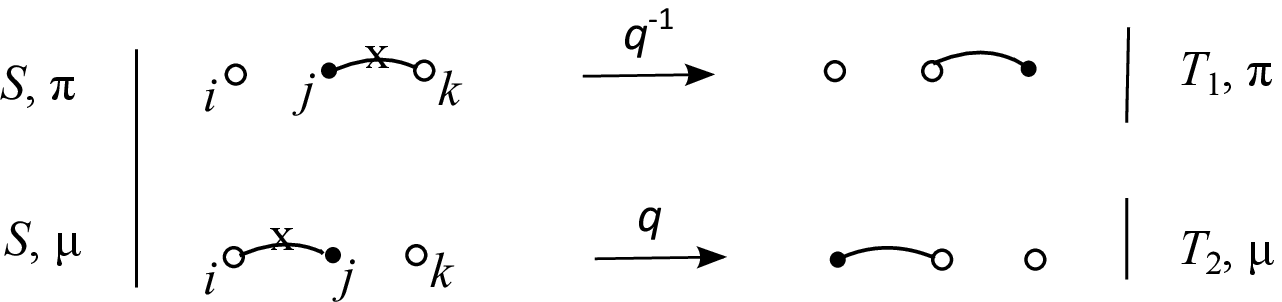}
\end{center}
\vspace{-0.3cm}

Next we demonstrate the following quantized counterpart of~\refeq{P4}:
   \begin{equation} \label{eq:qP4}
   [Xik][Xj\ell]=q^{-1}[Xij][Xk\ell]+q[Xi\ell][Xjk].
   \end{equation}
To see this, we use one-level diagrams and consider the column sets
$\Jw=\{i,k\}$ and $\Jb=\{j,\ell\}$. Then $\Mscr_{\emptyset,\emptyset,\Jw,\Jb}$
consists of two feasible matchings: $M=\{i\ell,jk\}$ and $N=\{ij,k\ell\}$.
Identity~\refeq{qP4} can be seen from the picture, where $S=(Xik,Xj\ell)$,
$T_1=(Xij,Xk\ell)$ and $T_2=(Xi\ell,Xjk)$.

\vspace{-0cm}
\begin{center}
\includegraphics{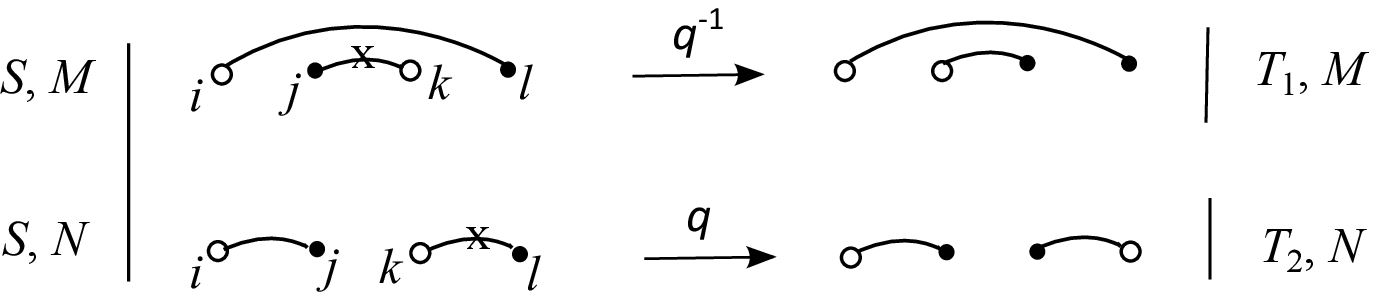}
\end{center}
\vspace{-0.3cm}

\medskip
\noindent\underline{\em Remark 4.} Note that, if wished, one can produce more
identities from~\refeq{qP3a} and~\refeq{qP3b}, using the fact that $Xij$ and
$Xk$ (as well as $Xjk$ and $Xi$) are weakly separated, and therefore their
corresponding flag $q$-minors quasicommute (see Section~\SSEC{quasicommut}). In
contrast, $Xj$ and $Xik$ are not weakly separated. Next, subtracting
from~\refeq{qP3b} identity~\refeq{qP3a} multiplied by $q$ results in the
identity of the form
   $$
   [Xik][Xj]=q[Xj][Xik]-(q-q^{-1})[Xij][Xk],
   $$
which is in spirit of \emph{commutation relations} for quantum minors studied
in~\cite{fio,good}.


\subsection{Dodgson's type identity.}  \label{ssec:dodg}
As one more simple illustration of our method, we consider a $q$-analogue of
the classical Dodgson's condensation formula for usual minors~\cite{dod}. It
can be stated as follows: for elements $i<k$ of $[m]$, a set $X\subseteq
[m]-\{i,k\}$, elements $i'<k'$ of $[n]$, and a set $X'\subseteq [n]-\{i',k'\}$
(with $|X'|=|X|$),
  \begin{equation} \label{eq:dodg}
[Xi|X'i'][Xk|X'k']= q[Xi|X'k'][Xk|X'i']+[Xik|X'i'k'][X|X'].
   \end{equation}
In this case we deal with the cortege $S=(I|J,I'|J')=(Xi|X'i',Xk|X'k')$ and its
refinement $(\Iw,\Ib,\Jw,\Jb)$ of the form $(i,k,i',k')$. The latter admits two
feasible matchings: $M=\{ik,i'k'\}$ and $N=\{ii',kk'\}$. Now~\refeq{dodg} can
be concluded by examining the picture below, where $T_1$ stands for
$(Xi|X'k',Xk|X'i')$, and $T_2$ for $(Xik|X'i'k', X|X')$:

\vspace{-0cm}
\begin{center}
\includegraphics{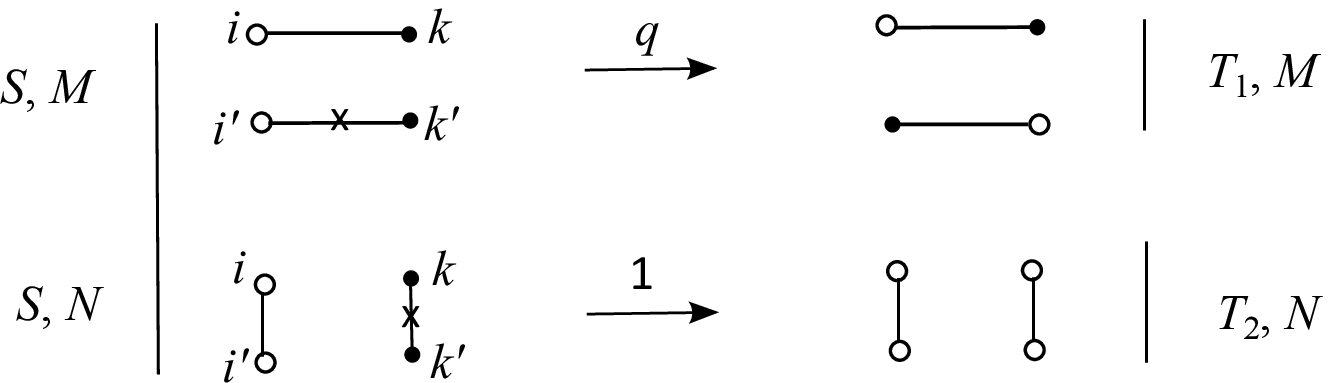}
\end{center}
\vspace{-0.3cm}


\subsection{Two general quadratic identities.}  \label{ssec:gen2}
Two representable quadratic identities of a general form were established for
quantum flag minors in~\cite{LR,TT}.

The first one considers column subsets $I,J\subset[n]$ with $|I|\le|J|$ and is
viewed as
   \begin{equation} \label{eq:R1}
[I][J]= \sum_{\mu\subseteq J-I,\, |\mu|=|J|-|I|}
(-q)^{Inv(J-\mu,\,\mu)-Inv(I,\,\mu)}
                      [I\cup\mu][J-\mu],
   \end{equation}
where $Inv(A,B)$ denotes the number of pairs $(a,b)\in A\times B$ with $a>b$.
Observe that~\refeq{qP3a} is a special case of~\refeq{R1} in which the roles of
$I$ and $J$ are played by $Xj$ and $Xik$, respectively. Indeed, in this case
$\mu$ ranges over the singletons $\{i\}$ and $\{k\}$, and we have
$Inv(Xk,i)-Inv(Xj,i)=0$ and $Inv(Xi,k)-Inv(Xj,k)=0$. (For brevity, we write
$Inv(\cdot,i')$ for $Inv(\cdot,\{i'\})$.)

The second one considers $I,J\subset[n]$ with $|I|-|J|\ge 2$ and is viewed as
   \begin{equation} \label{eq:R2}
\sum\nolimits_{a\in I-J} (-q)^{Inv(a,I-a)-Inv(a,J)} [Ja][I-a] =0
   \end{equation}
(where we write $Ja$ for $J\cup\{a\}$, and $I-a$ for $I-\{a\}$). A special case
is~\refeq{qP4} (with $I=Xjk\ell$ and $J=Xi$).

We explain how~\refeq{R1} and~\refeq{R2} can be proved for flag $q$-minors of
$\Path_G$ by use of our flow-matching method.
\medskip

\noindent\textbf{Proof of~\refeq{R1}.} ~The pair $(I,J)$ corresponds to the
cortege $S:=([p]\,|I,\;[p+k]\,|J)$ and its refinement
$R:=(\emptyset,\;Q:=\{p+1,\ldots, p+k\},\; \Iw:=I-J,\; \Jb:=J-I)$, where
$p:=|I|$ and $k:=|J|-|I|$. In its turn, each pair $(I\cup\mu)|(J-\mu)$
occurring in the R.H.S. of~\refeq{R1} corresponds to the cortege
$S_\mu:=([p+k]\,|(I_\mu:=I\cup\mu),\; [p]\,|(J_\mu:=J-\mu))$ and its refinement
$R_\mu:=(Q,\;\emptyset,\; I^\circ_\mu:=\Iw\cup\mu,\; J^\bullet_\mu:=\Jb-\mu)$.

So we deal with the set
  $$
  \Fscr:=\{S\}\cup\{S_\mu\colon \mu\subset \Jb, \,|\mu|=k\},
  $$
of corteges and the related set $\bfC(\Fscr)$ of configurations (of the form
$(S;M)$ or $(S_\mu;M)$), and our aim is to construct an involution
$\gamma:\bfC(\Fscr)\to \bfC(\Fscr)$ which is agreeable with matchings, signs
and $q$-factors figured in~\refeq{R1}. (Under reducing~\refeq{R1} to the
canonical form, $\Fscr$ splits into two families $\Iscr$ and $\Kscr$, and
$\gamma$ determines the $q$-balancedness for $\Iscr,\Kscr$ with corresponding
$\alpha,\beta$.)

Consider a refined cortege $R_\mu=(Q,\emptyset,I^\circ_\mu,J^\bullet_\mu)$ and
a feasible matching $M$ for it. Note that $M$ consists of $k=|Q|$ ~$RC$-couples
(connecting $Q$ and $I^\circ_\mu$) and $|J^\bullet_\mu|=|\Iw|$ ~$C$-couples
(connecting $I^\circ_\mu$ and $J^\bullet_\mu$). Two cases are possible.
  \smallskip

\noindent\underline{\emph{Case 1}}: each $C$-couple connects $J^\bullet_\mu$
and $\Iw$. Then all $RC$-couples in $M$ connect $Q$ and $\mu$. Therefore, the
exchange operation applied to $S_\mu$ using the set $\varPi$ of all
$RC$-couples of $M$ produces the ``initial'' cortege $S$ (corresponding to the
refinement $R=(\emptyset,Q, \Iw,\Jb)$). Clearly $M$ is a feasible matching for
$S$ and the exchange operation applied to $S$ using $\varPi$ returns $S_\mu$.
We link $(S;M)$ and $(S_\mu;M)$ by $\gamma$.

Note that for each $C$-couple $\pi=ij\in M-\varPi$ and for each $r\in \mu$,
either $r<i,j$ or $r>i,j$ (otherwise the $RC$-couple containing $r$ would
``cross'' $\pi$, contrary to the planarity requirement~\refeq{feasM}(ii) for
$M$). This implies $Inv(J_\mu,\mu)=Inv(I_\mu,\mu)$, whence the terms $[I][J]$
in the L.H.S. and $(-q)^0[I_\mu][J_\mu]$ in the R.H.S. of~\refeq{R1} are
$q$-balanced.
 \medskip

\noindent\underline{\emph{Case 2}}: there is a $C$-couple in $M$ connecting
$J^\bullet_\mu$ and $\mu$. Among such couples, choose the couple $\pi=ij$ with
$i<j$ such that: (a) $j-i$ is minimum, and (b) $i$ is minimum subject to~(a).
From~\refeq{feasM} and~(a) it follows that
  \begin{numitem1} \label{eq:betw_ij}
if a couple $\pi'\in M$ has an element (strictly) between $i$ and $j$, then
$\pi'$ connects $\Iw$ and $J^\bullet_\mu$, and the other element of $\pi'$ is
between $i$ and $j$ as well.
  \end{numitem1}

Let $S_{\mu'}$ be obtained by applying to $S_\mu$ the exchange operation using
the single couple $\pi$. Then $\mu'=\mu\triangle\pi$,
$I^\circ_{\mu'}=I^\circ_\mu\triangle\pi$ and
$J^\bullet_{\mu'}=J^\bullet_\mu\triangle\pi$. The matching $M$ is feasible for
$S_{\mu'}$, we are in Case~2 with $S_{\mu'}$ and $M$, and one can see that the
couple $\pi'\in M$ chosen for $S_{\mu'}$ according to the above rules~(a),(b)
coincides with $\pi$. Based on these facts, we link $(S_\mu; M)$ and
$(S_{\mu'}; M)$ by $\gamma$.

Now we compute and compare the numbers $a:=Inv(J^\bullet_{\mu'}=J-\mu',\,\mu')
- Inv(\Jb_{\mu}=J-\mu,\,\mu)$ and $b:=Inv(I,\mu')-Inv(I,\mu)$. Let $d$ be the
number of elements of $\Iw$ between $i$ and $j$ (recall that $\pi=ij$ and
$i<j$). Property~\refeq{betw_ij} ensures that the number of elements of
$J^\bullet_\mu$ (as well as of $J^\bullet_{\mu'}$) between $i$ and $j$ is equal
to $d$ too. Consider two possibilities.
  \smallskip

\noindent\underline{\emph{Subcase 2a}}: $i\in\mu$ (and $j\in J^\bullet_\mu$).
Then $i\in J^\bullet_{\mu'}$ and $j\in \mu'$. This implies that
$a=Inv(J^\bullet_{\mu'},j)- Inv(J^\bullet_\mu,i)=d+1$ and
$b=Inv(\Iw,j)-Inv(\Iw,i)=d$.
 \smallskip

\noindent\underline{\emph{Subcase 2b}}: $i\in J^\bullet_\mu$ (and $j\in \mu$).
Then $i\in\mu'$ and $j\in J^\bullet_{\mu'}$, yielding $a=-d-1$ and $b=-d$.
  \smallskip

Finally, let $(-q)^\alpha$ and $(-q)^\beta$ be the multipliers to the terms
$[I_\mu][J_\mu]$ and $[I_{\mu'}][J_{\mu'}]$ in~\refeq{R1}, respectively. Then
$\beta-\alpha=a-b$, which is equal to 1 in Subcase~2a and $-1$ in Subcase~2b.
In both cases this amounts to the value $\zeta^\circ-\zeta^\bullet$ for the
exchange operation applied to $S_\mu$ using $\pi$, and validity of~\refeq{R1}
follows from Theorem~\ref{tm:suff_q_bal}. \hfill\qed
\medskip

\noindent\underline{\em Remark 5.} Sometimes it is useful to consider the
identity formed by the corteges reversed to those in~\refeq{R1}; by
Proposition~\ref{pr:reversed}, it is viewed as
   $$
[J][I]= \sum_{\mu\subseteq J-I,\, |\mu|=|J|-|I|}
(-q)^{Inv(I,\,\mu)-Inv(J-\mu,\,\mu)} [J-\mu][I\cup\mu].
  $$

\noindent\textbf{Proof of~\refeq{R2}.} Let $p:=|J|$,\; $k:=|I|-|J|$,\;
$Q:=[p+k-1]-[p+1]$,\; $\Jw:=J-I$ and $\Ib:=I-J$. For $a\in \Ib$, the pair
$(Ja,I-a)$ in~\refeq{R2} corresponds to the cortege $S_a:=([p+1]\,|Ja,\;
[p+k-1]\,|(I-a))$ and its refinement $R_a:=(\emptyset,Q,\,J^\circ a,\,
I^\bullet_a:=\Ib-a)$ (we use the fact that $k\ge 2$).

We deal with the set $\Fscr:=\{S_a\colon a\in\Ib\}$ of corteges and the related
set $\bfC(\Fscr)$ of configurations $(S_a;M)$, and like the previous proof, our
aim is to construct an appropriate involution $\gamma:
\bfC(\Fscr)\to\bfC(\Fscr)$.

Consider a refined cortege $R_a=(\emptyset, Q,\Jw a,I^\bullet_a)$ and a
feasible matching $M$ for it. Take the couple in $M$ containing $a$, say,
$\pi=\{a,b\}$. Note that $\pi$ is a $C$-couple and $b\in I^\bullet_a$ (since
$a$ is white, and $Q$ and $I^\bullet_a$ are black). The exchange operation
applied to $S_a$ using $\pi$ produces the member $S_b$ of $\Fscr$, and we link
$S_a$ and $S_b$ by $\gamma$.

It remains to estimate the multipliers $(-q)^\alpha$ and $(-q)^\beta$ to the
terms $[Ja][I-a]$ and $[Jb][I-b]$ in~\refeq{R2}, respectively.

Let $d$ be the number of elements of $\Ib$ between $a$ and $b$. It is equal to
the number of elements of $\Jw$ between $a$ and $b$ (since, in view
of~\refeq{feasM}, the elements of $\Ib\cup\Jw$ between $a$ and $b$ must be
partitioned into $C$-couples in $M$). This implies that if $a<b$, then
$Inv(b,I-b)-Inv(a,I-a)=d+1$ and $Inv(b,J)-Inv(a,J)=d$. Therefore,
$\beta-\alpha=(d+1)-d=1$. And if $a>b$, then $Inv(b,I-b)-Inv(a,I-a)=-d-1$ and
$Inv(b,J)-Inv(a,J)=-d$, whence $\beta-\alpha=-1$. In both cases, $\beta-\alpha$
coincides with the corresponding value of $\zeta^\circ-\zeta^\bullet$, and the
result follows. \hfill\qed


\section{\Large Necessity of the $q$-balancedness}  \label{sec:necess}

In this section we show a converse assertion to Theorem~\ref{tm:suff_q_bal},
thus obtaining a complete characterization for the UQ identities on quantized
minors. This characterization, given in terms of the $q$-balancedness,
justifies the algorithm of recognizing UQ identities described in the end of
Section~\SEC{q_relat}. As before, we deal with homogeneous families of corteges
in $\Escr^{m,n}\times\Escr^{m,n}$.

\begin{theorem} \label{tm:nec_q_bal}
Let $\Kset$ be a field of characteristic zero and let $q\in\Kset^\ast$ be
transcendental over $\Qset$. Suppose that $\Iscr,\Kscr,\alpha,\beta$ (as in
Section~\SEC{q_relat}) are not $q$-balanced. Then there exists (and can be
explicitly constructed) an SE-graph $G$ for which relation~\refeq{gen_QR} is
violated.
 \end{theorem}
   \begin{proof}
~We essentially use an idea and construction worked out for the commutative
version in~\cite[Sec.~5]{DKK}.

Recall that the homogeneity of $\Fscr:=\Iscr\sqcup\Kscr$ means the existence of
$\Xr,\Yr\subseteq[m]$ and $\Xc,\Yc\subseteq[n]$ such that any cortege $(I|J,
I'|J')\in\Fscr$ satisfies
  \begin{equation} \label{eq:XrYr}
  I\cap I'=\Xr,\quad I\triangle I'=\Yr,\quad J\cap J'=\Xc,
          \quad J\triangle J'=\Yc
  \end{equation}
(cf.~\refeq{homogen}). For a perfect matching $M$ on $\Yr\sqcup\Yc$, let us
denote by $\Iscr_M$ the set of corteges $S=(I|J,I'|J')\in\Iscr$ for which $M$
is feasible (see~\refeq{feasM}), and denote by $\Kscr_M$ a similar set for
$\Kscr$. The $q$-balancedness of $\Iscr,\Kscr,\alpha,\beta$ would mean that,
for any $M\in\bfM(\Fscr)$, there exists a bijection
$\gamma_M:\Iscr_M\to\Kscr_M$ respecting~\refeq{q_balan}. That is, for any
$S=(I|J,I'|J')\in\Iscr_M$ and for $T=(K|L,K'|L')=\gamma_M(S)$, there holds
   \begin{equation} \label{eq:baz}
   \beta(T)-\alpha(S)=\zeta^\circ(\varPi_{S,T})-\zeta^\bullet(\varPi_{S,T}).
   \end{equation}
Here: $\varPi=\varPi_{S,T}$ is the subset of $M$ such that the refined cortege
$(\Kw,\Kb,\Lw,\Lb)$ is obtained from $(\Iw,\Ib,\Jw,\Jb)$ by the index exchange
operation using $\varPi$, and $\zeta^\circ(\varPi)$ (resp.
$\zeta^\bullet(\varPi)$) is the number of $R$- and $C$-couples
$\{i,j\}\in\varPi$ with $i<j$ and $i\in \Iw\cup\Jw$ (resp. $i\in\Ib\cup\Jb$).
The following assertion is crucial.

  \begin{prop} \label{pr:P1P2}
Let $M$ be a perfect planar matching on $\Yr\sqcup\Yc$. Then there exists (and
can be explicitly constructed) an SE-graph $G=(V,E)$ with the following
properties: for each cortege $S=(I|J,I'|J')\in\Escr^{m,n}\times \Escr^{m,n}$
satisfying~\refeq{XrYr},
  \begin{itemize}
\item[{\rm(P1)}] if $M$ is feasible for $S$, then $G$ has a unique $(I|J)$-flow and
a unique $(I'|J')$-flow;
\item[{\rm(P2)}] if $M$ is not feasible for $S$, then at least one of
$\Phi_G(I|J)$ and $\Phi_G(I'|J')$ is empty.
 \end{itemize}
  \end{prop}

We will prove this proposition later, and now, assuming that it is valid, we
complete the proof of the theorem.

Let $\Iscr,\Kscr,\alpha,\beta$ be not $q$-balanced. Then there exists a
matching $M\in\bfM(\Fscr)$ that admits no bijection $\gamma_M$ as above between
$\Iscr_M$ and $\Kscr_M$ (in particular, at least one of $\Iscr_M$ and $\Kscr_M$
is nonempty). We fix one $M$ of this sort and consider a graph $G$ as in
Proposition~\ref{pr:P1P2} for this $M$.

Our aim is to show that relation~\refeq{gen_QR} is violated for $q$-minors of
$\Path_{G}$ (yielding the theorem). Suppose, for a contradiction,
that~\refeq{gen_QR} is valid. By~(P2) in the proposition, we have
$[I|J][I'|J']=0$ for each cortege $(I|J,I'|J')\in \Fscr-\Fscr_M$, denoting
$\Fscr_M:=\Iscr_M\sqcup\Kscr_M$. On the other hand, (P1) implies that if
$(I|J,I'|J')\in \Fscr_M$, then
   $$
   [I|J][I'|J']=w(\phi_{I|J})\,w(\phi_{I'|J'}),
   $$
where $\phi_{I|J}$ (resp. $\phi_{I'|J'}$) is the unique $(I|J)$-flow (resp.
$(I'|J')$-flow) in $G$. Thus,~\refeq{gen_QR} can be rewritten as
  \begin{equation} \label{eq:QR_M}
  \sum\nolimits_{\Iscr_M} q^{\alpha(I|J,I'|J')} w(\phi_{I|J}) w(\phi_{I'|J'}) =
     \sum\nolimits_{\Kscr_M} q^{\beta(K|L,K'|L')} w(\phi_{K|L})w(\phi_{K'|L'}).
     \end{equation}

For each cortege $S=(I|J,I'|J')\in\Fscr_M$, the weight $Q(S):= w(\phi_{I|J})\,
w(\phi_{I'|J'})$ of the double flow $(\phi_{I|J},\phi_{I'|J'})$ is a monomial
in weights $w(e)$ of edges $e\in E$ (or a Laurent monomial in inner vertices of
$G$); cf.~\refeq{wP},\refeq{edge_weight},\refeq{w_phi}. For any two corteges in
$\Fscr_M$, one is obtained from the other by the index exchange operation using
a submatching of $M$, and we know from the description in Section~\SEC{double}
that if one double flow is obtained from another by the flow exchange
operation, then the (multi)sets of edges occurring in these double flows are
the same (cf. Lemma~\ref{lm:phi-psi}).

Thus, the (multi)set of edges occurring in the weight monomial $Q(S)$ is the
same for all corteges $S$ in $\Fscr_M$. Fix an arbitrary linear order $\xi$ on
$E$. Then the monomial $Q_\xi=Q_\xi(S)$ obtained from $Q(S)$ by a permutation
of the entries so as to make them weakly decreasing w.r.t. $\xi$ from left to
right is the same for all $S\in\Fscr_M$. Therefore, applying relations
(G1)--(G3) on vertices of $G$ (in Sect.~\SSEC{SE}), we observe that for
$S\in\Fscr_M$, the weight $Q(S)$ is expressed as
   \begin{equation} \label{eq:QS-Qsigma}
   Q(S)=q^{\rho(S)} Q_\xi
   \end{equation}
for some $\rho(S)\in\Zset$. Using such expressions, we rewrite~\refeq{QR_M} as
   $$
  \sum\nolimits_{S\in\Iscr_M} q^{\alpha(S)+\rho(S)} Q_\xi=
     \sum\nolimits_{T\in \Kscr_M} q^{\beta(T)+\rho(T)} Q_\xi,
   $$
obtaining
  \begin{equation} \label{eq:q_rho}
  \sum\nolimits_{S\in\Iscr_M} q^{\alpha(S)+\rho(S)}=
     \sum\nolimits_{T\in\Kscr_M} q^{\beta(T)+\rho(T)}.
  \end{equation}

Since $q$ is transcendental, the polynomials in $q$ in both sides
of~\refeq{q_rho} are equal. Then $|\Iscr_M|=|\Kscr_M|$ and there exists a
bijection $\tilde\gamma:\Iscr_M\to\Kscr_M$ such that
  \begin{equation} \label{eq:alpha_beta_rho}
  \alpha(S)+\rho(S)=\beta(\tilde\gamma(S))+\rho(\tilde\gamma(S)) \qquad
  \mbox{for each $S\in\Iscr_M$}.
  \end{equation}
This together with relations of the form~\refeq{QS-Qsigma} gives
  $$
q^{\alpha(S)}Q(S)=q^{\beta(\tilde\gamma(S))}Q(\tilde\gamma(S)).
  $$
Now, for $S=(I|J,I'|J')\in\Iscr_M$, let $T=(K|L,K'|L'):=\tilde\gamma(S)$ and
let $\varPi:=\varPi_{S,T}$. Using relation~\refeq{gen_exch} from
Corollary~\ref{cor:gen_exch}, we have
  \begin{multline*}
q^{\beta(T)-\alpha(S)}Q(T)=Q(S)=w(\phi_{I|J})\, w(\phi_{I'|J'}) \\
=q^{\zeta^\circ(\varPi)-\zeta^\bullet(\varPi)}w(\phi_{K|L})\,w(\phi_{K'|L'})=
q^{\zeta^\circ(\varPi)-\zeta^\bullet(\varPi)} Q(T),
  \end{multline*}
whence $\beta(T)-\alpha(S)=\zeta^\circ(\varPi)-\zeta^\bullet(\varPi)$. Thus,
the bijection $\gamma_M:=\tilde\gamma$ satisfies~\refeq{baz}. A contradiction.
    \end{proof}

\noindent\textbf{Proof of Proposition~\ref{pr:P1P2}.} ~We utilize the
construction of a graph (which need not be an SE-graph) with properties~(P1)
and~(P2) from~\cite{DKK}; denote this graph by $H=(Z,U)$. We first outline
essential details of that construction and then explain how to turn $H$ into an
equivalent SE-graph $G$. A series of transformations of $H$ that we apply to
obtain $G$ consists of subdividing some edges $e=(u,v)$ (i.e., replacing $e$ by
a directed path from $u$ to $v$) and parallel shifting some sets of vertices
and edges in the plane (preserving the planar structure of the graph). Such
transformations maintain properties~(P1) and~(P2), whence the result will
follow.

Let $\Yr\cup \Xr=\{1,2,\ldots,k\}$ and $\Yc\cup \Xc=\{1',2',\ldots,k'\}$.
Denote the sets of $R$-, $C$-, and $RC$-couples in $M$ by $\Mr,\, \Mc$, and
$\Mrc$, respectively. An $R$-couple $\pi=\{i,j\}$ with $i<j$ is denoted by
$ij$, and we denote by $\prec$ the natural partial order on $R$-couples where
$\pi'\prec\pi$ if $\pi'=pr$ is an $R$-couple with $i<p<r<j$. And similarly for
$C$-couples. When $\pi'\prec \pi$ and there is no $\pi''$ between $\pi$ and
$\pi'$ (i.e., $\pi'\prec\pi''\prec\pi$), we say that $\pi'$ is an immediate
successor of $\pi$ and denote the set of these by $\Succ(\pi)$. Also for
$\pi=ij\in \Mr$ and $d\in\Xr$, we say that $d$ is \emph{open} for $\pi$ if
$i<d<j$ and there is no $\pi'=pr\prec\pi$ with $p<d<r$, and denote the set of
these by $\Open(\pi)$. And similarly for couples in $\Mc$ and elements of
$\Xc$.

A current graph and its ingredients are identified with their images in the
plane, and any edge in it is represented by a (directed) straight-line segment.
We write $(x_v,y_v)$ for the coordinates of a point $v$, and say that an edge
$e=(u,v)$ \emph{points down} if $y_u>y_v$.

The initial graph $H$ has the following features (seen from the construction
in~\cite{DKK}).
\smallskip

(i) The ``sources'' $1,\ldots,k$ (``sinks'' $1',\ldots,k'$) are disposed in
this order from left to right in the upper (resp. lower) half of a
circumference $O$, and the graph $H$ is drawn within the circle (disk) $O^\ast$
surrounded by $O$. (Strictly speaking, the construction of $H$ in~\cite{DKK} is
a mirror reflection of that we describe; the latter is more convenient for us
and does not affect the result.)
\smallskip

(ii) Each couple $\pi=ij\in\Mr\cup\Mc$ is extended to a chord between the
points $i$ and $j$, which is subdivided into a path $L_\pi$ whose edges are
alternately forward and backward. Let $R_\pi$ denote the region in $O^\ast$
bounded by $L_\pi$ and the paths $L_{\pi'}$ for $\pi'\in\Succ(\pi)$. Then each
edge $e$ of $H$ (regarded as a line-segment) having a point in the interior of
$R_\pi$ connects a vertex in $L_\pi$ with either a vertex in $L_{\pi'}$ for
some $\pi'\in\Succ(\pi)$ or some vertex $d\in\Open(\pi)$. Moreover, $e$ is
directed to $L_\pi$ if $\pi\in\Mr$, and from $L_\pi$ if $\pi\in\Mc$.
\smallskip

(iii) Let $R^\ast$ be the region in $O^\ast$ bounded by the paths $L_\pi$ for
all maximal $R$- and $C$-couples $\pi$. Then any edge $e$ of $H$ having a point
in the interior of $R^\ast$ points down. Also if such an $e$ has an incident
vertex $v$ lying on $L_\pi$ for a maximal $R$-couple (resp. $C$-couple) $\pi$,
then $e$ leaves (resp. enters) $v$.
 \smallskip

Using these properties, we transform $H$, step by step, keeping notation
$H=(Z,U)$ for a current graph, and $O^\ast$ for a current region (which becomes
a deformed circle) containing $H$. Iteratively applied steps~(S1) and (S2) are
intended to obtain a graph whose all edges point down. \smallskip

(S1) Choose $\pi=ij\in\Mr$ and let $\bar R_\pi$ be the ``upper part'' of
$O^\ast$ bounded by $L_\pi$ (then $\bar R_\pi$ contains $L_\pi$, the paths
$L_{\pi'}$ for all $\pi'\prec\pi$, and the elements $d\in\Xr$ with $i<d<j$). We
shift $\bar R_\pi-L_\pi$ upward by a sufficiently large distance $\lambda>0$.
More precisely, each vertex $v\in Z$ lying in $\bar R_\pi-L_\pi$ is replaced by
vertex $v'$ with $x_{v'}=x_v$ and $y_{v'}=y_v+\lambda$. Each edge $(u,w)\in U$
of the old graph induces the corresponding edge of the new one, namely: edge
$(u',w')$ if both $u,w$ lie in $\bar R_\pi-L_\pi$; edge $(u,w)$ if $u,w\not\in
\bar R_\pi-L_\pi$; and edge $(u',w)$ if $u\in\bar R_\pi-L_\pi$ and $w\in
L_\pi$. (Case $u\in O^\ast-\bar R_\pi$ and $w\in\bar R_\pi$ is impossible.)
Accordingly, the region $O^\ast$ is enlarged by shifting the part $\bar R_\pi$
by $(0,\lambda)$ and filling the gap between $L_\pi$ and $L_\pi+ (0,\lambda)$
by the corresponding parallelogram.
\smallskip

One can realize that upon application of (S1) to all $R$-couples, the following
property is ensured: for each $\pi\in\Mr$, all initial edges incident to
exactly one vertex on $L_\pi$ turn into edges pointing down. Moreover, since
$L_\pi$ is alternating and there is enough space (from below and from above) in
a neighborhood of the current $L_\pi$, we can deform $L_\pi$ into a zigzag path
with all edges pointing down (by shifting each inner vertex $v$ of $L_\pi$ by a
vector $(0,\eps)$ with an appropriate (positive or negative) $\eps\in\Rset$).
\smallskip

(S2) We choose $\pi\in\Mc$ and act similarly to (S1) with the differences that
now $\bar R_\pi$ denotes the ``lower part'' of $O^\ast$ bounded by $L_\pi$ and
that $\bar R_\pi$ is shifted downward (by a sufficiently large $\lambda>0$).
\smallskip

Upon termination of the process for all $R$- and $C$-couples, all edges of the
current graph $H$ (which is homeomorphic to the initial one) point down, as
required. Moreover, $H$ has one more useful property: the sources $1,\ldots,k$
are ``seen from above'' and the sinks $1',\dots,k'$ are ``seen from below''.
Hence we can add to $H$ ``long'' vertical edges $h_1,\ldots,h_k$ entering the
vertices $1,\ldots, k$, respectively, and ``long'' vertical edges
$h_{1'},\ldots,h_{k'}$ leaving the vertices $1',\ldots,k'$, respectively,
maintaining the planarity of the graph. In the new graph one should transfer
each source $i$ into the tail of $h_i$, and each sink $i'$ into the head of
$h_{i'}$. One may assume that the new sources (sinks) lie within one horizontal
line $L$ (resp. $L'$), and that the rest of the graph lies between $L$ and
$L'$.

Now we get rid of the edges $(u,v)$ such that $x_u>x_v$ (i.e. ``pointing to the
left''), by making the linear transformation $v\mapsto v'$ for the points $v$
in $H$, defined by $x_{v'}=x_v-\lambda y_v$ and $y_{v'}=y_v$ with a
sufficiently large $\lambda>0$.

Thus, we eventually obtain a graph $H$ (homeomorphic to the initial one)
without edges pointing up or to the left. Also the sources and sinks are
properly ordered from left to right in the horizontal lines $L$ and $L'$,
respectively. Now it is routine to turn $H$ into an SE-graph $G$ as required in
the proposition. \hfill\qed
 \medskip

The transformation of $H$ into $G$ as in the proof is illustrated in the
picture; here $\Xr=\{4\}$, ~$\Yr=\{1,2,3\}$, ~$\Xc=\emptyset$,
~$\Yc=\{1',\ldots,5'\}$, and $M=\{12,1'4',2'3',35'\}$.

\vspace{0cm}
\begin{center}
\includegraphics{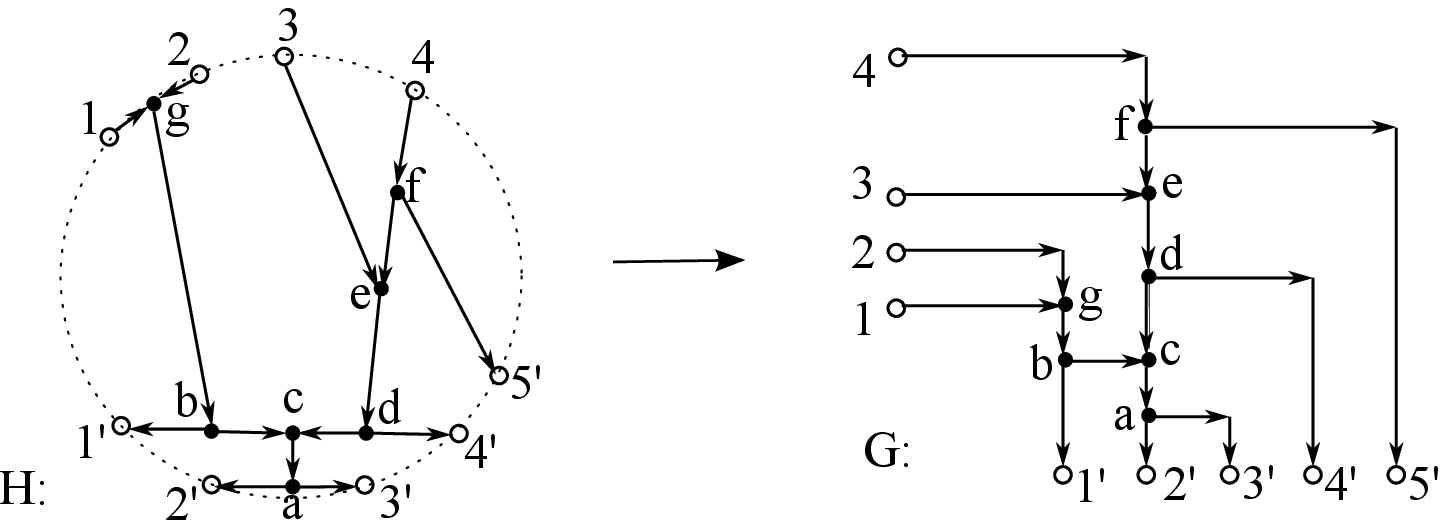}
\end{center}
\vspace{0cm}


\section{\Large Concluding remarks and additional results}  \label{sec:concl}

\subsection{An open question.}
~It looks reasonable to ask: how narrow is the class of UQ identities for
minors of $q$-matrices compared with the class of those in the commutative
version. We know that the latter class is formed by balanced families
$\Iscr,\Kscr$, whereas the former one is characterized via a stronger property
of $q$-balancedness. So we can address the problem of characterizing the set of
(homogeneous) balanced families $\Iscr,\Kscr\subset
\Escr^{m,n}\times\Escr^{m,n}$ that admit functions $\alpha:\Iscr\to\Zset$ and
$\beta:\Kscr\to\Zset$ such that the quadruple $\Iscr,\Kscr,\alpha,\beta$ is
$q$-balanced.

In an algorithmic setting, we deal with the following problem~$(\ast)$: given
$\Iscr,\Kscr$ (as above), decide whether or not there exist corresponding
$\alpha,\beta$ (as above). Concerning algorithmic complexity aspects, note that
the number $|\bfC(\Iscr)|+|\bfC(\Kscr)|$ of configurations for $\Iscr,\Kscr$
may be exponentially large compared with the number $|\Iscr|+|\Kscr|$ of
corteges (since a cortege of size $N$ may have $2^{O(N)}$ feasible matchings).
In light of this, it is logically reasonable to regard as the input of
problem~($\ast$) just the set $\bfC(\Iscr)\sqcup \bfC(\Kscr)$ rather than
$\Iscr\sqcup \Kscr$ (and measure the input size of~($\ast$) accordingly). We
conjecture that problem~($\ast$) specified in this way is NP-hard and,
moreover, it remains NP-hard even in the flag case.

\subsection{Non-quasicommuting flag minors.}
~The simplest example of balanced $\Iscr,\Kscr$ for which problem~($\ast$) has
answer ``not'' arises in the flag case with $\Iscr,\Kscr$ consisting of single
corteges. That is, we deal with quantized flag minors $[I]=[A|I]$ and
$[J]=[B|J]$, where $A:=\{1,\ldots,|I|\}$ and $B:=\{1,\ldots,|J|\}$, and
consider the (trivially balanced) one-element families $\Iscr=\{S:=(A|I,B|J)\}$
and $\Kscr=\{T:=(B|J,A|I)\}$. By Leclerc--Zelevinsky's theorem
(Theorem~\ref{tm:LZ}), $[I]$ and $[J]$ quasicommute if and only if the sets
$I,J$ are weakly separated. We have explained how to obtain ``if'' part of this
theorem by use of the flow-matching method, and now we explain how to use this
method to show, relatively easily, ``only if'' part (which has a rather
sophisticated proof in~\cite{LZ}).

So, assuming that $I,J$ are not weakly separated, our aim is to show that there
do not exist $\alpha(S),\beta(T)\in\Zset$ such that the equality
  \begin{equation} \label{eq:bTaS}
   \beta(T)-\alpha(S)=\zeta^\circ (S;M)-\zeta^\bullet (S;M)
   \end{equation}
holds for all feasible matching $M$ for $S$. The crucial observation is that
  \begin{numitem1} \label{eq:1match}
  $I,J\subset[n]$ are weakly separated if and only if $S$ has exactly one
feasible matching
   \end{numitem1}
(where ``only if'' part, mentioned in~\SSEC{quasicommut}, is trivial). In fact,
we need a sharper version of ``if'' part of~\refeq{1match}: when
$I,J\subset[n]$ are not weakly separated, there exist $M,M'\in\Mscr(S)$ such
that
   \begin{equation} \label{eq:MMp}
   \zeta^\circ (S;M)-\zeta^\bullet (S;M)\ne \zeta^\circ (S;M')-\zeta^\bullet
   (S;M').
   \end{equation}
Then the fact that the exchange operation applied to $S$ using $M$ results in
$T$, and similarly for $M'$, implies that~\refeq{bTaS} cannot hold
simultaneously for both $M$ and $M'$.

To construct the desired $M$ and $M'$, we argue as follows. Let for
definiteness $|I|\ge |J|$ and let $\Iw:=I-J$ and $\Jb:=J-I$. From the property
that $I,J$ are not weakly separated one can conclude that there are $a,b\in[n]$
with $a<b$ such that the sets $\tilde\Iw:=\{i\in \Iw\colon a\le i\le b\}$ and
$\tilde J^\bullet:=\{j\in\Jb\colon a\le j\le b\}$ satisfy $|\tilde
\Iw|-1=|\tilde J^\bullet|=:k$, and $\tilde\Iw$ has a partition into nonempty
sets $I_1,I_2$ satisfying $I_1<\tilde J^\bullet<I_2$. Let
  $$
  I_1=(i_1<i_2<\ldots<i_p),\quad I_2=(i_{p+1}<\ldots<i_{k+1}),\quad
  \tilde J^\bullet=(j_1<\ldots<j_k)
  $$
(then $i_p<j_1$ and $j_k<i_{p+1}$). Choose an arbitrary matching
$M\in\Mscr(S)$, and consider the set $\varPi$ of couples in $M$ containing
elements of $\tilde J^\bullet$; let $\varPi=\{\pi_1,\ldots,\pi_k\}$, where
$j_\ell\in\pi_\ell$. Each $\pi_\ell$ is a $C$-couple (since it cannot be an
$RC$- couple, in view of $B-A=\emptyset$), and the feasibility
condition~\refeq{feasM} for $M$ implies that only two cases are possible: (a)
$p$ couples in $\varPi$ meet $I_1$ and the remaining $k-p$ couples meet $I_2$,
and (b) $p-1$ couples in $\varPi$ meet $I_1$ and the remaining $k-p+1$ couples
meet $I_2$.

In case (a), we have $\pi_\ell=\{j_\ell,i_{p-\ell+1}\}$ for $\ell=1,\ldots,p$,
and $\pi_\ell=\{j_\ell,i_{\ell}\}$ for $\ell=p+1,\ldots,k$. An especial role is
played by the couple in $M$ containing the last element $i_{k+1}$ of $I_2$,
say, $\pi=\{i_{k+1},d\}$ (note that $d$ belongs to either $A-B$ or $\Jb-\tilde
J^\bullet$). We modify $M$ by replacing the couple $\pi$ by $\pi':=\{i_1,d\}$,
and replacing $\pi_p=\{j_p,i_1\}$ by $\pi'_p:=\{j_p,i_{k+1}\}$, forming
matching $M'$. The picture illustrates the case $k=3$, $p=2$ and $d\in A-B$.

\vspace{0cm}
\begin{center}
\includegraphics{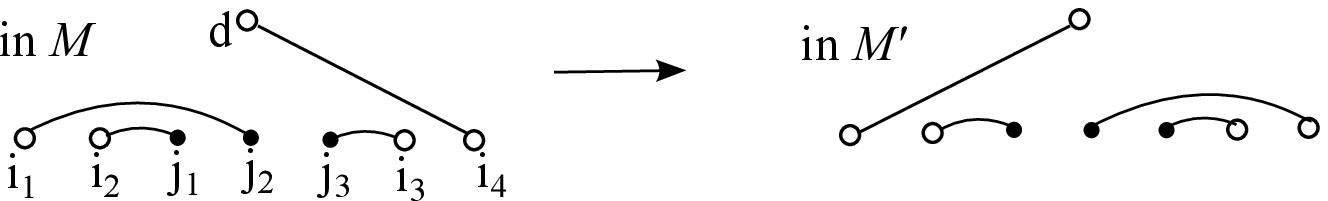}
\end{center}
\vspace{0cm}

One can see that $M'$ is feasible for $S$. Moreover, $M$ and $M'$
satisfy~\refeq{MMp}. Indeed, $\pi_p$ contributes one unit to $\zeta^\circ(S;M)$
while $\pi'_p$ contributes one unit to $\zeta^\bullet(S;M')$, the contributions
from $\pi$ and from $\pi'$ are the same, and the rests of $M$ and $M'$
coincide.

Thus, in case~(a), the one-element families $\{S\}$ and $\{T\}$ along with any
numbers $\alpha(S),\beta(T)$ are not $q$-balanced. Then
relation~\refeq{quasiIJ} (with any $c$) is impossible by
Theorem~\ref{tm:nec_q_bal}. In case~(b), the argument is similar. This yields
the necessity (``only if'' part) in Theorem~\ref{tm:LZ}. \hfill\qed

\subsection{Quasicommuting general minors.}
~It is tempting to ask: can one characterize the set of quasicommuting quantum
minors in a general (non-flag) case? Such a characterization can be obtained,
without big efforts, by use of the flow-matching method, yielding a
generalization of Theorem~\ref{tm:LZ}.
  \begin{theorem} \label{tm:generalLZ}
Let $(I|J),(I'|J)\in\Escr^{m,n}\times \Escr^{m,n}$ and let $|I|\ge|I'|$. The
following statements are equivalent:

{\rm (i)} the minors $[I|J]$ and $[I'|J']$ quasicommute, i.e., $[I|J][I'|J'] =
q^c [I'|J'][I|J]$ for some $c\in\Zset$;

{\rm(ii)} the cortege $S=(I|J,I'|J')$ admits exactly one feasible matching;

{\rm(iii)} the sets $I,I'$ are weakly separated, the sets $J,J'$ are weakly
separated, and for the refinement $(\Iw,\Ib,\Jw,\Jb)$ of $S$, one of the
following takes place:

{\rm(a)} $|\Ib| |\Jb|=0$; or

{\rm(b)} both sets $\Ib,\Jb$ are nonempty, and either $\Iw<\Ib$ and $\Jb<\Jw$,
or $\Ib<\Iw$ and $\Jw<\Jb$.

Also in case~(iii) the number $c$ is computed as follows: if $\Ib=\emptyset$,
~$\Jb=J_1\cup J_2$ and $J_1<\Jw<J_2$, then $c=|J_2|-|J_1|$; (symmetrically) if
$\Jb=\emptyset$, ~$\Ib=I_1\cup I_2$ and $I_1<\Iw<I_2$, then $c=|I_2|-|I_1|$;
~if $\Iw<\Ib$ and $\Jb<\Jw$, then $c=|\Ib|-|\Jb|$; and (symmetrically) if
$\Ib<\Iw$ and $\Jw<\Jb$, then $c=|\Jb|-|\Ib|$.
  \end{theorem}
  \begin{proof}
~Implication (ii)$\to$(i) is proved as in Section~\SSEC{quasicommut}, and
(iii)$\to$(ii) is easy.

To show (i)$\to$(iii), we use the fact that $|\Iw|-|\Ib|=|\Jw|-|\Jb|\ge 0$
(cf.~\refeq{balancIJ}) and observe that a feasible matchings for $S$ can be
constructed by the following procedure~(P) consisting of three steps. First,
choose an arbitrary maximal feasible set $\Mr$ of $R$-couples in
$\Yr:=\Iw\cup\Ib$. Here the feasibility means that the elements of each couple
have different colors and there are neither couples $\{i,j\}$ and $\{p,r\}$
with $i<p<j<r$, nor a couple $\{i,j\}$ and an element $d\in \Yr-\cup(\pi\in
\Mr)$ with $i<d<j$; cf.~\refeq{feasM}. Second, choose an arbitrary maximal
feasible set $\Mc$ of $C$-couples in $\Yc:=\Jw\cup\Jb$. Third, when $|I|>|I'|$,
the remaining elements of $\Yr\sqcup\Yc$ (which are all white) are coupled by a
unique set $\Mrc$ of $RC$-couples. Then $M:=\Mr\cup\Mc\cup\Mrc$ is a feasible
matching for $S$.

Suppose that (iii) is false and consider possible cases.

1) Let $J,J'$ be not weakly separated. Then we construct $\Mr,\Mc,\Mrc$ by
procedure~(P) and work with the matching $\tilde M:=\Mc\cup\Mrc$ in a similar
way as in the above proof for the flag case (with non-weakly-separated column
sets). This transforms $\tilde M$ into $\tilde M'$, and we obtain two different
feasible matchings $M:=\tilde M\cup\Mr$ and $M':=\tilde M'\cup\Mr$ for $S$
satisfying~\refeq{MMp}. This leads to a contradiction with (i) (as well
as~(ii)) in the theorem. When $I,I'$ are not weakly separated, the argument is
similar.

2) Assuming that $I,I'$ are weakly separated, and similarly for $J,J'$, let
both $\Ib \Jb$ be nonnempty. Then $\Iw,\Jw$ are nonempty as well, and for the
matching $M$ formed by procedure~(P),  $\Mr$ covers $\Ib$ and $\Mc$ covers
$\Jb$.

Denote by $a,a'$ (resp. $b,b'$) the minimal and maximal elements in $\Yr$
(resp. $\Yc$), respectively. Suppose that both $a,b$ are black. Then we can
transform $M$ into $M'$ by replacing the $R$-couple containing $a$, say, $ad$,
and the $C$-couple containing $b$, say, $bf$, by the two $RC$-couples $ab$ and
$df$. It is easy to see that $M'$ is feasible and $M,M'$ satisfy~\refeq{MMp}
(since under the transformation $M\to M'$ the value $\zeta^\circ-\zeta^\bullet$
decreases by two), whence~(i) is false. When both $a',b'$ are black, we act
similarly. So we may assume that each pair $\{a,b\}$ and $\{a',b'\}$ contains a
white element. The case $a\in\Iw$ and $b\in\Jw$ is possible only if
$|\Iw|=|\Ib|$ (taking into account that $|\Iw|\ge |\Ib|$ and that $\Iw,\Ib$, as
well as $\Jw,\Jb$, are weakly separated), implying $|\Jw|=|\Jb|$. But then
$\Mr$ covers $\Iw$ and $\Mc$ covers $\Jw$; so we can construct a feasible
matching $M'\ne M$ as in the previous case (after changing the colors
everywhere). And similarly when both $a',b'$ are white.

Thus, we may assume that $a,b$ have different colors, and so are $a',b'$.
Suppose that $a,a'\in \Iw$ and $b,b'\in \Jb$ (the case $a,a'\in \Ib$ and
$b,b'\in \Jw$ is similar). This is possible only if $|\Iw|=|\Ib|$ (since
$|I|\ge|I'|$, and $I,I'$ are weakly separated). Then the feasible matching $M$
constructed by~(P) consists of only $R$- and $C$-couples. Take the $R$-couple
in $M$ containing $a$ and the $C$-couple containing $b'$, say, $\pi=\{a,i\}$
and $\pi'=\{j,b'\}$; then both $a,j$ are white and both $i,b'$ are black.
Replace $\pi,\pi'$ by the $RC$-couples $\{a,j\}$ and $\{i,b'\}$. This gives a
feasible matching $M'\ne M$ satisfying~\refeq{MMp}.

The remaining situation is just as in~(a) or~(b) of~(iii), yielding
(i)$\to$(iii).
 \end{proof}

\noindent\underline{\em Remark 6.} Note that the situation when
$(\Iw,\Ib,\Jw,\Jb)$ has only one feasible matching can also be interpreted as
follows. Let us change the colors of all elements in the upper half of the
circumference $O$ (i.e., $\Iw$ becomes black and $\Ib$ becomes white). Then the
quantities of white and black elements in $O$ are equal and the elements of
each color go in succession cyclically.
\smallskip

\noindent\underline{\em Remark 7.} When minors $[I|J]$ and $[I'|J']$
quasicommute with $c=0$, we obtain the situation of ``purely commuting''
quantum minors, such as those discussed in Sect.~\SSEC{com_min}. The last
assertion in Theorem~\ref{tm:generalLZ} enables us to completely characterize
the set of corteges $(I|J,I'|J')$ determining commuting $q$-minors, as follows.

\begin{prop} \label{pr:gen_commute}
~$[I|J][I'|J']=[I'|J'][I|J]$ holds if and only if the refinement
$(\Iw,\Ib,\Jw,\Jb)$ satisfies at least one of the following:

(C1) $|\Iw|=|\Jw|$ (as well as $|\Ib|=|\Jb|$) and either $\Iw<\Ib$ and
$\Jb<\Jw$, or, symmetrically, $\Ib<\Iw$ and $\Jw<\Jb$;

(C2) assuming for definiteness that $|I|\ge|I'|$, either $\Ib=\emptyset$ and
$\Jb$ has a partition $J_1\cup J_2$ such that $|J_1|=|J_2|$ and $J_1<\Jw<J_2$,
or, symmetrically, $\Jb=\emptyset$ and $\Ib$ has a partition $I_1\cup I_2$ such
that $|I_1|=|I_2|$ and $I_1<\Iw<I_2$.
  \end{prop}

Cases (C1) and (C2) are illustrated in the picture by two level diagrams.

\vspace{0cm}
\begin{center}
\includegraphics{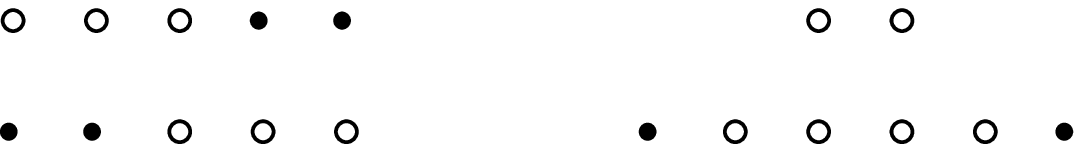}
\end{center}
\vspace{-0.3cm}

\subsection{Rotations.}
~Return to a general UQ identity~\refeq{gen_QR}. In Sect.~\SEC{q_relat} we
demonstrated two transformations of $q$-balanced $\Iscr,\Kscr,\alpha,\beta$
that preserve the $q$-balancedness (namely, the ones of \emph{reversing} and
\emph{transposing}, which result in $\Iscr^\reve,\Kscr^\reve,-\alpha,-\beta$
and $\Iscr^\top,\Kscr^\top,\alpha,\beta$, respectively.) Now we demonstrate one
more interesting (and nontrivial) transformation of $\Iscr,\Kscr,\alpha,\beta$
(in Theorem~\ref{tm:rotation}).

First, for corresponding $\Xr,\Yr\subset [m]$ and $\Xc,\Yc\subset[n]$
(cf.~\refeq{XrYr}), let $\Yr=(i_1<\cdots<i_k)$ and $\Yc=(j_1<\cdots<j_{k'})$.
Choose $g,h\in\Zset$ such that
  \begin{gather}
\mbox{$g+h\le k$ \;\;if $g,h\ge 0$}; \qquad \mbox{$|g|+|h'|\le k'$
           \;\;if $g,h\le 0$}; \label{eq:gh} \\
\mbox{$g\le k$ and $|h|\le k'$ \;\;if $g\ge 0\ge h$};
 \qquad \mbox{$|g|\le k'$ and $h\le k$ \;\;if $g\le 0\le h$}.
 \nonumber
  \end{gather}

Assuming that the numbers $i_1$, $m-i_k$, $j_1$, $n-j_{k'}$ are large enough,
we take sets $A,B\subset[m]$ and $A',B'\subset[n]$ such that $|A|=|A'|=|g|$,
$|B|=|B'|=|h'|$, $(A\cup B)\cap \Xr=\emptyset$, $(A'\cup B')\cap
\Xc=\emptyset$, and
 \begin{numitem1} \label{eq:AApBBp}
(a) $A=\{i_1,\ldots,i_g\}$ and $A'<\Yc$ if $g\ge 0$;

(a') $A<\Yr$ and $A'=\{j_1,\ldots, j_{|g|}\}$ if $g\le 0$;

(b) $B=\{i_{k-h+1}, \ldots, i_k\}$ and $B'>\Yc$ if $h\ge 0$;

(b') $B>\Yr$ and $B'=\{j_{k'-|h|+1},\ldots, j_{k'}\}$ if $h\le 0$.
  \end{numitem1}

Let $\xi$ ($\eta$) be the \emph{order-reversing} bijection between $A$ and $A'$
(resp. $B$ and $B'$), i.e., $\ell$-th element of $A$ is bijective to
$(|g|+1-\ell)$-th element of $A'$, and similarly for $\eta$.
  \smallskip

Second, we transform each cortege $S=(I|J,I'|J')\in\Iscr\cup\Kscr$ into cortege
$S_{g,h}=(\tilde I|\tilde J,\tilde I'|\tilde J')$ such that $\tilde I\cap
\tilde I'=\Xr$, \;$\tilde J\cap\tilde J'=\Xc$, and the refinement $(\tilde\Iw,
\tilde \Ib,\tilde \Jw,\tilde \Jb)$ of $S_{g,h}$ is expressed via the refinement
$(\Iw,\Ib,\Jw,\Jb)$ of $S$ as follows:
  \begin{itemize}
\item[(i)] $\tilde \Iw\cup \tilde \Ib=\Yr\pm A\pm B=:\Yr_{g,h}$ \; and $\tilde
\Jw\cup \tilde \Jb=\Yc\pm A'\pm B'=:\Yc_{g,h}$ (where we write $P+Q$ for $P\cup
Q$ in case $P\cap Q=\emptyset$, and write $P-Q$ for $P\setminus Q$ in case
$P\supseteq Q$);
\item[(ii)] If $i\in\Iw$ ($i\in\Ib$) is not in $A\cup B$, then $i\in\tilde \Iw$
(resp. $i\in\tilde \Ib$), and symmetrically, if $j\in\Jw$ ($j\in\Jb$) is not in
$A'\cup B'$, then $j\in\tilde \Jw$ (resp. $j\in\tilde \Jb$);
\item[(iii)] If $i\in\Iw$ ($i\in\Ib$) is in $A\cup B$, then the element
bijective to $i$ (by $\xi$ or $\eta$) belongs to $\tilde \Jb$ (resp.
$\tilde\Jw$); and symmetrically, if $j\in\Jw$ ($j\in\Jb$) is in $A'\cup B'$,
then the element bijective to $j$ belongs to $\tilde \Ib$ (resp. $\tilde\Iw$).
  \end{itemize}
(In other words, $\xi$ and $\eta$ change the colors of elements occurring in
$A,B,A',B'$). We call $\Yr_{g,h}, \Yc_{g,h}, S_{g,h}$ the
$(g,h)$-\emph{rotations} of $\Yr,\Yc,S$, respectively. Accordingly, we say that
$\Iscr_{g,h}^\circal:=\{S_{g,h}\colon S\in\Iscr\}$ is the $(g,h)$-rotation of
$\Iscr$, and similarly for $\Kscr$.

(This terminology is justified by the observation that if $g=-h$, then each
cortege $S$ is transformed as though being rotated (by $|g|$ positions
clockwise or counterclockwise) on the circular diagram on $\Yr\sqcup\Yc$;
thereby each element moving across the middle horizontal line of the diagram
changes its color.)

Third, extend $\xi$ and $\eta$ to the bijection $\rho:\Yr\sqcup\Yc\to
\Yr_{g,h}\sqcup \Yc_{g,h}$ so that $\rho$ be identical on $\Yr-(A\cup B)$ and
on $\Yc-(A'\cup B')$. Then a perfect matching $M$ on $\Yr\sqcup \Yc$ induces
the perfect matching $\{\rho(\pi)\colon \pi\in M\}$ on $\Yr_{g,h}\sqcup
\Yc_{g,h}$, denoted as $M_{g,h}$. An important property (which is easy to
check) is that
  \begin{numitem1} \label{eq:Mgh}
if $M$ is a feasible matching for $S\in\Iscr\cup\Kscr$, then $M_{g,h}$ is a
feasible matching for $S_{g,h}$, and vice versa.
  \end{numitem1}

An example of rotation of $S$ with $M\in\Mscr(S)$ is illustrated in the picture
where $k=5$, $k'=3$, $g=2$ and $h=-1$.

\vspace{0cm}
\begin{center}
\includegraphics[scale=0.8]{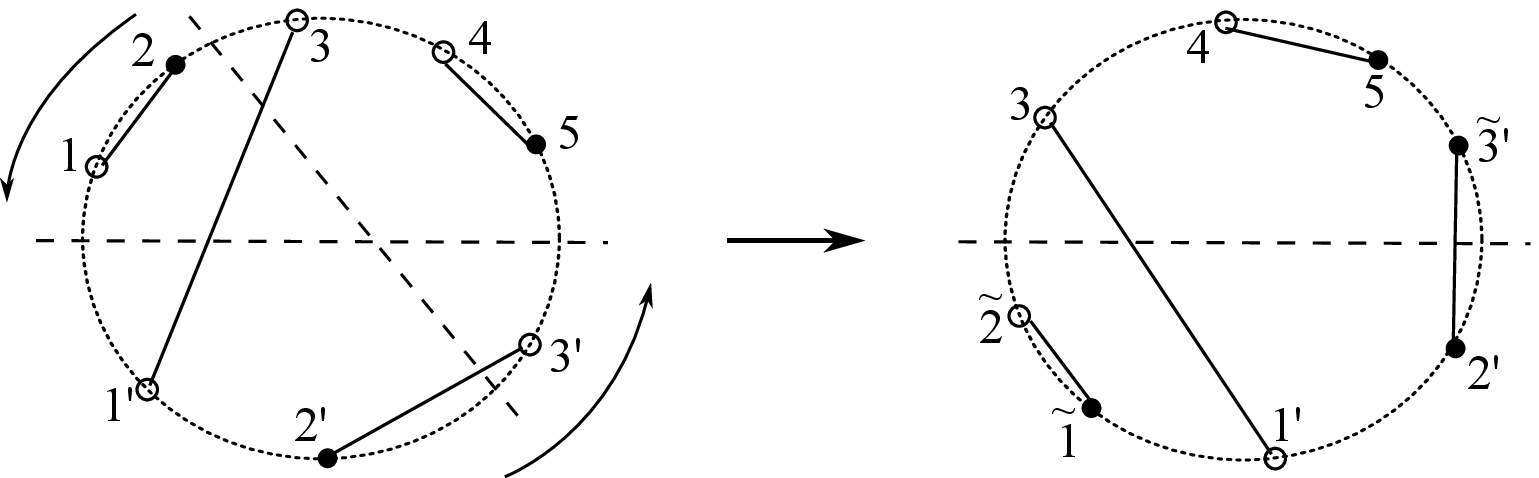}
\end{center}
\vspace{-0.3cm}

Fourth, for $S=(I|J,I'|J')$, define $\omega(S):=\delta_S(A)
+\delta_S(A')+\delta_S(B)+\delta_S(B')$, where
  \begin{gather}
  \delta_S(A):=|A\cap \Iw|, \qquad \delta_S(B):=-|B\cap \Iw|, \label{eq:deltaAB}\\
  \delta_S(A'):=|A'\cap \Jw|, \qquad \delta_S(B'):=-|B'\cap \Jw|. \nonumber
  \end{gather}

  %
  \begin{theorem} \label{tm:rotation}
Let $\Iscr,\Kscr,\alpha,\beta$ be $q$-balanced and let $g,h,$ be as
in~\refeq{gh}. Define $\alpha_{g,h}(S_{g,h}):=\alpha(S)+\omega(S)$ for
$S\in\Iscr$, and $\beta_{g,h}(T_{g,h}):=\beta(T)+\omega(T)$ for $T\in\Kscr$.
Then $\Iscr^\circal_{g,h},\Kscr^\circal_{g,h},\alpha_{g,h},\beta_{g,h}$ are
$q$-balanced.
  \end{theorem}
  \begin{proof}
Let $\gamma:\bfC(\Iscr)\to\bfC(\Kscr)$ be a bijection providing the
$q$-balancedness of $\Iscr,\Kscr,\alpha,\beta$. By~\refeq{Mgh}, $\gamma$
induces a bijection $\gamma_{g,h}:\bfC(\Iscr^\circal_{g,h})\to
\bfC(\Kscr^\circal_{g,h})$. More precisely, for configurations $(S;M)\in
\bfC(\Iscr)$ and $(T;M)=\gamma(S;M)$, ~$\gamma_{g,h}$ maps the configuration
$(S_{g,h};M_{g,h})$ to $(T_{g,h};M_{g,h})$. We assert that $\gamma_{g,h}$
satisfies the corresponding equality of the form
  \begin{equation} \label{eq:balanc_gh}
  \beta_{g,h}(T_{g,h})-\alpha_{g,h}(S_{g,h}) = \zeta^\circ(S_{g,h};\rho(\varPi))
   -\zeta^\bullet(S_{g,h}; \rho(\varPi))
   \end{equation}
(cf.~\refeq{q_balan}), yielding the result; here, as before, $\varPi$ is the
set of couples in $M$ having different colorings in the refinements of $S$ and
$T$.

For additivity reasons, it suffices to show~\refeq{balanc_gh} when $|g|+|h|=1$.
We will abbreviate corresponding $S_{g,h},T_{g,h},M_{g,h}$ as $S',T',M'$. (So
$T'$ is obtained from $S'$ by the exchange operation using
$\rho(\varPi)\subseteq M'$.) Let $d$ denote the (only) element of $\Yr\sqcup
\Yc$ that is not in $\Yr_{g,h}\sqcup\Yc_{g,h}$, and $\pi=\{d,f\}$ the couple in
$M$ containing $d$. Also we define $\Delta:=\zeta^\circ(S;\varPi)
-\zeta^\bullet(S; \varPi)$ and $\Delta':=\zeta^\circ(S';\rho(\varPi))
-\zeta^\bullet(S'; \rho(\varPi))$.

Our aim is to show that $\omega(T)-\omega(S)=\Delta'-\Delta$;
then~\refeq{balanc_gh} would immediately follow from~\refeq{q_balan}. One can
see that if $\pi\notin \varPi$, then $\Delta'=\Delta$, and
$\delta_S(D)=\delta_T(D)$ holds for $D=A,A',B,B'$ (cf.~\refeq{deltaAB}),
implying $\omega(S)=\omega(T)$. So we may assume that $\pi\in\varPi$. Consider
possible cases (where $S=(I|J,I'|J')$ and $T=(K|L,K'|L')$).
\smallskip

\noindent\emph{Case 1.} ~Let $g=1$. Then $d=i_1$. First suppose that $d\in\Iw$.
Then $\omega(S)=\delta_S(A)=1$ and $\omega(T)=\delta_T(A)=0$ (since the
exchange operation changes the color of $d$, i.e., $d\in \Kb$). If $\pi$ is an
$R$-couple for $S$, then $\pi$ contributes 1 to $\Delta$ (since $d$ is white
and $d<f$), and $\rho(\pi)$ contributes 0 to $\Delta'$ (since $\rho(\pi)$ is an
$RC$-couple for $S'$). Hence $\omega(T)-\omega(S)=-1=\Delta'-\Delta$, as
required. And if $\pi$ is an $RC$-couple for $S$, then $\pi$ contributes 0 to
$\Delta$ and $\rho(\pi)$ contributes $-1$ to $\Delta'$ (since $\rho(\pi)$ is a
$C$-couple for $S'$, $\rho(d)$ is black, $\rho(f)=f$ is white, and
$\rho(d)<f$), giving again $\Delta'-\Delta=-1$.

When $d\in\Ib$, we argue ``symmetrically'' (as though the roles of $S$ and $T$,
as well as $\zeta^\circ$ and $\zeta^\bullet$, are exchanged). Briefly, one can
check that: $\omega(S)=0$ and $\omega(T)=1$; if $\pi$ is an $R$-couple, then
$\pi$ contributes $-1$ to $\Delta$, and $\rho(\pi)$ contributes 0 to $\Delta'$;
and if $\pi$ is an $RC$-couple then $\pi$ contributes 0 to $\Delta$ and
$\rho(\pi)$ contributes 1 to $\Delta'$. Thus, every time we obtain
$\omega(T)-\omega(S)=1=\Delta'-\Delta$, as required.
  \smallskip

\noindent\emph{Case 2.} ~Let $h=1$. Then $d=i_k$. Suppose that $d\in \Iw$. Then
$\omega(S)=\delta_S(B)=-1$ and $\omega(T)=\delta_T(B)=0$. If $\pi$ is an
$R$-couple for $S$, then $\pi$ contributes $-1$ to $\Delta$ (since $d$ is white
and $d>f$) and $\rho(\pi)$ contributes 0 to $\Delta'$ (since $\rho(\pi)$ is an
$RC$-couple). And if $\pi$ is an $RC$-couple for $S$, then $\pi$ contributes 0
to $\Delta$ and $\rho(\pi)$ contributes 1 to $\Delta'$ (since $\rho(\pi)$ is a
$C$-couple for $S'$, $\rho(d)$ is black, and $\rho(d)>f$). In both cases, we
obtain $\omega(T)-\omega(S)=1=\Delta'-\Delta$, as required. When $d\in\Ib$, we
argue ``symmetrically''.
  \smallskip

Finally, the cases $g=-1$ and $h=-1$ are ``transposed'' to Cases 1 and 2,
respectively, and~\refeq{balanc_gh} follows by using relation~\refeq{transQR}.
  \end{proof}

\smallskip

\noindent\emph{Acknowledgements.} ~We thank Gleb Koshevoy for pointing out to
us paper~\cite{cast1}.



\newpage

\appendix


\section{\Large Appendix: Commutation properties of paths and a proof
of Theorem~\ref{tm:Linds}}  \label{sec:two_paths}

This Appendix contains auxiliary lemmas that are used in the proof of
Theorem~\ref{tm:Linds} given in this section as well, and in the proof of
Theorem~\ref{tm:single_exch} given in Appendix~B. These lemmas deal with
special pairs $P,Q$ of paths in an SE-graph $G=(V,E;R,C)$ and compare the
weights $w(P)w(Q)$ and $w(Q)w(P)$. Similar or close statements for Cauchon
graphs are given in~\cite{cast1,cast2}, and our method of proof is somewhat
similar and rather straightforward as well.

We first specify some terminology, notation and conventions.

When it is not confusing, vertices, edges, paths and other objects in $G$ are
identified with their corresponding images in the plane. We assume that the set
$R=\{r_1,\ldots,r_m\}$ of sources and the set $C=\{c_1,\ldots,c_n\}$ of sinks
lie on the coordinate rays $(0,\Rset_{\ge 0})$ and $(\Rset_{\ge 0},0)$,
respectively (then $G$ is disposed within the nonnegative quadrant
$\Rset^2_{\ge 0}$). The coordinates of a point $v$ in $\Rset^2$ (e.g., a vertex
$v$ of $G$) are denoted as $(\alpha(v),\beta(v))$. It is convenient to assume
that two vertices $u,v\in V$ have the same first (second) coordinate if and
only if they belong to a vertical (resp. horizontal) path in $G$, in which case
$u,v$ are called \emph{V-dependent} (resp. \emph{H-dependent}); for we always
can slightly perturb $G$ to ensure such a property, without affecting the graph
structure in essence. When $u,v$ are V-dependent, i.e., $\alpha(u)=\alpha(v)$,
we say that $u$ is \emph{lower} than $v$ (and $v$ is \emph{higher} than $u$) if
$\beta(u)< \beta(v)$. (In this case the commutation relation $uv=qvu$ takes
place.)

Let $P$ be a path in $G$. We denote: the first and last vertices of $P$ by
$s_P$ and $t_P$, respectively; the \emph{interior} of $P$ (the set of points of
$P-\{s_P,t_P\}$ in $\Rset^2$) by $\Inter(P)$; the set of horizontal edges of
$P$ by $E^H_P$; and the projection $\{\alpha(x)\;\colon x\in P\}$ by
$\alpha(P)$. Clearly if $P$ is directed, then $\alpha(P)$ is the interval
between $\alpha(s_P)$ and $\alpha(t_P)$.

For a directed path $P$, the following are equivalent: $P$ is non-vertical;
$E^H_P\ne \emptyset$; and $\alpha(s_P)\ne\alpha(t_P)$. We will refer to such a
$P$ as a \emph{standard} (rather than non-vertical) path.

For a standard path $P$, we will take advantage from a compact expression for
the weight $w(P)$. We call a vertex $v$ of $P$ \emph{essential} if either $P$
makes a turn at $v$ (changing the direction from horizontal to vertical or
back), or $v=s_P\not\in R$ and the first edge of $P$ is horizontal, or $v=t_P$
and the last edge of $P$ is horizontal. If $u_0,u_1,\ldots,u_k$ is the sequence
of essential vertices of $P$ in the natural order, then the weight of $P$ can
be expressed as
  \begin{equation} \label{eq:wP2}
  w(P)=u_0^{\sigma_0}u_1^{\sigma_1}\ldots u_k^{\sigma_k},
  \end{equation}
where $\sigma_i=1$ if $P$ makes a \horvert-turn at $u_i$ or if $i=k$, while
$\sigma_i=-1$ if $P$ makes a \verthor-turn at $u_i$ or if $i=0$ and $u_0$ is
the beginning of $P$. (Compare with~\refeq{telescop} where a path from $R$ to
$C$ is considered.) It is easy to see that if $P$ does not begin in $R$, then
its essential vertices are partitioned into H-dependent pairs.

Throughout the rest of the paper, for brevity, we denote $q^{-1}$ by $\bar q$,
and for an inner vertex $v\in W$ regarded as a generator, we may denote
$v^{-1}$ by $\bar v$.


\subsection{Auxiliary lemmas.} \label{ssec:aux}
~These lemmas deal with \emph{weakly intersecting} directed paths $P$ and $Q$,
which means that
  \begin{equation} \label{eq:pathsPQ}
P\cap Q=\{s_P,t_P\}\cap \{s_Q,t_Q\};
  \end{equation}
in particular, $\Inter(P)\cap\Inter(Q)=\emptyset$. For such $P,Q$, we say that
$P$ is \emph{lower} than $Q$ if there are points $x\in P$ and $y\in Q$ such
that $\alpha(x)=\alpha(y)$ and $\beta(x)<\beta(y)$ (then there are no $x'\in P$
and $y'\in Q$ with $\alpha(x')=\alpha(y')$ and $\beta(x')>\beta(y')$).

For paths $P,Q$, we define the value $\varphi=\varphi(P,Q)$ by the relation
  $$
  w(P)w(Q)=\varphi w(Q)w(P).
  $$
Obviously, $\varphi(P,Q)=1$ when $P$ or $Q$ is a V-path. In the lemmas below we
default assume that both $P,Q$ are standard.

  \begin{lemma} \label{lm:varphi=1}
Let $\{\alpha(s_P),\alpha(t_P)\}\cap\{\alpha(s_Q),\alpha(t_Q)\}\cap
\Rset_{>0}=\emptyset$. Then $\varphi(P,Q)=1$.
  \end{lemma}
  \begin{proof}
~Consider an essential vertex $u$ of $P$ and an essential vertex $v$ of $Q$.
Then for any $\sigma,\sigma'\in\{1,-1\}$, we have $u^\sigma
v^{\sigma'}=v^{\sigma'} u^\sigma$ unless $u,v$ are dependent.

Suppose that $u,v$ are V-dependent. From hypotheses of the lemma it follows
that at least one of the following is true:
$\alpha(s_P)<\alpha(u)<\alpha(t_P)$, or $\alpha(s_Q)<\alpha(v)<\alpha(t_Q)$.
For definiteness assume the former. Then there is another essential vertex $z$
of $P$ such that $\alpha(z)=\alpha(u)=\alpha(v)$. Moreover, $P$ makes a
\horvert-turn an one of $u,z$, and a \verthor-turn at the other. Since $P\cap
Q=\emptyset$ (in view of~\refeq{pathsPQ}), the vertices $u,z$ are either both
higher or both lower than $v$. Let for definiteness $u,z$ occur in this order
in $P$; then $w(P)$ contains the terms $u,\bar z$. Let $w(Q)$ contain the term
$v^\sigma$ and let $uv^\sigma=\rho v^\sigma u$, where $\sigma\in\{1,-1\}$ and
$\rho\in\{q,\bar q\}$. Then $\bar z v^\sigma= \bar \rho v^\sigma \bar z$,
implying $u\bar z v^\sigma =v^\sigma u\bar z$. Hence the contributions to
$w(P)w(Q)$ and $w(Q)w(P)$ from the pairs using terms $u,z,v$ (namely
$\{u,v^\sigma\}$ and $\{\bar z,v^\sigma\}$) are equal.

Next suppose that $u,v$ are H-dependent. One may assume that
$\alpha(u)<\alpha(v)$. Then $Q$ contains one more essential vertex $y\ne v$
with $\beta(y)=\beta(v)=\beta(u)$. Also $\alpha(u)<\alpha(v)$ and $P\cap
Q=\emptyset$ imply $\alpha(u)<\alpha(y)$. Let for definiteness
$\alpha(y)<\alpha(v)$. Then $w(Q)$ contains the terms $\bar y,v$, and we can
conclude that the contributions to $w(P)w(Q)$ and $w(Q)w(P)$ from the pairs
using terms $u,y,v$ are equal (using the fact that
$\alpha(u)<\alpha(y),\alpha(v)$).

These reasonings imply $\varphi(P,Q)=1$.
  \end{proof}

  \begin{lemma} \label{lm:asP=asQ}
Let $\alpha(s_P)=\alpha(s_Q)>0$ and $\alpha(t_P)\ne\alpha(t_Q)$. Let $P$ be
lower than $Q$. Then $\varphi(P,Q)=q$.
  \end{lemma}
  \begin{proof}
~Let $u$ and $v$ be the first essential vertices in $P$ and $Q$, respectively.
Then $\alpha(u)=\alpha(s_P)=\alpha(s_Q)=\alpha(v)$ (in view of
$\alpha(s_P)=\alpha(s_Q)>0$). Since $P$ is lower than $Q$, we have
$\beta(u)\le\beta(v)$. Moreover, this inequality is strong (since
$\beta(u)=\beta(v)$ is impossible in view of~\refeq{pathsPQ} and the obvious
fact that $u,v$ are the tails of first H-edges in $P,Q$, respectively).

Now arguing as in the above proof, we can conclude that the discrepancy between
$w(P)w(Q)$ and $w(Q)w(P)$ can arise only due to swapping the vertices $u,v$.
Since $u$ gives the term $\bar u$ in $w(P)$, and $v$ the term $\bar v$ in
$w(Q)$, the contribution from these vertices to $w(P)w(Q)$ and $w(Q)w(P)$ are
expressed as $\bar u\bar v$ and $\bar v \bar u$, respectively. Since
$\beta(u)<\beta(v)$, we have $\bar u \bar v=q\bar v \bar u$, and the result
follows.
  \end{proof}

  \begin{lemma} \label{lm:atP=atQ}
Let $\alpha(t_P)=\alpha(t_Q)$ and let either $\alpha(s_P)\ne\alpha(s_Q)$ or
$\alpha(s_P)=\alpha(s_Q)=0$. Let $P$ be lower than $Q$. Then $\varphi(P,Q)=q$.
  \end{lemma}
 \begin{proof}
~We argue in spirit of the proof of Lemma~\ref{lm:asP=asQ}. Let $u$ and $v$ be
the last essential vertices in $P$ and $Q$, respectively. Then
$\alpha(u)=\alpha(t_P)=\alpha(t_Q)=\alpha(v)$. Also $\beta(u)<\beta(v)$ (since
$P$ is lower than $Q$, and in view of~\refeq{pathsPQ} and the fact that $u,v$
are the heads of H-edges in $P,Q$, respectively). The condition on
$\alpha(s_P)$ and $\alpha(s_Q)$ imply that the discrepancy between $w(P)w(Q)$
and $w(Q)w(P)$ can arise only due to swapping the vertices $u,v$ (using
reasonings as in the proof of Lemma~\ref{lm:varphi=1}). Observe that $w(P)$
contains the term $u$, and $w(Q)$ the term $v$. So the generators $u,v$
contribute $uv$ to $w(P)w(Q)$, and $vu$ to $w(Q)w(P)$. Now $\beta(u)<\beta(v)$
implies $uv=qvu$, and the result follows.
 \end{proof}

  \begin{lemma} \label{lm:1atP=asQ}
Let $\alpha(t_P)=\alpha(s_Q)$ and $\beta(t_P)\ge\beta(s_Q)$. Then
$\varphi(P,Q)=q$.
  \end{lemma}
 \begin{proof}
~Let $u$ be the last essential vertex in $P$ and let $v,z$ be the first and
second essential vertices of $Q$, respectively (note that $z$ exists because of
$0<\alpha(s_Q)<\alpha(t_Q)$). Then
$\alpha(u)=\alpha(t_P)=\alpha(s_Q)=\alpha(v)<\alpha(z)$. Also
$\beta(u)\ge\beta(t_P) \ge \beta(s_Q)\ge \beta(v)=\beta(z)$. Let $Q'$ and $Q''$
be the parts of $Q$ from $s_Q$ to $z$ and from $z$ to $t_Q$, respectively. Then
$\alpha(P)\cap\alpha(Q'')=\emptyset$, implying $\varphi(P,Q'')=1$ (using
Lemma~\ref{lm:varphi=1} when $Q''$ is standard). Hence
$\varphi(P,Q)=\varphi(P,Q')$.

To compute $\varphi(P,Q')$, consider three possible cases.

(a) Let $\beta(u)>\beta(v)$. Then $u,v$ form the unique pair of dependent
essential vertices for $P,Q'$. Note that $w(P)$ contains the term $u$, and
$w(Q')$ contains the term $\bar v$. Since $\beta(u)>\beta(v)$, we have $u\bar
v=q\bar vu$, implying $\varphi(P,Q')=q$.

(b) Let $u=v$ and let $u$ be the unique essential vertex of $P$ (in other
words, $P$ is an H-path with $s_P\in R$). Note that $u=v$ and
$\beta(t_P)\ge\beta(s_Q)$ imply $t_P=u=v=s_Q$. Also $\alpha(u)<\alpha(z)$ and
$\beta(u)=\beta(z)$; so $u,z$ are dependent essential vertices for $P,Q'$ and
$uz=qzu$. We have $w(P)=u$ and $w(Q')=\bar u z$ (in view of $u=v$). Then $u\bar
u z=\bar u uz=q\bar u z u$ gives $\varphi(P,Q')=q$.

(c) Now let $u=v$ and let $y$ be the essential vertex of $P$ preceding $u$.
Then $t_P=u=v=s_Q$, ~$\beta(y)=\beta(u)=\beta(z)$, and
$\alpha(y)<\alpha(u)<\alpha(z)$. Hence $y,u,z$ are dependent, $w(P)$ contains
$\bar yu$, and $w(Q')=\bar uz$. We have
  $$
  \bar y u\bar u z= \bar y\bar u uz=(q\bar u\bar y)(qzu)
          =q^2 \bar u(\bar q z\bar y) u=q\bar u z\bar y u,
  $$
again obtaining $\varphi(P,Q')=q$.
 \end{proof}

   \begin{lemma} \label{lm:2atP=asQ}
Let $\alpha(t_P)=\alpha(s_Q)$ and  $\beta(t_P)<\beta(s_Q)$. Then
$\varphi(P,Q)=\bar q$.
  \end{lemma}
 \begin{proof}
~Let $u$ be the last essential vertex of $P$, and $v$ the first essential
vertex of $Q$. Then $\alpha(u)=\alpha(t_P)=\alpha(s_Q)=\alpha(v)$, and
$\beta(t_P)<\beta(s_Q)$ together with~\refeq{pathsPQ} implies
$\beta(u)<\beta(v)$. Also $w(P)$ contains $u$ and $w(Q)$ contains $\bar v$. Now
$u\bar v=\bar q \bar v u$ implies $\varphi(P,Q)=\bar q$.
 \end{proof}


\subsection{Proof of Theorem~\ref{tm:Linds}.} \label{ssec:q_determ}
~It can be conducted as a direct extension of the proof of a similar
Lindstr\"om's type result given by Casteels~\cite[Sec.~4]{cast1} for Cauchon
graphs. To make our description more self-contained, we outline the main
ingredients of the proof, leaving the details where needed to the reader.

Let $(I|J)\in\Escr^{m,n}$, $I=\{i(1)<\cdots <i(k)\}$ and
$J=\{j(1)<\cdots<j(k)\}$. Recall that an $(I|J)$-flow in an SE-graph $G$ (with
$m$ sources and $n$ sinks) consists of pairwise disjoint paths $P_1,\ldots,
P_k$ from the source set $R_I=\{r_{i(1)},\ldots,r_{i(k)}\}$ to the sink set
$C_J=\{c_{j(1)},\ldots,c_{j(k)}\}$, and (because of the planarity of $G$) we
may assume that each $P_d$ begins at $r_{i(d)}$ and ends at $c_{j(d)}$.
Besides, we are forced to deal with an arbitrary \emph{path system}
$\Pscr=(P_1,\ldots, P_k)$ in which for $i=1,\ldots,k$, ~$P_d$ is a directed
path in $G$ beginning at $r_{i(d)}$ and ending at $c_{j(\sigma(d))}$, where
$\sigma(1), \ldots,\sigma(k)$ are different, i.e., $\sigma=\sigma_\Pscr$ is a
permutation on $[k]$. (In particular, $\sigma_\Pscr$ is identical if $\Pscr$ is
a flow.)

We naturally partition the set of all path systems for $G$ and $(I|J)$ into the
set $\Phi(I|J)$ of $(I|J)$-flows and the rest $\Psi(I|J)$ (consisting of those
path systems that contain intersecting paths). The following property easily
follows from the planarity of $G$ (cf.~\cite[Lemma~4.2]{cast1}):
   \begin{numitem1} \label{eq:PiPi+1}
For any $\Pscr=(P_1,\ldots,P_k)\in\Psi(I|J)$, there exist two
\emph{consecutive} intersecting paths $P_d,P_{d+1}$.
  \end{numitem1}

The $q$-\emph{sign} of a permutation $\sigma$ is defined by
   $$
   \Sign_q(\sigma):=(-q)^{\ell(\sigma)},
   $$
where $\ell(\sigma)$ is the length of $\sigma$ (see Sect.~\SEC{prelim}).

Now we start computing the $q$-minor $[I|J]$ of the matrix $\Path_G$ with the
following chain of equalities:
   \begin{eqnarray*}
   [I|J]&=& \sum\nolimits_{\sigma\in S_k} \Sign_q(\sigma)
        \left( \prod\nolimits_{d=1}^{k} \Path_G(i(d)|j(\sigma(d))\right) \\
   &=& \sum\nolimits_{\sigma\in S_k} \Sign_q(\sigma)
       \left( \prod\nolimits_{d=1}^{k} \left(
                     \sum(w(P)\;\colon P\in \Phi_G(i(d)|j(\sigma(d))\right)\right) \\
   &=&\sum(\Sign_q(\sigma_\Pscr)w(\Pscr)\;\colon \Pscr\in\Phi(I|J)\cup\Psi(I|J)) \\
   &=&\sum(w(\Pscr)\;\colon \Pscr\in\Phi(I|J))
                        +\sum(\Sign_q(\sigma_\Pscr)w(\Pscr)\;\colon
                        \Pscr\in\Psi(I|J)).
   \end{eqnarray*}

Thus, we have to show that the second sum in the last row is zero. It will
follow from the existence of an involution $\eta:\Psi(I|J)\to\Psi(I|J)$ without
fixed points such that for each $\Pscr\in\Psi(I|J)$,
  \begin{equation} \label{eq:invol}
  \Sign_q(\sigma_\Pscr)w(\Pscr)=-\Sign_q(\sigma_{\eta(\Pscr)}) w(\eta(\Pscr)).
  \end{equation}

To construct the desired $\eta$, consider $\Pscr=(P_1,\ldots,P_k)\in\Psi(I|J)$,
take the minimal $i$ such that $P_i$ and $P_{i+1}$ meet, take the last common
vertex $v$ of these paths, represent $P_i$ as the concatenation $K\circ L$, and
$P_{i+1}$ as $K'\circ L'$, so that $t_K=t_{K'}=s_L=s_{L'}=v$, and exchange the
portions $L,L'$ of these paths, forming $Q_i:=K\circ L'$ and $Q_{i+1}:=K'\circ
L$. Then we assign $\eta(\Pscr)$ to be obtained from $\Pscr$ by replacing
$P_i,P_{i+1}$ by $Q_i,Q_{i+1}$. It is routine to check that $\eta$ is indeed an
involution (with $\eta(\Pscr)\ne\Pscr$) and that
   \begin{equation} \label{eq:ell+1}
  \ell(\sigma_{\eta(\Pscr)})=\ell(\sigma_\Pscr)+1,
  \end{equation}
assuming w.l.o.g. that $\sigma(i)<\sigma(i+1)$. On the other hand, applying to
the paths $K,L,K',L'$ Lemmas~\ref{lm:asP=asQ} and~\ref{lm:1atP=asQ}, one can
obtain
   \begin{multline*} 
   \quad w(P_i)w(P_{i+1})=w(K)w(L)w(K')w(L')=qw(K)w(L)w(L')w(K')  \\
  =q^2 w(K)w(L')w(L)w(K')=qw(K)w(L')w(K')w(L)=qw(Q_i)w(Q_{i+1}),\quad
   \end{multline*}
whence $w(\Pscr)=qw(\eta(\Pscr))$. This together with~\refeq{ell+1} gives
   %
   \begin{equation*}
  \Sign_q(\sigma_\Pscr)w(\Pscr)+\Sign_q(\sigma_{\eta(\Pscr)}) w(\eta(\Pscr))
  =(-q)^{\ell(\sigma_\Pscr)}q w(\eta(\Pscr))+(-q)^{\ell(\sigma_\Pscr)+1}
  w(\eta(\Pscr))   =0,
  \end{equation*}
yielding~\refeq{invol}, and the result follows. \hfill\qed


\section{\Large Appendix: Proof of Theorem~\ref{tm:single_exch}}  \label{sec:exchange}

Using notation as in the hypotheses of this theorem, we first consider the case
when

\begin{itemize}
\item[(C):] $\pi=\{f,g\}$ is a $C$-couple in $M(\phi,\phi')$ with $f<g$ and $f\in
J$.
  \end{itemize}
(Then $f\in\Jw$ and $g\in\Jb$.) We have to prove that
  \begin{equation} \label{eq:caseC}
  w(\phi)w(\phi')=qw(\psi)w(\psi')
  \end{equation}
The proof is given throughout Sects.~\SSEC{seglink}--\SSEC{degenerate}. The
other possible cases in Theorem~\ref{tm:single_exch} will be discussed in
Sect.~\SSEC{othercases}.

 \subsection{Snakes and links.} \label{ssec:seglink}
Let $Z$ be the exchange path determined by $\pi$ (i.e., $Z=P(\pi)$ in notation
of Sect.~\SEC{double}). It connects the sinks $c_f$ and $c_g$, which may be
regarded as the first and last vertices of $Z$, respectively. Then $Z$ is
representable as a concatenation $Z=\bar Z_1\circ Z_2\circ\bar Z_3\circ \ldots
\circ \bar Z_{k-1}\circ Z_k$, where $k$ is even, each $Z_i$ with $i$ odd (even)
is a directed path concerning $\phi$ (resp. $\phi'$), and $\bar Z_i$ stands for
the path reversed to $Z_i$. More precisely, let $z_0:=c_f$, ~$z_k:=c_g$, and
for $i=1,\ldots,k-1$, denote by $z_i$ the common endvertex of $Z_i$ and
$Z_{i+1}$. Then each $Z_i$ with $i$ odd is a directed path from $z_i$ to
$z_{i-1}$ in $\langle E_\phi- E_{\phi'}\rangle$, while each $Z_i$ with $i$ even
is a directed path from $z_{i-1}$ to $z_i$ in $\langle E_{\phi'}-
E_{\phi}\rangle$.

We refer to $Z_i$ with $i$ odd (even) as a \emph{white} (resp. \emph{black})
\emph{snake}.

Also we refer to the vertices $z_1,\ldots,z_{k-1}$ as the \emph{bends} of $Z$.
A bend $z_i$ is called a \emph{peak} (a \emph{pit}) if both path $Z_i,Z_{i+1}$
leave (resp. enter) $z_i$; then $z_1,z_3,\ldots,z_{k-1}$ are the peaks, and
$z_2,z_4,\ldots,z_{k-2}$ are the pits. Note that some peak $z_i$ may coincide
with some pit $z_j$; in this case we say that $z_i,z_j$ are \emph{twins}.

The rests of flows $\phi$ and $\phi'$ consist of directed paths that we call
\emph{white} and \emph{black links}, respectively. More precisely, the white
(black) links correspond to the connected components of the subgraph $\phi$
(resp. $\phi'$) from which the interiors of all snakes are removed. So a link
connects either (a) a source and a sink (being a component of $\phi$ or
$\phi'$), or (b) a source and a pit, or (c) a peak and a sink, or (d) a peak
and a pit. We say that a link is \emph{unbounded} in case (a),
\emph{semi-bounded} in cases (b),(c), and \emph{bounded} in case (d). Note that
  \begin{numitem1} \label{eq:4paths}
a bend $z_i$ occurs as an endvertex in exactly four paths among snakes and
links, namely: either in two snakes and two links (of different colors), or in
four snakes $Z_i,Z_{i+1},Z_j,Z_{j+1}$ (when $z_i,z_j$ are twins).
  \end{numitem1}

We denote the sets of snakes and links (for $\phi,\phi',\pi$) by $\Sscr$ and
$\Lscr$, respectively; the corresponding subsets of white and black elements of
these sets are denoted as $\Sscr^\circ,\; \Sscr^\bullet,\; \Lscr^\circ,\;
\Lscr^\bullet$.

The picture below illustrates an example. Here $k=10$, the bends
$z_1,\ldots,z_9$ are marked by squares, the white and black snakes are drawn by
thin and thick solid zigzag lines, respectively, the white links
($L_1,\ldots,L_7$) by short-dotted lines, and the black links
($M_1,\ldots,M_6$) by long-dotted lines.

\vspace{0cm}
\begin{center}
\includegraphics{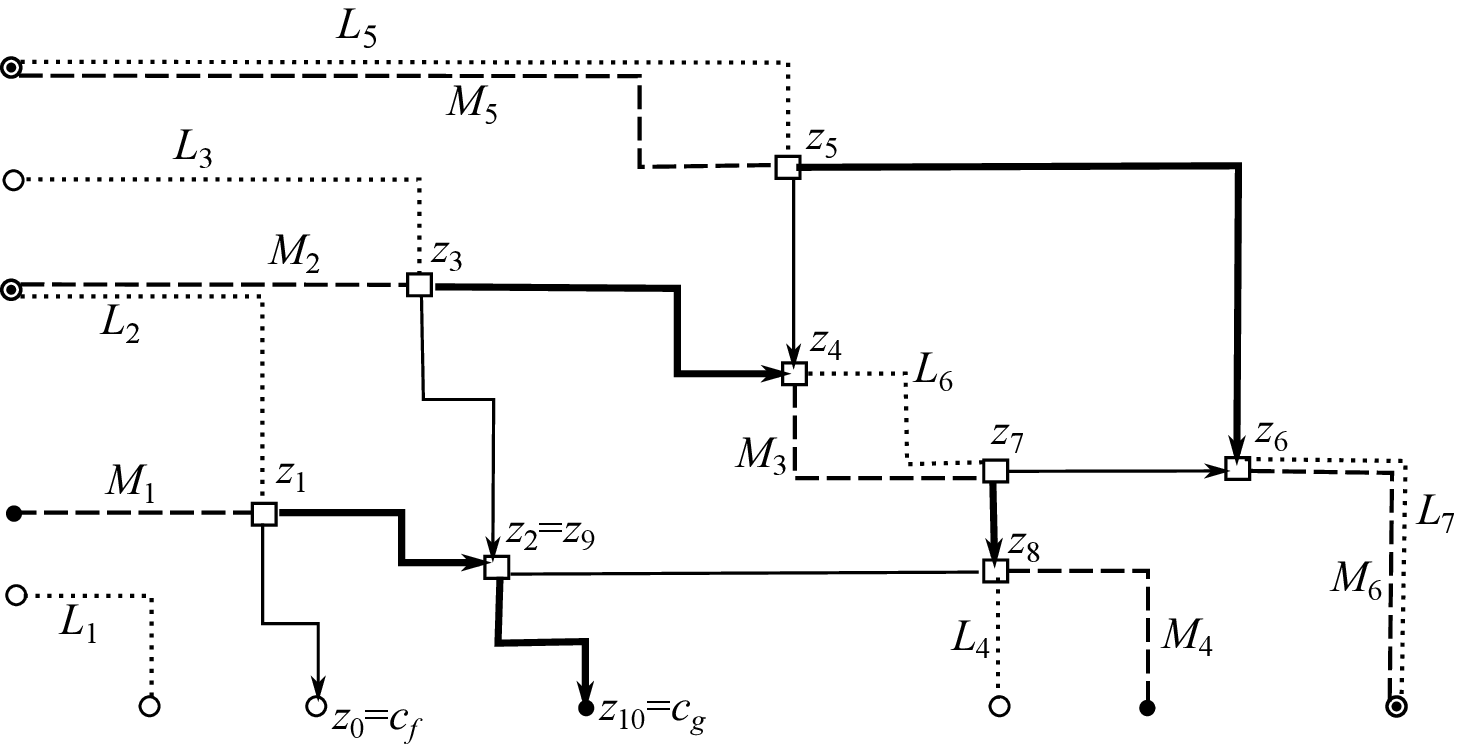}
\end{center}
\vspace{0cm}

The weight $w(\phi)w(\phi')$ of the double flow $(\phi,\phi')$ can be written
as the corresponding ordered product of the weights of snakes and links; let
$\Nscr$ be the string (sequence) of snakes and links in this product. The
weight of the double flow $(\psi,\psi')$ uses a string consisting of the same
snakes and links but occurring in another order; we denote this string by
$\Nscr^\ast$.

We say that two elements among snakes and links are \emph{invariant} if they
occur in the same order in $\Nscr$ and $\Nscr^\ast$, and \emph{permuting}
otherwise. In particular, two links of different colors are invariant, whereas
two snakes of different colors are always permitting.

For example, observe that the string $\Nscr$ for the above illustration is
viewed as
  $$
L_1L_2Z_1L_3Z_3Z_9L_4L_5Z_5L_6Z_7L_7M_1Z_2Z_{10}M_2Z_4M_3Z_8M_4M_5Z_6M_6,
  $$
whereas $\Nscr^\ast$ is viewed as
  $$
L_1L_2Z_2Z_{10}L_3Z_4L_6Z_8L_4L_5Z_6L_7M_1Z_1M_2Z_3Z_9M_4M_5Z_5M_3Z_7M_6.
  $$

For $A,B\in\Sscr\cup\Lscr$, we write $A\prec B$ (resp. $A\prec^\ast B$) if $A$
occurs in $\Nscr$ (resp. in $\Nscr^\ast$) earlier than $B$. We define
$\varphi_{A,B}=\varphi_{B,A}:=1$ if $A,B$ are invariant, and define
$\varphi_{A,B}=\varphi_{B,A}$ by the relation
   \begin{equation} \label{eq:phiAB}
   w(A)w(B)=\varphi_{A,B} w(B)w(A).
   \end{equation}
if $A,B$ are permuting and $A\prec B$. Note that $\varphi_{A,B}$ is defined
somewhat differently than $\varphi(P,Q)$ in Sect.~\SSEC{aux}.

For $A,B\in\Sscr\cup\Lscr$, we may use notation $(A,B)$ when $A,B$ are
permuting and $A\prec B$ (and may write $\{A,B\}$ when their orders by $\prec$
and $\prec^\ast$ are not important for us).

Our goal is to prove that in case~(C),
  \begin{equation}\label{eq:Pi=q}
  \prod(\varphi_{A,B}\;\colon A,B\in \Sscr\cup\Lscr)=q,
  \end{equation}
whence~\refeq{caseC} will immediately follow.

We first consider the \emph{non-degenerate} case. This means the following
restriction:
  \begin{numitem1} \label{eq:nondegenerate}
all coordinates $\alpha(z_1),\ldots,\alpha(z_{k-1}),
\alpha(c_1),\ldots,\alpha(c_n)$ of bends and sinks are different.
  \end{numitem1}

The proof of~\refeq{Pi=q} subject to~\refeq{nondegenerate} will consist of
three stages I, II, III where we compute the total contribution from the pairs
of links, the pairs of snakes, and the pairs consisting of one snake and one
link, respectively. As a consequence, the following three results will be
obtained (implying~\refeq{Pi=q}).
  \begin{prop} \label{pr:link-link}
In case~\refeq{nondegenerate}, the product $\varphi^I$ of the values
$\varphi_{A,B}$ over links $A,B\in\Lscr$ is equal to 1.
  \end{prop}
  \begin{prop} \label{pr:seg-seg}
In case~\refeq{nondegenerate}, the product $\varphi^{II}$ of the values
$\varphi_{A,B}$ over snakes $A,B\in\Sscr$ is equal to q.
  \end{prop}
  \begin{prop} \label{pr:seg-link}
In case~\refeq{nondegenerate}, the product $\varphi^{III}$ of the values
$\varphi_{A,B}$ where one of $A,B$ is a snake and the other is a link is equal
to 1.
  \end{prop}

These propositions are proved in Sects.~\SSEC{prop1}--\SSEC{prop3}. Sometimes
it will be convenient for us to refer to a white (black) snake/link concerning
$\phi,\phi',\pi$ as a $\phi$-snake/link (resp. a $\phi'$-snake/link), and
similarly for $\psi,\psi',\pi$.

 \subsection{Proof of Proposition~\ref{pr:link-link}.} \label{ssec:prop1}
Under the exchange operation using $Z$, any $\phi$-link becomes a $\psi$-link
and any $\phi'$-link becomes a $\psi'$-link. The white links occur in $\Nscr$
earlier than the black links, and similarly for  $\Nscr^\ast$. Therefore, if
$A,B$ are permuting links, then they are of the same color. This implies that
$A\cap B=\emptyset$. Also each endvertex of any link either is a bend or
belongs to $R\cup C$. Then~\refeq{nondegenerate} implies that the sets
$\{\alpha(s_A),\alpha(t_A)\}\cap \Rset_{>0}$ and
$\{\alpha(s_B),\alpha(t_B)\}\cap \Rset_{>0}$ are disjoint. Now
Lemma~\ref{lm:varphi=1} gives $\varphi_{A,B}=1$, and the proposition follows.
\hfill\qed

 \subsection{Proof of Proposition~\ref{pr:seg-seg}.} \label{ssec:prop2}

Consider two snakes $A=Z_i$ and $B=Z_j$, and let $A\prec B$. If $|i-j|>1$ then
$A\cap B=\emptyset$ and, moreover, $\{\alpha(s_A),\alpha(t_A)\}\cap
\{\alpha(s_B),\alpha(t_B)\}=\emptyset$ (since $Z$ is simple and in view
of~\refeq{nondegenerate}). This gives $\varphi_{A,B}=1$, by
Lemma~\ref{lm:varphi=1}.

Now let $|i-j|=1$. Then $A,B$ have different colors; hence $A$ is white and $B$
is black (in view of $A\prec B$). So $i$ is odd, and two cases are possible:
\smallskip

\noindent\underline{\emph{Case 1}:} ~$j=i+1$ and $z_i$ is a peak:
$z_i=s_A=s_B$;
  \smallskip

\noindent\underline{\emph{Case 2}:} ~$j=i-1$ and $z_{i-1}$ is a pit:
$z_{i-1}=t_A=t_B$.
  \smallskip

Cases 1,2 are divided into two subcases each.
 \smallskip

\noindent\underline{\emph{Subcase 1a}:} ~$j=i+1$ and $A$ is lower than $B$.
 \smallskip

\noindent\underline{\emph{Subcase 1b}:} ~$j=i+1$ and $B$ is lower than $A$.
 \smallskip

\noindent\underline{\emph{Subcase 2a}:} ~$j=i-1$ and $A$ is lower than $B$.
 \smallskip

\noindent\underline{\emph{Subcase 2b}:} ~$j=i-1$ and $B$ is lower than $A$.
 \smallskip

\noindent(The term \emph{lower} is explained in Sect.~A.) Subcases~1a--2b are
illustrated in the picture:

\vspace{-0cm}
\begin{center}
\includegraphics{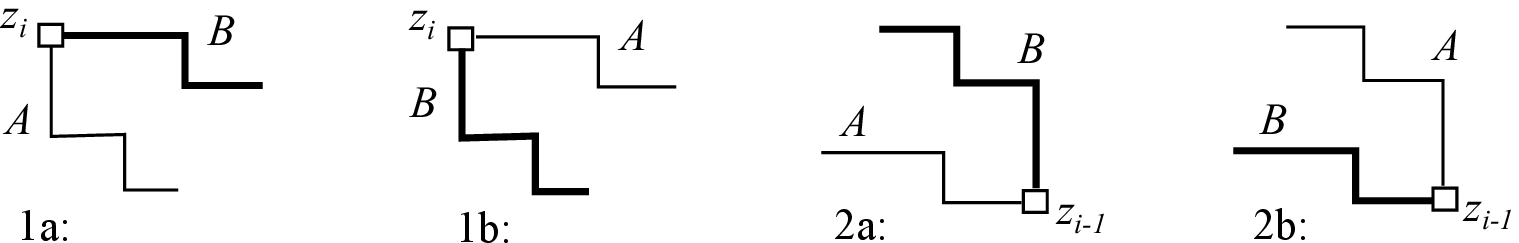}
\end{center}
\vspace{0cm}

Under the exchange operation using $Z$, any snake changes its color; so $A,B$
are permuting. Applying to $A,B$ Lemmas~\ref{lm:asP=asQ} and~\ref{lm:atP=atQ},
we obtain $\varphi_{A,B}=q$ in Subcases~1a,2a, and $\varphi_{A,B}=\bar q$ in
Subcases~1b,2b.

It is convenient to associate with a bend $z$ the number $\gamma(z)$ which is
equal to $+1$ if, for the corresponding pair $A\in\Sscr^\circ$ and
$B\in\Sscr^\bullet$ sharing $z$, ~$A$ is lower than $B$ (as in Subcases~1a,2a),
and equal to $-1$ otherwise (as in Subcases~1b,2b). Define
  \begin{equation} \label{eq:gammaZ}
  \gamma_Z:=\sum(\gamma(z)\;\colon z\;\; \mbox{a bend of}\;\; Z).
  \end{equation}
Then $\varphi^{II}=q^{\gamma_Z}$. Thus, $\varphi^{II}=q$ is equivalent to
   \begin{equation} \label{eq:gamma=1}
   \gamma_Z=1.
   \end{equation}

To show~\refeq{gamma=1}, we are forced to deal with a more general setting.
More precisely, let us turn $Z$ into simple cycle $D$ by combining the directed
path $Z_1$ (from $z_1$ to $z_0=c_f$) with the horizontal path from $c_f$ to
$c_g$ (to create the latter, we formally add to $G$ the horizontal edges
$(c_j,c_{j+1})$ for $j=f,\ldots,g-1$). The resulting directed path $\tilde Z$
from $z_1$ to $c_g=z_k$ is regarded as the new white snake replacing $Z_1$.
Then $\tilde Z_1$ shares the end $z_k$ with the black path $Z_k$; so $z_k$ is a
pit of $D$, and $\tilde Z$ is lower than $Z_k$. Thus, compared with $Z$, the
cycle $D$ acquires an additional bend, namely, $z_k$. We have $\gamma(z_k)=1$,
implying $\gamma_D=\gamma_Z+1$. Then~\refeq{gamma=1} is equivalent to
$\gamma_D=2$.

On this way, we come to a new (more general) setting by considering an
arbitrary simple (non-directed) cycle $D$ rather than a special path $Z$.
Moreover, instead of an SE-graph as before, we can work with a more general
directed planar graph $G$ in which any edge $e=(u,v)$ points arbitrarily within
the south-east sector, i.e., satisfies $\alpha(u)\le \alpha(v)$ and
$\beta(u)\ge \beta(v)$. We call $G$ of this sort a \emph{weak SE-graph}.

So now we are given a colored simple cycle $D$ in $G$, i.e., $D$ is
representable as a concatenation $\bar D_1\circ D_2\circ\ldots \circ \bar
D_{k-1}\circ D_k$, where each $D_i$ is a directed path in $G$; a path
(\emph{snake}) $D_i$ with $i$ odd (even) is colored white (resp. black). Let
$d_1,\ldots,d_k$ be the sequence of bends in $D$, i.e., $d_i$ is a common
endvertex of $D_{i}$ and $D_{i+1}$ (letting $D_{k+1}:=D_1$). We assume that $D$
is oriented according to the direction of $D_i$ with $i$ even. When this
orientation is clockwise (counterclockwise) around a point in the open bounded
region $O_D$ of the plane surrounded by $D$, we say that $D$ is
\emph{clockwise} (resp. \emph{counterclockwise}). Then the cycle arising from
the above path $Z$ is clockwise.

Our goal is to prove the following
  \begin{lemma} \label{lm:gammaD}
Let $D$ be a colored simple cycle in a weak SE-graph $G$. If $D$ is clockwise
then $\gamma_D=2$. If $D$ is counterclockwise then $\gamma_D=-2$.
  \end{lemma}
  \begin{proof}
~We use induction on the number $\eta(D)$ of bends of $D$. It suffices to
consider the case when $D$ is clockwise (since for a counterclockwise cycle
$D'=\bar D'_1\circ D'_2\circ\ldots \circ \bar D'_{k-1}\circ D'_k$, the reversed
cycle $\bar D'=\bar D'_k\circ D'_{k-1}\circ\ldots \circ \bar D'_2\circ D'_1$ is
clockwise, and it is easy to see that $\gamma_{\bar D'}=-\gamma_{D'}$).

W.l.o.g., one may assume that the coordinates $\beta(d_i)$ of all bends $d_i$
are different (as we can make, if needed, a due small perturbation on $D$,
which does not affect $\gamma$).

If $\eta(D)=2$, then $D=\bar D_1\circ D_2$, and the clockwise orientation of
$D$ implies that the path $D_1$ is lower than $D_2$. So
$\gamma(d_1)=\gamma(d_2)=1$, implying $\gamma_D=2$.

 Now assume that $\eta(D)>2$. Then at least one of the following is true:
 \smallskip

(a) there exists a peak $d_i$ such that the horizontal line through $d_i$ meets
$D$ on the left of $d_i$, i.e., there is a point $x$ in $D$ with
$\alpha(x)<\alpha(d_i)$ and $\beta(x)=\beta(d_i)$;

(b) there exists a pit $d_i$ such that the horizontal line through $d_i$ meets
$D$ on the right of $d_i$.
 \smallskip

(This can be seen as follows. Let $d_j$ be a peak with $\beta(d_j)$ maximum.
Then the clockwise orientation of $D$ implies that $D_{j+1}$ lies on the right
from $D_j$.  If $\beta(d_{j-1})< \beta(d_{j+1})$, then, by easy topological
reasonings, either the pit $d_{j+1}$ is as required in~(b) (when $d_{j+2}$ is
on the right from $D_{j+1}$), or the peak $d_{j+2}$ is as required in~(a) (when
$d_{j+2}$ is on the left from $D_{j+1}$), or both. And if $\beta(d_{j-1})>
\beta(d_{j+1})$, then $d_{j-1}$ is as in~(b).)

We may assume that case~(a) takes place (for case~(b) is symmetric to~(a), in a
sense). Choose the point $x$ as in~(a) with $\alpha(x)$ maximum and draw the
horizontal line-segment $L$ connecting the points $x$ and $d_i$. Then the
interior of $L$ does not meet $D$. Two cases are possible:
  \smallskip

(I) $\Inter(L)$ is contained in the region $O_D$; or
 \smallskip

(O) $\Inter(L)$ is outside $O_D$.
\smallskip

Since $x$ cannot be a bend of $D$ (in view of $\beta(x)=\beta(d_i)$ and
$\beta(d_i)\ne\beta(d_{i'})$ for any $i'\ne i$), $x$ is an interior point of
some snake $D_j$; let $D'_j$ and $D''_j$ be the parts of $D_j$ from $s_{D_j}$
to $x$ and from $x$ to $t_{D_j}$, respectively. Using the facts that $D$ is
oriented clockwise and this orientation is agreeable with the forward
(backward) direction of each black (resp. white) snake, one can realize that
  \begin{numitem1} \label{eq:casesIO}
(a) in case (I), ~$D_j$ is white and $\gamma(d_i)=-1$ (i.e., for the white
snake $D_i$ and black snake $D_{i+1}$ that share the peak $d_i$, ~$D_{i+1}$ is
lower than $D_i$); and (b) in case~(O), ~$D_j$ is black and $\gamma(d_i)=1$
(i.e., $D_i$ is lower than $D_{i+1}$)
  \end{numitem1}

See the picture (where the orientation of $D$ is indicated):

\vspace{-0.3cm}
\begin{center}
\includegraphics{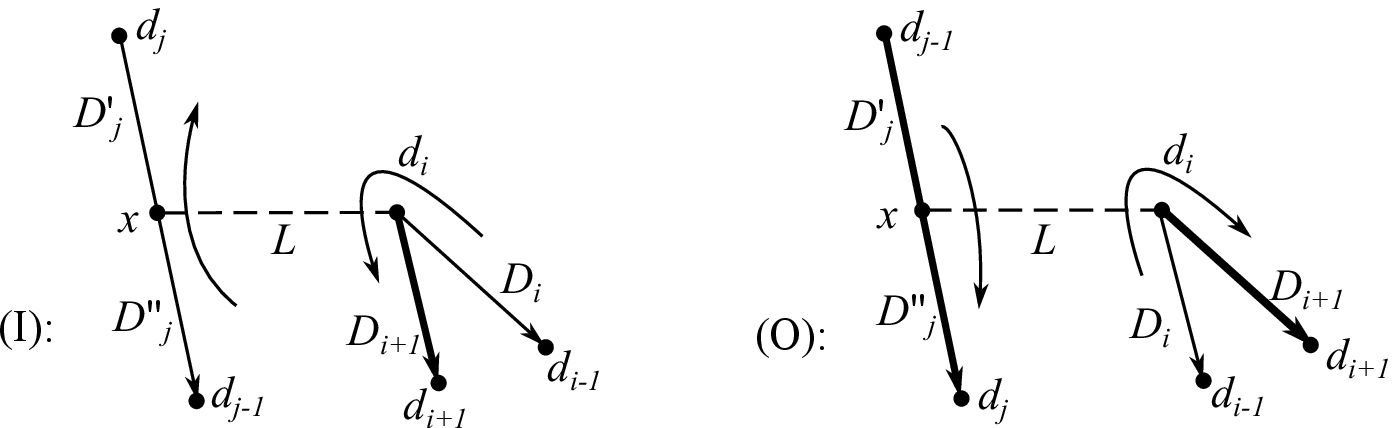}
\end{center}
\vspace{0cm}

The points $x$ and $d_i$ split the cycle (closed curve) $D$ into two parts
$\zeta',\zeta''$, where the former contains $D'_j$ (and $D_i$) and the latter
does $D''_j$ (and $D_{i+1}$).

We first examine case (I). The line $L$ divides the region $O_D$ into two parts
$O'$ and $O''$ lying above and below $L$, respectively. Orienting the curve
$\zeta'$ from $x$ to $d_i$ and adding to it the segment $L$ oriented from $d_i$
to $x$, we obtain closed curve $D'$ surrounding $O'$. Note that $D'$ is
oriented clockwise around $O'$. We combine the paths $D'_j$, $L$ (from $x$ to
$d_i$) and $D_i$ into one directed path $A$ (going from $s_{D'_j}=s_{D_j}=d_j$
to $t_{D_i}=d_{i-1}$). Then $D'$ turns into a correctly colored simple cycle in
which $A$ is regarded as a white snake and the white/black snakes structure of
the rest preserves (cf.~\refeq{casesIO}(a)).

In its turn, the curve $\zeta''$ oriented from $d_{i}$ to $x$ plus the segment
$L$ (oriented from $x$ to $d_i$) form closed curve $D''$ that surrounds $O''$
and is oriented clockwise as well. We combine $L$ and $D_{i+1}$ into one black
snake $B$ (going from $x$ to $d_{i+1}$). Then $D''$ becomes a correctly colored
cycle, and $x$ is a peak in it. (The point $x$ turns into a vertex of $G$.) We
have $\gamma(x)=1$ (since the white $D''_j$ is lower than the black $B$).

The creation of $D',D''$ from $D$ in case (I) is illustrated in the picture:

\vspace{-0.0cm}
\begin{center}
\includegraphics{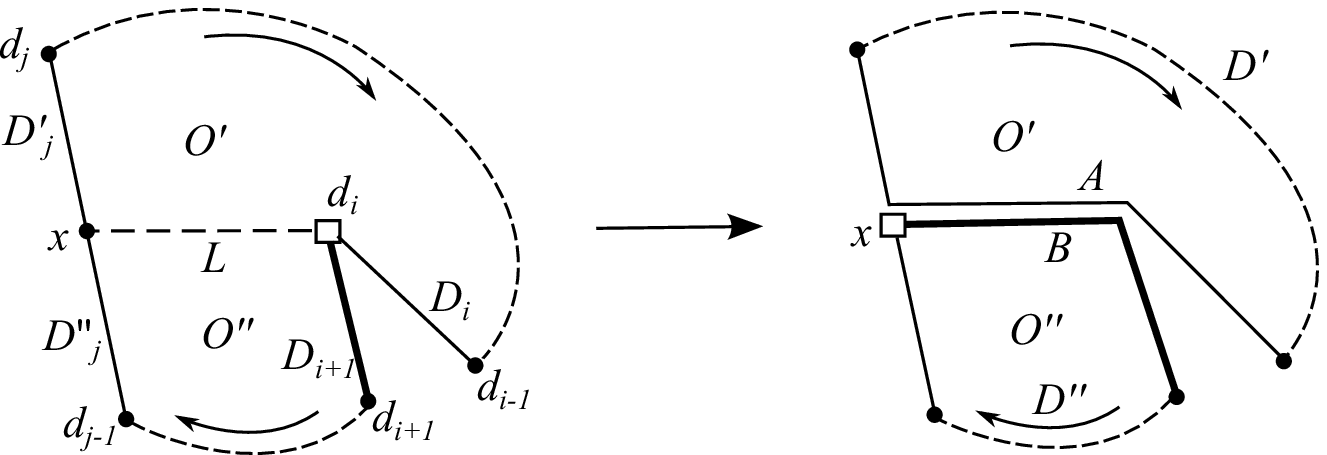}
\end{center}
\vspace{-0.1cm}

We observe that, compared with $D$, the pair $D',D''$ misses the bend $d_i$
(with $\gamma(d_i)=-1$) but acquires the bend $x$ (with $\gamma(x)=1$). Then
   \begin{equation}  \label{eq:DD'D''}
   \eta(D)=\eta(D')+\eta(D''),
   \end{equation}
implying $\eta(D'),\eta(D'')<\eta(D)$. Therefore, we can apply induction. This
gives $\gamma_{D'}=\gamma_{D''}=2$. Now, by reasonings above,
  $$
  \gamma_D=\gamma_{D'}+\gamma_{D''}+\gamma(d_i)-\gamma(x)=2+2-1-1=2,
  $$
as required.

Next we examine case~(O). The curve $\zeta'$ (containing $D'_j$) passes through
the black snake $D_{i+1}$, and the curve $\zeta''$ (containing $D''_j$) through
the white snake $D_i$. Adding to each of $\zeta',\zeta''$ a copy of $L$, we
obtain closed curves $D',D''$, respectively, each inheriting the orientation of
$D$. They become correctly colored simple cycles when we combine the paths
$D'_j,L,D_{i+1}$ into one black snake (from $d_{j-1}$ to $d_{i+1}$) in $D'$,
and combine the paths $L,D_i$ into one white snake (from the new bend $x$ to
$d_i$) in $D''$. Let $O',O''$ be the bounded regions in the plane surrounded by
$D',D''$, respectively. It is not difficult topological exercise to see that
two cases are possible:
  \smallskip

  (O1) ~$O'$ includes $O''$ (and $O_D$);
  \smallskip

  (O2) ~$O''$ includes $O'$ (and $O_D$).

These cases are illustrated in the picture:

\vspace{-0.0cm}
\begin{center}
\includegraphics{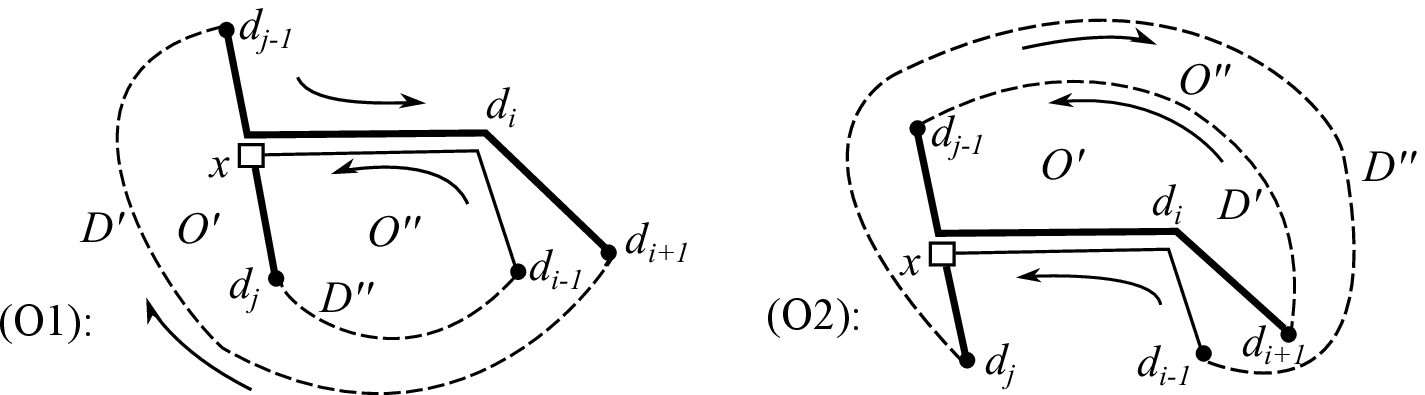}
\end{center}
\vspace{0cm}

Then in case~(O1), ~$D'$ is clockwise and $D''$ is counterclockwise, whereas in
case~(O2) the behavior is converse. Also $\gamma(d_i)=1$ and $\gamma(x)=-1$.
Similar to case~(I), relation~\refeq{DD'D''} is true and we can apply
induction. Then in case~(O1), we have $\gamma_{D'}=2$ and $\gamma_{D''}=-2$,
whence
   $$
  \gamma_D=\gamma_{D'}+\gamma_{D''}+\gamma(d_i)-\gamma(x)=2-2+1-(-1)=2.
  $$
And in case~(O2), we have $\gamma_{D'}=-2$ and $\gamma_{D''}=2$, whence
   $$
  \gamma_D=\gamma_{D'}+\gamma_{D''}+\gamma(d_i)-\gamma(x)=-2+2+1-(-1)=2.
  $$

Thus, in all cases we obtain $\gamma_D=2$, yielding the lemma.
 \end{proof}

This completes the proof of Proposition~\ref{pr:seg-seg}. \hfill\qed

 \subsection{Proof of Proposition~\ref{pr:seg-link}.} \label{ssec:prop3}

Consider a link $L$. By Lemma~\ref{lm:varphi=1}, for any snake $P$,
~$\varphi_{L,P}\ne 1$ is possible only if $L$ and $P$ have a common endvertex
$v$. Note that $v\notin R\cup C$. In particular, it suffices to examine only
bounded and semi-bounded links.

First assume that $s_L\notin R$. Then there are exactly two snakes containing
$s_L$, namely, a white snake $A$ and a black snake $B$ such that $s_L=t_A=t_B$.
If $L$ is white, then $A$ and $L$ belong to the same path in $\phi$; therefore,
$A\prec L\prec B$. Under the exchange operation $A$ becomes black, $B$ becomes
white, and $L$ continues to be white. Then $B,L$ belong to the same path in
$\psi$; this implies $B\precast L\precast A$. So both pairs $(A,L)$ and $(L,B)$
are permuting. Lemma~\ref{lm:1atP=asQ} gives $\varphi_{A,L}=q$ and
$\varphi_{L,B}=\bar q$, whence $\varphi_{A,L}\varphi_{L,B}=1$.

Now let $L$ be black. Then $A\prec B\prec L$ and $B\precast A\precast L$. So
both pairs $\{A,L\}$ and $\{B,L\}$ are invariant, whence
$\varphi_{A,L}=\varphi_{B,L}=1$.

The end $t_L$ is examined in a similar way. Assuming $t_L\notin C$, there are
exactly two snakes, a white snake $A'$ and a black snake $B'$, that contain
$t_L$, namely: $t_L=s_{A'}=s_{B'}$. If $L$ is white, then $L\prec A'\prec B'$
and $L\precast B'\precast A'$. Therefore, $\{L,A\}$ and $\{L,B'\}$ are
invariant, yielding $\varphi_{L,A'}=\varphi_{L,B'}=1$. And if $L$ is black,
then $A'\prec L\prec B'$ and $B'\precast L\precast A'$. So both $(A',L)$ and
$(L,B')$ are permuting, and we obtain from Lemma~\ref{lm:1atP=asQ} that
$\varphi_{A',L}=\bar q$ and $\varphi_{L,B'}=q$, yielding
$\varphi_{A',L}\varphi_{L,B'}=1$.

These reasonings prove the proposition. \hfill\qed

 \subsection{Degenerate case.} \label{ssec:degenerate}

We have proved relation~\refeq{Pi=q} in a non-degenerate case, i.e., subject
to~\refeq{nondegenerate}, and now our goal is to prove~\refeq{Pi=q} when the
set
  $$
  \Zscr:=\{z_1,\ldots,z_{k-1}\}\cup \{c_j\colon j\in J\cup J'\}
  $$
contains distinct elements $u,v$ with $\alpha(u)=\alpha(v)$. We say that such
$u,v$ form a \emph{defect pair}. A special defect pair is formed by twins
$z_i,z_j$ (bends satisfying $i\ne j$, ~$\alpha(z_i)=\alpha(z_j)$ and
$\beta(z_i)=\beta(z_j)$). Another special defect pair is of the form
$\{s_P,t_P\}$ when $P$ is a \emph{vertical} snake or link, i.e.,
$\alpha(s_P)=\alpha(t_P)$.

We will show~\refeq{Pi=q} by induction on the number of defect pairs.

Let $a$ be the \emph{minimum} number such that the set $X:=\{u\in \Zscr\;\colon
\alpha(u)=a\}$ contains a defect pair. We denote the elements of $X$ as
$v_0,v_1,\ldots,v_r$, where for each $i$, ~$v_{i-1}$ is \emph{higher} than
$v_i$, which means that either $\beta(v_{i-1})>\beta(v_i)$, or $v_{i-1},v_i$
are twins and $v_{i-1}$ is a pit (and $v_{i}$ is a peak) in the exchange path
$Z$. The highest element $v_0$ is also denoted by $u$.

In order to conduct induction, we deform the graph $G$ within a sufficiently
narrow vertical strip $S=[a-\eps,a+\eps]\times \Rset$ (where $0<\eps<
\min\{|\alpha(z)-a|\colon z\in \Zscr-X\}$) to get rid of the defect pairs
involving $u$ in such a way that the configuration of snakes/links in the
arising graph $\tilde G$ remains ``equivalent'' to the initial one. More
precisely, we shift the bend $u$ at a small distance ($<\eps$) to the left,
keeping the remaining elements of $\Zscr$; then the bend $u'$ arising in place
of $u$ satisfies $\alpha(u')<\alpha(u)$ and $\beta(u')=\beta(u)$. The
snakes/links with an endvertex at $u$ are transformed accordingly; see the
picture for an example.

\vspace{-0cm}
\begin{center}
\includegraphics{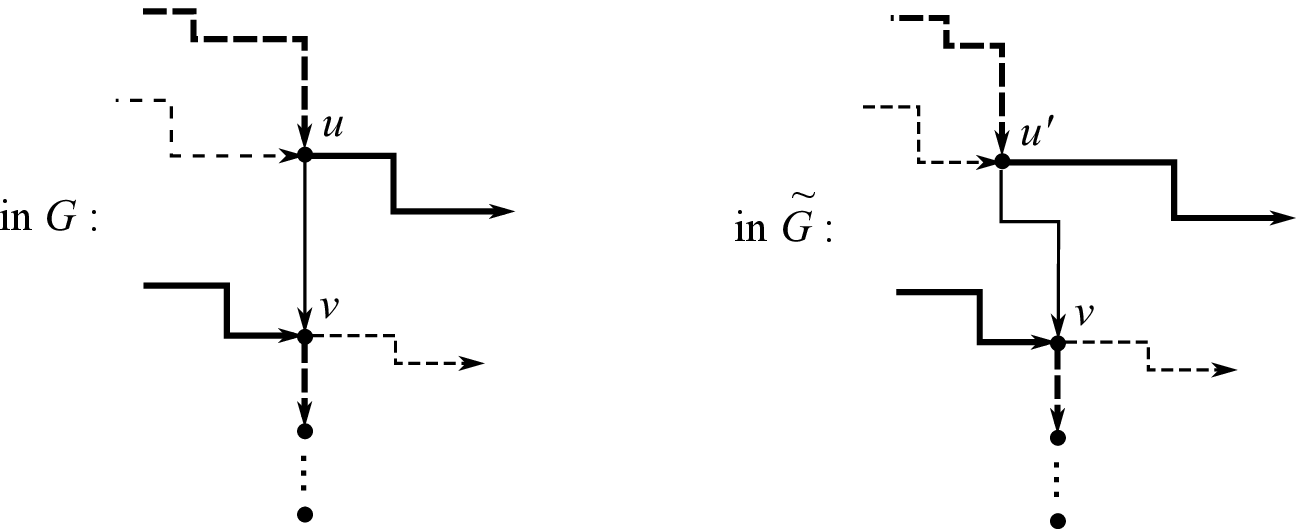}
\end{center}
\vspace{-0.2cm}

Let $\varPi$ and $\tilde\varPi$ denote the L.H.S. value in~\refeq{Pi=q} for the
initial and deformed configurations, respectively. Under the deformation, the
number of defect pairs becomes smaller, so we may assume by induction that
$\tilde\varPi=q$. Thus, we have to prove that
   \begin{equation} \label{eq:varPi}
 \varPi=\tilde\varPi.
  \end{equation}

We need some notation and conventions. For $v\in X$, the set of (initial)
snakes and links with an endvertex at $v$ is denoted by $\Pscr_v$. For
$U\subseteq X$, ~$\Pscr_U$ denotes $\cup(\Pscr_v\;\colon v\in U)$.
Corresponding objects for the deformed graph $\tilde G$ are usually denoted
with tildes as well; e.g.: for a path $P$ in $G$, its image in $\tilde G$ is
denoted by $\tilde P$; the image of $\Pscr_v$ is denoted by $\tilde \Pscr_{v}$
(or $\tilde \Pscr_{\tilde v}$), and so on. The set of standard paths in
$\Pscr_U$ (resp. $\tilde\Pscr_U$) is denoted by $\Pscr^{\rm st}_U$ (resp.
$\tilde\Pscr^{\rm st}_U$). Define
  \begin{equation} \label{eq:u_X-u}
  \varPi_{u,X-u}:=\prod(\varphi_{P,Q}\colon P\in\Pscr_u,\;
Q\in\Pscr_{X-u}).
  \end{equation}
A similar product for $\tilde G$ (i.e., with $\tilde\Pscr_u$ instead of
$\Pscr_u$) is denoted by $\tilde\varPi_{u,X-u}$ .

Note that~\refeq{varPi} is equivalent to
  \begin{equation} \label{eq:varPiX}
  \varPi_{u,X-u}=\tilde\varPi_{u,X-u}.
    \end{equation}
This follows from the fact that for any paths $P,Q\in\Sscr\cup\Lscr$ different
from those involved in~\refeq{u_X-u}, the values $\varphi_{P,Q}$ and
$\varphi_{\tilde P,\tilde Q}$ are equal. (The only nontrivial case arises when
$P,Q\in\Pscr_u$ and $Q$ is vertical (so $\tilde Q$ becomes standard). Then
$t_Q=v_1$. Hence $Q\in \Pscr_{X-u}$, the pair $P,Q$ is involved in
$\varPi_{u,X-u}$, and the pair $\tilde P,\tilde Q$ in $\tilde\varPi_{u,X-u}$.)

To simplify our description technically, one trick will be of use. Suppose that
for each standard path $P\in\Pscr^{\rm st}_X$, we choose a point (not
necessarily a vertex) $v_P\in\Inter(P)$ in such a way that
$\alpha(s_P)<\alpha(v_P)<\alpha(t_P)$, and the coordinates $\alpha(v_P)$ for
all such paths $P$ are different. Then $v_P$ splits $P$ into two subpaths
$P',P''$, where we denote by $P'$ the subpath connecting $s_P$ and $v_P$ when
$\alpha(s_P)=a$, and connecting $v_P$ and $t_P$ when $\alpha(t_P)=a$, while
$P''$ is the rest. This provides the following property: for any
$P,Q\in\Pscr^{\rm st}_X$, ~$\varphi_{P',Q''}=\varphi_{Q',P''}=1$ (in view of
Lemma~\ref{lm:varphi=1}). Hence
$\varphi_{P,Q}=\varphi_{P',Q'}\varphi_{P'',Q''}$. Also $P''=\tilde P''$. It
follows that~\refeq{varPiX} would be equivalent to the equality
  $$
  \prod(\varphi_{P',Q'}\colon P\in\Pscr_u,\;Q\in\Pscr_{X-\{u\}})
  =\prod(\varphi_{\tilde P',\tilde Q'}\colon
  P\in\Pscr_u,\;Q\in\Pscr_{X-\{u\}}).
  $$

In light of these observations, it suffices to prove~\refeq{varPiX} in the
special case when
  \begin{numitem1} \label{eq:assumption}
any $P\in\Pscr_u$ and $Q\in\Pscr_{X-u}$ satisfy
$\{\alpha(s_P),\alpha(t_P)\}\cap \{\alpha(s_Q),\alpha(t_Q)\}=\{a\}$.
  \end{numitem1}

For $i=0,\ldots,r$, we denote by $A_i,B_i,K_i,L_i$, respectively, the white
snake, black snake, white link, and black link that have an endvertex at $v_i$.
Note that if $v_{i-1},v_i$ are twins, then the fact that $v_{i-1}$ is a pit
implies that $A_{i-1},B_{i-1}$ are the snakes entering $v_{i-1}$, and $A_i,B_i$
are the snakes leaving $v_i$; for convenience, we formally define $K_{i-1},K_i,
L_{i-1},L_i$ to be the same trivial path consisting of the single vertex $v_i$.
Note that if $v_r\in C$, then some paths among $A_k,B_k,K_k,L_k$ vanish (e.g.,
both snakes and one link).

When vertices $v_i$ and $v_{i+1}$ are connected by a (vertical) path in
$\Sscr\cup \Lscr$, we denote such a path by $P_i$ and say that the vertex $v_i$
is \emph{open}; otherwise $v_i$ is said to be closed. Note that $v_i,v_{i+1}$
can be connected by either one snake, or one link, or two links (namely,
$K_i,L_i$); in the latter case, $P_i$ is chosen arbitrarily among them. In
particular, if $v_i,v_{i+1}$ are twins, then $v_i$ is open and the role of
$P_i$ is played by any of the trivial links $K_i,L_i$. Obviously, in a sequence
of vertical paths $P_i,P_{i+1},\ldots,P_{j}$, the snakes and links alternate.
One can see that if $P_i$ is a white snake, i.e., $P_i=A_i=A_{i+1}=:A$, then
both black snakes $B_i,B_{i+1}$ are standard, and we have $v_i=s_{B_i}$ and
$v_{i+1}=t_{B_{i+1}}$. See the left fragment of the picture:

\vspace{-0.3cm}
\begin{center}
\includegraphics{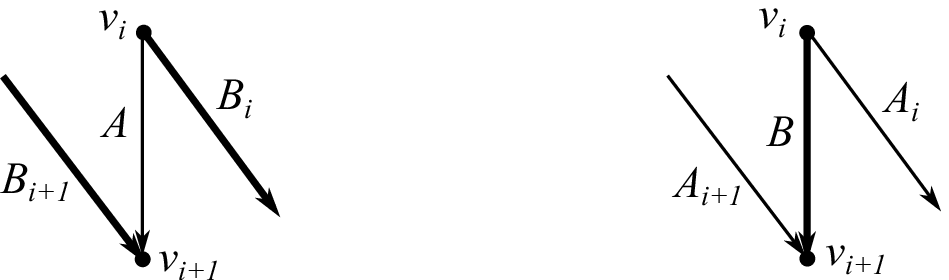}
\end{center}
\vspace{0cm}

Symmetrically, if $P_i$ is a black snake: $B_i=B_{i+1}=:B$, then the white
snakes $A_i,A_{i+1}$ are standard, $v_i=s_{A_i}$ and $v_{i+1}=t_{A_{i+1}}$; see
the right fragment of the above picture.

In its turn, if $P_i$ is a nontrivial white link, i.e., $P_i=K_i=K_{i+1}$, then
two cases are possible: either the black links $L_i,L_{i+1}$ are standard,
$v_i=s_{L_i}$ and $v_{i+1}=t_{L_{i+1}}$, or $L_i=L_{i+1}=P_i$. And if $P_i$ is
a black link, the behavior is symmetric. See the picture:

\vspace{-0.3cm}
\begin{center}
\includegraphics{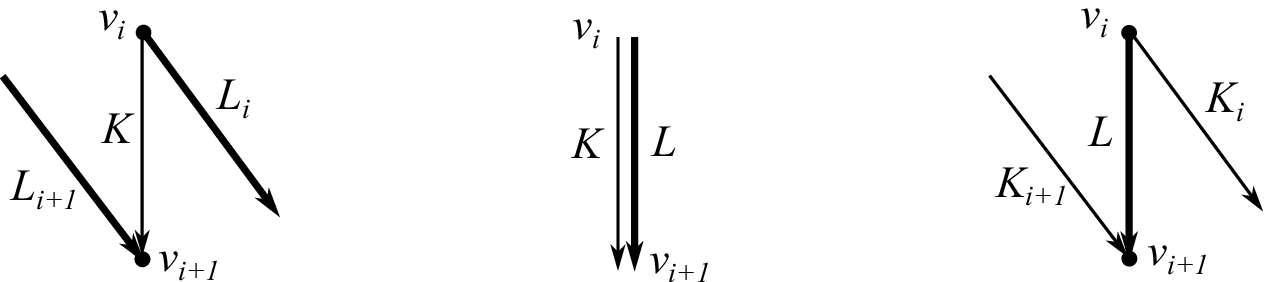}
\end{center}
\vspace{0cm}

Now we are ready to start proving equality~\refeq{varPiX}. Note that the
deformation of $G$ preserves both orders $\prec$ and $\precast$.

We say that paths $P,P'\in\Pstan_X$ are \emph{separated} (from each other) if
they are not contained in the same path of any of the flows
$\phi,\phi',\psi,\psi'$. The following observation will be of use:
  \begin{numitem1}  \label{eq:monochPQ}
if $P,P'\in\Pstan_X$ have the same color, are separated, and $P'$ is lower than
$P$, then $P'\prec P$; and similarly w.r.t. the order $\precast$ (concerning
$\psi,\psi'$).
  \end{numitem1}
Indeed, suppose that $P,P'$ are white, and let $Q$ and $Q'$ be the components
of the flow $\phi$ containing $P$ and $P'$, respectively. Since $P,P'$ are
separated, the paths $Q,Q'$ are different. Moreover, the fact that $P'$ is
lower than $P$ implies that $Q'$ is lower than $Q$ (since $Q,Q'$ are disjoint).
Then $Q'$ precedes $Q$ in $\phi$, yielding $P'\prec P$, as required. When
$P,P'$ concern one of $\phi',\psi,\psi'$, the argument is similar.
 \smallskip

In what follows we will use the abbreviated notation $A,B,K,L$ for the paths
$A_0,B_0,K_0,L_0$ (respectively) having an endvertex at $u=v_0$. Also for
$R\in\Pscr_{X-u}$, we denote the product $\varphi_{A,R}\varphi_{B,R}
\varphi_{K,R}\varphi_{L,R}$ by $\varPi(R)$, and denote by $\tilde \varPi(R)$ a
similar product for the paths $\tilde A,\tilde B,\tilde K, \tilde L, \tilde R$
(concerning the deformed graph $\tilde G$). One can see that $\varPi_{u,X-u}$
(resp. $\tilde \varPi_{u,X-u}$) is equal to the product of the values
$\varPi(R)$ (resp. $\tilde\varPi(R)$) over $R\in\Pscr_{X-u}$.

To show~\refeq{varPiX}, we examine several cases. First we consider
 \smallskip

\noindent\underline{\emph{Case (R1)}:} ~the vertex $u$ is closed; in other
words, all paths $A,B,K,L$ are standard.

  \begin{prop}  \label{pr:caseR1}
In case~(R1), ~$\varPi(R)=\tilde\varPi(R)=1$ holds for any $R\in\Pscr_{X-u}$.
As a consequence, \refeq{varPiX} is valid.
  \end{prop}
  \begin{proof}
~Let $R\in \Pscr_{v_p}$ for $p\ge 1$. Observe that~\refeq{assumption} together
with the fact that the vertex $u$ moves under the deformation of $G$ implies
that $\{\alpha(s_{\tilde P}),\alpha(t_{\tilde P})\} \cap \{\alpha(s_{\tilde
R}),\alpha(t_{\tilde R})\}=\emptyset$ holds for any $P\in \Pscr_u$. This gives
$\tilde\varPi(R)=1$, by Lemma~\ref{lm:varphi=1}.

Next we show the equality $\varPi(R)=1$. One may assume that $R$ is standard
(otherwise the equality is trivial). Since $u$ is closed, $A,B,K,L$ are
separated from $R$.

Note that $A,B,K,L,R$ are as follows: either (a) $t_A=t_B=s_K=s_L$ or (b)
$s_A=s_B=t_K=t_L$, and either (c) $\alpha(s_R)=a$ or (d) $\alpha(t_R)=a$. Let
us examine the possible cases when the combination of~(a) and~(d) takes place.
  \smallskip

1) Let $R$ be a white link, i.e., $R=K_p$. Since $R$ is white and lower than
$A,B,K,L$, we have $R\prec A,B,K,L$ (cf.~\refeq{monochPQ}). The exchange
operation preserves the color of $R$. Then $R\precast A,B,K,L$. Therefore, all
pairs $\{P,R\}$ with $P\in\Pscr_u$ are invariant, and $\varPi(R)=1$ is trivial.
 \smallskip

2) Let $R=L_p$. Since $R$ is black, we have $A,K\prec R\prec B,L$. The exchange
operation changes the colors of $A,B$ and preserves the ones of $K,L,R$. Hence
$B,K\precast R\precast A,L$, giving the permuting pairs $(A,R)$ and $(R,B)$.
Lemma~\ref{lm:atP=atQ} applied to these pairs implies $\varphi_{A,R}=\bar q$
and $\varphi_{R,B}=q$. Then $\varPi(R)=\varphi_{A,R}\varphi_{R,B}=\bar q q=1$.
\smallskip

3) Let $R=A_p$. Then $R\prec A,B,K,L$ and $B,K\precast R\precast A,L$ (since
the exchange operation changes the colors of $A,B,R$). This gives the permuting
pairs $(R,B)$ and $(R,K)$. Then $\varphi_{R,B}=q$, by Lemma~\ref{lm:atP=atQ},
and $\varphi_{R,K}=\bar q$ by Lemma~\ref{lm:2atP=asQ}, and we have
$\varPi(R)=\varphi_{R,B}\varphi_{R,K}=1$.
 \smallskip

4) Let $R=B_p$. (In fact, this case is symmetric to the previous one, as it is
obtained by swapping $(\phi,\phi')$ and $(\psi,\psi')$. Yet we prefer to give
details.) We have $A,K\prec R\prec B,L$ and $R\precast A,B,K,L$, giving the
permuting pairs $(A,R)$ and $(K,R)$. Then $\varphi_{A,R}=\bar q$, by
Lemma~\ref{lm:atP=atQ}, and $\varphi_{K,R}=q$, by Lemma~\ref{lm:2atP=asQ},
whence $\varPi(R)=1$.
 \smallskip

The other combinations, namely,~(a) and~(c), ~(b) and~(c), ~(b) and~(d), are
examined in a similar way (where we appeal to appropriate lemmas from
Sect.~\SEC{two_paths}), and we leave this to the reader as an exercise.
\end{proof}

Next we consider
\smallskip

\noindent\underline{\emph{Case (R2)}:} ~$u$ is open; in other words, at least
one path among $A,B,K,L$ is vertical (going from $u=v_0$ to $v_1$).
  \smallskip

It falls into several subcases examined in propositions below.

 \begin{prop} \label{pr:caseR2_sep}
In case~(R2), let $R\in\Pstan_{X-u}$ be separated from $A,B,K,L$. Then
$\varPi(R)=\tilde\varPi(R)$.
  \end{prop}
  \begin{proof}
~We first assume that $u$ and $v_1$ are connected by exactly one path $P_0$
(which is one of $A,B,K,L$) and give a reduction to the previous proposition,
as follows.

Suppose that we replace $P_0$ by a standard path $P'$ of the same color and
type (snake or link) such that $s_{P'}=u$ (and $\alpha(t_{P'})>a$). Then the
set $\Pscr'_u:=(\{A,B,K,L\}-\{P_0\})\cup\{P'\}$ becomes as in case~(R1), and by
Proposition~\ref{pr:caseR1}, the corresponding product $\Pi'(R)$ of values
$\varphi_{R,P}$ over $P\in\Pscr'_u$ is equal to 1. (This relies on the fact
that $R$ is separated from $A,B,K,L$.)

Now compare the effects from $P'$ and $\tilde P_0$. These paths have the same
color and type, and both are separated from, and higher than $R$. Also
$\alpha(s_{P'})=\alpha(t_{\tilde P_0})=a$ (since $s_{P'}=u$ and $t_{\tilde
P_0}=v_1$). Then using appropriate lemmas from Sect.~\SEC{two_paths}, one can
conclude that $\{\varphi_{R,P'},\varphi_{R,\tilde P_0}\}=\{q,\bar q\}$.
Therefore,
   $$
   \tilde\varPi(R)=\varphi_{R,\tilde P_0}=\varPi'(R)\varphi^{-1}_{R,P'} =\varPi(R).
   $$

Now let $u$ and $v_1$ be connected by two paths, namely, by $K,L$. We can again
appeal to Proposition~\ref{pr:caseR1}. Consider $\Pscr''_u:=\{A,B,K'',L''\}$,
where $K'',L''$ are standard links (white and black, respectively) with
$s_{K''}=s_{L''}=u$. Then $\varPi''(R):= \varPi( \varphi_{R,P}\colon
P\in\Pscr''_u)=1$ and $\{\varphi_{R,K''}, \varphi_{R,\tilde
K}\}=\{\varphi_{R,L''}, \varphi_{R,\tilde L}\}=\{q,\bar q\}$, and we obtain
  $$
  \tilde\varPi(R)=\varphi_{R,\tilde K}\varphi_{R,\tilde L}=
  \varPi''(R)\varphi^{-1}_{R,K''}\varphi^{-1}_{R,L''} =\varphi_{R,A}\varphi_{R,B}
   =\varPi(R),
   $$
as required.
 \end{proof}

 \begin{prop} \label{pr:caseR2_nonsep}
In case~(R2), let $R$ be a standard path in $\Pscr_{v_p}$ with $p\ge 1$. Let
$R$ be not separated from at least one of $A,B,K,L$. Then
$\varPi(R)=\tilde\varPi(R)$.
  \end{prop}
  \begin{proof}
We first assume that $P_0$ is the unique vertical path connecting $u$ and $v_1$
(in particular, $u$ and $v_1$ are not twins). Then $R$ is not separated from
$P_0$.

Suppose that $P_0$ and $R$ are contained in the same path of the flow $\phi$;
equivalently, both $P_0,R$ are white and $P_0\prec R$. Then neither $\psi$ nor
$\psi'$ has a path containing both $P_0,R$ (this is easy to conclude from the
fact that one of $R$ and $P_{p-1}$ is a snake and the other is a link).
Consider four possible cases for $P_0,R$.

(a) Let both $P_0,R$ be links, i.e., $P_0=K$ and $R=K_p$. Then $A,K\prec
K_p\prec B,L$ and $K_p\precast B,K,A,L$ (since $K\precast K_p$ is impossible by
the above observation). This gives the permuting pairs $(A,K_p)$ and $(\tilde
K,K_p)$, yielding $\varphi_{A,K_p}=\varphi_{\tilde K,K_p}$.

(b) Let $P_0=K$ and $R=A_p$. Then $A,K\prec A_p\prec B,L$ and $B,K\precast
A_p\precast A,L$. This gives the permuting pairs $(A,A_p)$ and $(A_p,B)$,
yielding $\varphi_{A,A_p}\varphi_{\tilde A_p,B}=1=\varphi_{\tilde K,A_p}$.

(c) Let $P_0=A$ and $R=K_p$. Then $K,A\prec K_p\prec L,B$ and $K_p\precast
K,B,L,A$. This gives the permuting pairs $(K,K_p)$ and $(\tilde A,K_p)$,
yielding $\varphi_{K,K_p}=\varphi_{\tilde A,K_p}$.

(d) Let $P_0=A$ and $R=A_p$. Then $K,A\prec A_p\prec L,B$ and $K,B\precast
A_p\precast L,A$. This gives the permuting pairs $(\tilde A,A_p)$ and
$(A_p,B)$, yielding $\varphi_{\tilde A,A_p}=\varphi_{A_p,B}$.

In all cases, we obtain $\varPi(R)=\tilde\varPi(R)$.

When $P_0,R$ are contained in the same path in $\phi'$ (i.e., $P_0,R$ are black
and $P_0\prec R$), we argue in a similar way. The cases with $P_0,R$ contained
in the same path of $\psi$ or $\psi'$ are symmetric.

A similar analysis is applicable (yielding $\varPi(R)=\tilde\varPi(R)$) when
$u$ and $v_1$ are connected by two vertical paths (namely, $K,L$) and exactly
one relation among $K\prec R$, $L\prec R$, $K\precast R$ and $L\precast R$
takes place (equivalently: either $K,R$ or $L,R$ are separated, not both).

Finally, let $u$ and $v_1$ be connected by both $K,L$, and assume that $K,R$
are not separated, and similarly for $L,R$. An important special case is when
$p=1$ and $u,v_1$ are twins. Note that from the assumption it easily follows
that $R$ is a snake. If $R$ is the white snake $A_p$, then we have $A,K\prec
A_p\prec B,L$ and $B,K,A,L\precast A_p$. This gives the permuting pairs
$(A_p,B)$ and $(A_p,\tilde L)$, yielding $\varphi_{A_p,B}=\varphi_{A_p,\tilde
L}$ (since $\alpha(t_B)=\alpha(t_{\tilde L}$)). The case with $R=B_p$ is
symmetric. In both cases, $\varPi(R)=\tilde\varPi(R)$.
  \end{proof}

 \begin{prop} \label{pr:caseR2_P0}
Let $R=P_0$ be the unique vertical path connecting $u$ and $v_1$. Then
$\varPi(R)=\tilde\varPi(R)=1$.
  \end{prop}
  \begin{proof}
~The equality $\varPi(R)=1$ is trivial. To see $\tilde\varPi(R)=1$, consider
possible cases for $R$. If $R=K$, then $\tilde A\prec \tilde K\prec \tilde
B,\tilde L$ and $\tilde B\precast \tilde K\precast \tilde A,\tilde L$, giving
the permuting pairs $(\tilde A,\tilde K)$ and $(\tilde K,\tilde B)$ (note that
$t_{\tilde A}=t_{\tilde B}=s_{\tilde K}=\tilde u$). If $R=L$, then $\tilde
A,\tilde K,\tilde B\prec\tilde L$ and $\tilde B,\tilde K,\tilde A\precast
\tilde L$; so all pairs involving $\tilde L$ are invariant. If $R=A$, then
$\tilde K\prec\tilde A\prec \tilde L,\tilde B$ and $\tilde K,\tilde B,\tilde
L\precast \tilde A$, giving the permuting pairs $(\tilde A,\tilde L)$ and
$(\tilde A,\tilde B)$ (note that $s_{\tilde A}=s_{\tilde B}=t_{\tilde L}=\tilde
u$). And the case $R=B$ is symmetric to the previous one.

In all cases, using appropriate lemmas from Sect.~\SEC{two_paths} (and relying
on the fact that all paths $\tilde A,\tilde B,\tilde K,\tilde L$ are standard),
one can conclude that $\tilde\varPi(R)=1$.
  \end{proof}

 \begin{prop} \label{pr:caseR2_KL}
Let both $K,L$ be vertical. Then $\varPi(K)\varPi(L)=
\tilde\varPi(K)\tilde\varPi(L)=1$.
  \end{prop}
  \begin{proof}
The equality $\varPi(K)\varPi(L)=1$ is trivial. To see
$\tilde\varPi(K)\tilde\varPi(L)=1$, observe that $\tilde A\prec\tilde K\prec
\tilde B\prec \tilde L$ and $\tilde B\precast \tilde K\precast \tilde A\precast
\tilde L$. This gives the permuting pairs $(\tilde A,\tilde K)$ and $(\tilde
K,\tilde B)$. By Lemma~\ref{lm:1atP=asQ}, $\varphi_{\tilde A,\tilde K}=q$ and
$\varphi_{\tilde K\tilde B}=\bar q$, and the result follows.
  \end{proof}

Taken together, Propositions~\ref{pr:caseR2_sep}--\ref{pr:caseR2_KL} embrace
all possibilities in case~(R2). Adding to them Proposition~\ref{pr:caseR1}
concerning case~(R1), we obtain the desired relation~\refeq{varPiX} in a
degenerate case.

This completes the proof of Theorem~\ref{tm:single_exch} in case~(C), namely,
relation~\refeq{caseC}.  \hfill\qed\qed

 \subsection{Other cases.} \label{ssec:othercases}
Let $(I|J),(I'|J'),\phi,\phi',\psi,\psi'$ and $\pi=\{f,g\}$ be as in the
hypotheses of Theorem~\ref{tm:single_exch}. We have proved this theorem in
case~(C), i.e., when $\pi$ is a $C$-couple with $f<g$ and $f\in J$ (see the
beginning of Sect.~\SEC{exchange}). In other words, the exchange path
$Z=P(\pi)$, used to transform the initial double flow $(\phi,\phi')$ into the
new double flow $(\psi,\psi')$, connects the sinks $c_f$ and $c_g$ covered by
the ``white flow'' $\phi$ and the ``black flow'' $\phi'$, respectively.

The other possible cases in the theorem are as follows:
  \smallskip

(C1) ~$\pi$ is a $C$-couple with $f<g$ and $f\in J'$;
\smallskip

(C2) ~$\pi$ is an $R$-couple with $f<g$ and $f\in I$;
\smallskip

(C3) ~$\pi$ is an $R$-couple with $f<g$ and $f\in I'$;
\smallskip

(C4) ~$\pi$ is an $RC$-couple with $f\in I$ and $g\in J$;
\smallskip

(C5) ~$\pi$ is an $RC$-couple with $f\in I'$ and $g\in J'$.
\smallskip

Case~(C1) is symmetric to~(C). This means that if double flows $(\phi,\phi')$
and $(\psi,\psi')$ are obtained from each other by applying the exchange
operation using $\pi$ (which, in particular, changes the ``colors'' of both $f$
and $g$), and if one double flow is subject to~(C), then the other is subject
to~(C1). Rewriting $w(\phi)w(\phi')=qw(\psi)w(\psi')$ as
$w(\psi)w(\psi')=q^{-1}w(\phi)w(\phi')$, we just obtain the required equality
in case~(C1) (where $(\psi,\psi')$ and $(\phi,\phi')$ play the roles of the
initial and updated double flows, respectively).

For a similar reasons, case~(C3) is symmetric to~(C2), and~(C5) is symmetric
to~(C4). So it suffices to establish the desired equalities merely in
cases~(C2) and~(C4).

To do this, we appeal to reasonings similar to those in
Sects.~\SSEC{prop1}--\SSEC{degenerate}. More precisely, one can check that the
descriptions in Sects.~\SSEC{prop1} and \SSEC{prop3} (concerning link-link and
snake-link pairs in $\Nscr$) remain applicable and
Propositions~\ref{pr:link-link} and~\ref{pr:seg-link} are directly extended to
cases~(C2) and~(C4). The method of getting rid of degeneracies developed in
Sect.~\SSEC{degenerate} does work, without any troubles, for~(C2) and~(C4) as
well.

As to the method in Sect~\SSEC{prop2} (concerning snake-snake pairs in
case~(C)), it should be modified as follows. We use terminology and notation
from Sects.~\SSEC{seglink} and~\SSEC{prop2} and appeal to
Lemma~\ref{lm:gammaD}.

When dealing with case~(C2), we represent the exchange path $Z=P(\pi)$ as a
concatenation $Z_1\circ \bar Z_2\circ Z_3\circ\cdots \circ\bar Z_k$, where each
$Z_i$ with $i$ odd (even) is a snake contained in the black flow $\phi'$ (resp.
the white flow $\phi$). Then $Z_1$ begins at the source $r_g$ and $Z_k$ begins
at the source $r_f$. An example with $k=6$ is illustrated in the left fragment
of the picture:

\vspace{0cm}
\begin{center}
\includegraphics{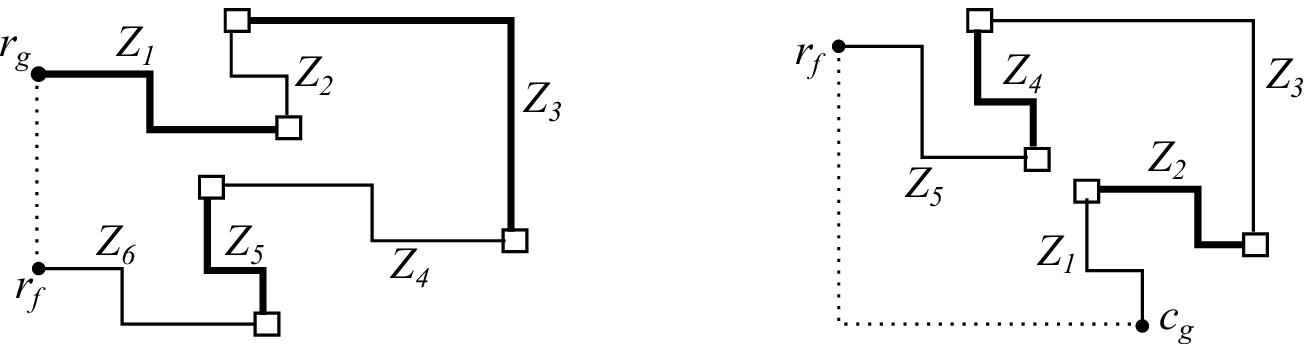}
\end{center}
\vspace{0cm}

The common vertex (bend) of $Z_i$ and $Z_{i+1}$ is denoted by $z_i$. As before,
we associate with a bend $z$ the number $\gamma(z)$ (equal to 1 if, in the pair
of snakes sharing $z$, the white snake is lower that the black one, and $-1$
otherwise), and define $\gamma_Z$ as in~\refeq{gammaZ}. We turn $Z$ into simple
cycle $D$ by combining the directed path $Z_k$ (from $r_f$ to $z_{k-1}$) with
the vertical path from $r_g$ to $r_f$, which is formally added to $G$. (In the
above picture, this path is drawn by a dotted line.) Then, compared with $Z$,
the cycle $D$ has an additional bend, namely, $r_g$. Since the extended white
path $\tilde Z_k$ is lower than the black path $Z_1$, we have $\gamma(r_g)=1$,
and therefore $\gamma_D=\gamma_Z+1$.

One can see that the cycle $D$ is oriented clockwise (where, as before, the
orientation is according to that of black snakes). So $\gamma_D=2$, by
Lemma~\ref{lm:gammaD}, implying $\gamma_Z=1$. This is equivalent to the
``snake-snake relation'' $\varphi^{II}=q$, and as a consequence, we obtain the
desired equality
  $$
  w(\phi)w(\phi')=qw(\psi)w(\psi').
  $$

Finally, in case~(C4), we represent the exchange path $Z$ as the corresponding
concatenation $\bar Z_1\circ Z_2\circ \bar Z_3\circ\cdots \circ Z_{k-1}\circ
\bar Z_k$ (with $k$ odd), where the first white snake $Z_1$ ends at the sink
$c_g$ and the last white snake $Z_k$ begins at the source $r_f$. See the right
fragment of the above picture, where $k=5$. We turn $Z$ into simple cycle $D$
by adding a new ``black snake'' $Z_{k+1}$ beginning at $r_f$ and ending at
$c_g$ (it is formed by the vertical path from $r_f$ to $(0,0)$, followed by the
horizontal path from $(0,0)$ to $c_g$; see the above picture). Compared with
$Z$, the cycle $D$ has two additional bends, namely, $r_f$ and $c_g$. Since the
black snake $Z_{k+1}$ is lower than both $Z_1$ and $Z_k$, we have
$\gamma(r_f)=\gamma(c_g)=-1$, whence $\gamma_D=\gamma_Z-2$. Note that the cycle
$D$ is oriented counterclockwise. Therefore, $\gamma_D=-2$, by
Lemma~\ref{lm:gammaD}, implying $\gamma_Z=0$. As a result, we obtain the
desired equality $w(\phi)w(\phi')=w(\psi)w(\psi')$.

This completes the proof of Theorem~\ref{tm:single_exch}.

\end{document}